\documentclass[12pt]{article}
\usepackage[russian]{babel}
\usepackage[cp1251]{inputenc}
\usepackage{latexsym}
\usepackage{amsmath}
\usepackage{amssymb}
\usepackage{stackrel}
\usepackage{color}
\usepackage{graphicx}
\usepackage{hyperref}

\usepackage[left=2cm,right=2cm,
    top=2cm,bottom=2cm,bindingoffset=0cm]{geometry}

\begin{document}
\pagestyle{empty}
\begin{center}
\bf\large Convergence of (generalized) power series solutions \\ of functional equations
\end{center}
\bigskip

\centerline{R.\,Gontsov, I.\,Goryuchkina}
\bigskip

\centerline{\bf\large Extended abstract}
\bigskip

The article is planned to be published in a Russian journal having an English version. Thus before the text in Russian we put an extended abstract in English containing all the main statements.
\bigskip

{\bf Abstract.} Solutions of nonlinear functional equations are generally not expressed as a finite number of combinations and compositions of elementary and known special functions. One of the approaches to study them is, firstly, to find {\it formal} solutions (that is, series whose terms are described and ordered in some way but which don't converge {\it apriori}) and, secondly, to study the convergence or summability of these formal solutions (the existence and uniqueness of actual solutions with the given asymptotic expansion in a certain domain). In this paper we deal only with the {\it convergence} of formal functional series having the form of an infinite sum of power functions with (complex, in general) power exponents and satisfying analytical functional equations of the following three types: a differential, $q$-difference or Mahler equation.

\bigskip
\centerline{\bf Table of content}
\begin{enumerate} \item Introduction \item Theorems on existence and uniqueness of holomorphic solutions of functional equations
\begin{enumerate}
  \item[2.1] An analytical differential equation
  \item[2.2] An analytical $q$-difference equation
  \begin{enumerate}
    \item[2.2.1] An example: the linearization of a local diffeomorphism
  \end{enumerate}
  \item[2.3] An analytical Mahler equation
  \begin{enumerate}
    \item[2.3.1] An example: the reduction of a germ of a holomorphic function to its main part
  \end{enumerate}
\end{enumerate}
\item Generalized power series that formally satisfy an analytical functional equation
\begin{enumerate}
\item[3.1] The ring of generalized power series
\item[3.2] The substitution of a generalized power series into an analytical functional equation. A formal solution
\item[3.3] Lemmas on reducing functional equations to a special form
\item[3.4] The structure of the set of power exponents of a generalized power series solution of an analytical functional equation
\item[3.5] Representing the elements of the algebra ${\mathbb C}[[x^{\Gamma}]]$ in the form of multi-variable Taylor series
\end{enumerate}

\item A condition of convergence for a generalized power series satisfying an ordinary differential equation
\item A condition of convergence for a generalized power series satisfying a $q$-difference equation
\begin{enumerate}
\item[5.1] The proof of Lemma 7
\item[5.2] Weakening the conditions of Theorem 6 in the special cases of the location of the generators of the semigroup $\Gamma$ on the complex plane
\item[5.3] Examples
\end{enumerate}
\item A condition of convergence for a generalized power series satisfying a Mahler equation
\item Remarks on functional equations of a more general form than $q$-difference or Mahler ones
\begin{enumerate}
  \item[7.1] Existence and uniqueness of a local holomorphic solution
  \item[7.2] A formal solution in the form of a generalized power series and its convergence
\end{enumerate}

\end{enumerate}
\bigskip

In the extended abstract we keep the numeration of formulae as in the original manuscript attached below.
\bigskip

{\bf Some notations:}
\bigskip

$\mathcal{G}$ is the ring of all generalized power series of the form $(6)$.
\medskip

$o(x^\nu)\in\cal G$ is a generalized power series starting with a power exponent whose real part is strictly greater than ${\rm Re}\,\nu$.
\medskip

$\delta$ is the differential operator defined as $\displaystyle\delta=x\frac{d}{dx}$.
\medskip

$\sigma$ is the $q$-difference operator locally defined on functions as $\sigma: y(x)\rightarrow y(qx)$, $q\neq0,1$, and naturally extended on $\cal G$.
\medskip

$\mu$ is the Mahler operator locally defined on functions as $\mu: y(x)\rightarrow y(x^\ell)$, $\ell\in\mathbb{Z}_{\geqslant 2}$, and naturally extended on $\cal G$.
\medskip

${\mathbb C}[[x^\Gamma]]\subset\mathcal{G}$ is the ring of generalized power series whose power exponents lie in the finitely generated additive semigroup $\Gamma$.
\bigskip

\centerline{\bf Main statements}
\bigskip

{\bf Lemma 1} (section 3). {\it Let $F(x, y_0, y_1,\dots,y_n)$ be a holomorphic function near $0\in\mathbb{C}^{n+2}$,
$\varphi\in\mathcal{G}$, and $\Delta$ one of the operators $\delta$, $\sigma$, $\mu$. Then $F(x,\varphi,\Delta\varphi,\dots,\Delta^n\varphi)\in\mathcal{G}$.}
\medskip

This lemma explains the correctness of the notion of a {\it formal generalized power series solution} of the equation $F(x,y,\Delta y,\dots,\Delta^ny)=0$.
\bigskip

{\bf Lemma 2} (section 3, lemmas 2+4+5). {\it Let the generalized power series $(6)$ satisfy equation $(3)$ or equation $(4)$ and
$$
F'_{y_j}(x,\varphi,\Delta\varphi,\ldots,\Delta^n\varphi)=A_jx^{\nu}+o(x^\nu)\in{\cal G},
$$
where $\Delta=\delta$ or $\sigma$, respectively, the number $\nu\in\mathbb C$ is the same for all $j=0,1,\ldots,n$, and at least one of the $A_j$'s is nonzero. Then there exists
$N\in\mathbb N$ such that a transformation
$$
y=\sum_{k=1}^Nc_kx^{\lambda_k}+x^{\lambda_N}u
$$
reduces
\begin{itemize}
\item[\rm i)] equation $(3)$ to the equation
$$L(\delta+\lambda_N)u=M(x,u,\delta u,\ldots,\delta^nu);\eqno(25)
$$
\item[\rm ii)] equation $(4)$ to the equation
$$L(q^{\lambda_N}\sigma)u=M(x,u,\sigma u,\ldots,\sigma^nu),\eqno(26)
$$
\end{itemize}
where $L(z)=A_nz^n+\ldots+A_1z+A_0\not\equiv0$ is a polynomial of degree $\leqslant n$,
$$
M(x,u_0,u_1,\ldots,u_n)=\sum_{{\bf p}\in{\mathbb Z}_+^{n+1}}a_{\bf p}(x)u_0^{p_0}u_1^{p_1}\ldots u_n^{p_n}
$$
is a holomorphic function in variables $u_j$ that is represented by a series with coefficients $a_{\bf p}\in\cal G$, which are generalized power series with a nonzero radius of convergence. Moreover, power exponents of the latter belong to a finitely generated additive semigroup whose generators are linearly independent over $\mathbb Z$ and have positive real parts.}
\bigskip

{\bf Lemma 3} (section 3, lemmas 3+4+5). {\it Let the generalized power series $(6)$ satisfy equation $(5)$ and
$$
F'_{y_0}(x,\varphi,\mu\varphi,\ldots,\mu^n\varphi)=A_0x^{\nu_0}+o(x^{\nu_0})\in{\cal G}, \qquad A_0\ne0.
$$
Then there exists $N\in\mathbb N$ such that the transformation
$$
y=\sum_{k=1}^Nc_kx^{\lambda_k}+x^{\lambda_N}u
$$
reduces equation $(5)$ to the equation
$$
u=M(x,u,\mu u,\ldots,\mu^nu),\eqno(29)
$$
where the function $M$ is the same as in Lemma $2$.}
\bigskip

{\bf Theorem 4} (section 3). {\it Let the generalized power series $(6)$ satisfy
\begin{itemize}
\item[--]  equation $(3)$ or equation $(4)$ and the condition of Lemma $2$ or
\item[--]  equation $(5)$ and the condition of Lemma $3$.
\end{itemize}
Then $\varphi\in{\mathbb C}[[x^\Gamma]]\subset\cal G$, where $\Gamma$ is a finitely generated additive semigroup whose generators are linearly independent over $\mathbb Z$ and have positive real parts.}
\bigskip

Linear independence of the generators of $\Gamma$ over $\mathbb Z$ (denote them by $\rho_1,\ldots,\rho_\tau$) allows to define a linear map
$\iota: {\mathbb C}[[x^\Gamma]]\rightarrow{\mathbb C}[[x_1,\ldots,x_\tau]]$ from the ring ${\mathbb C}[[x^\Gamma]]$ to the ring of formal Taylor series in $\tau$ independent variables (more precisely, we consider this map between the corresponding maximal ideals of these rings consisting of series without constant term). This map is an isomorphism.
\bigskip

{\bf Lemma 6} (section 3). {\it For any series of the form $(35)$ and generalized power series $\varphi_0,\varphi_1,\ldots,\varphi_n\in{\mathbb C}[[x^\Gamma]]$ there is an equality
$$
\iota\bigl(F(x,\varphi_0,\varphi_1,\ldots,\varphi_n)\bigr)=\sum_{{\bf p}\in{\mathbb Z}_+^{n+1}}\iota(b_{\bf p})\iota(\varphi_0)^{p_0}\iota(\varphi_1)^{p_1}\ldots\iota(\varphi_n)^{p_n}.\eqno(36)
$$
}
\medskip

{\bf Theorem 5} (section 4). {\it Let the generalized power series $(6)$ satisfies equation $(3)$ and
$$
F'_{y_j}(x,\varphi,\delta\varphi,\ldots,\delta^n\varphi)=A_jx^{\nu}+o(x^\nu)\in{\cal G},
$$
where $\nu\in\mathbb C$ is the same for all $j=0,1,\ldots,n$ and $A_n\ne0$. Then the series $\varphi$ uniformly converges in each sector $S$ of a small radius with the vertex at the origin and of the opening less than $2\pi$, defining there a holomorphic function.}
\bigskip

Let the generalized power series $(6)$ satisfies equation $(4)$ and
$$
F'_{y_j}(x,\varphi,\sigma\varphi,\ldots,\sigma^n\varphi)=A_jx^{\nu}+o(x^\nu)\in{\cal G},
$$
where $\nu\in\mathbb C$ is the same for all $j=0,1,\ldots,n$ and at least one of the $A_j$'s is nonzero. Then by Theorem 4 the series $\varphi$ is an element of the ring ${\mathbb C}[[x^\Gamma]] \subset\cal G$. The next two statements, Theorems 6 and 6bis, provide sufficient conditions for its convergence.
\bigskip

{\bf Theorem 6} (section 5). {\it Let $A_n\ne0$ and $A_0\ne0$, and for each root $z=a$ of the polynomial $(z-1)L(z)$, $L(z)=A_nz^n+\ldots+A_1z+A_0$, the following diophantine condition hold:
$$
|(m_1\rho_1+\ldots+m_\tau\rho_\tau)\ln q-\ln a-2\pi m{\rm i}|>c\,|m_1+\ldots+m_\tau|^{-\gamma} \quad \forall\; m_i\in{\mathbb Z}_+,\; m\in{\mathbb Z} \eqno(46)
$$
$($with the exception of $m_1=\ldots=m_\tau=0)$, where $c$ and $\gamma$ are some positive constants. Then the series $\varphi$ uniformly converges in each sector $S$ of a small radius with the vertex at the origin and of the opening less than $2\pi$, defining there a holomorphic function.}
\medskip

{\bf Theorem 6 bis} (section 5, page 38) {\it The statement of Theorem $6$ is correct in the following particular cases:
\smallskip

{\rm a)} $A_0\ne0$, and all the $q^{\rho_i}$'s are located {\rm strictly inside} the unite circle $\{|z|=1\}$;

{\rm b)} $A_n\ne0$, and all the $q^{\rho_i}$'s are located {\rm strictly outside} the unite circle;

{\rm c)} all the $q^{\rho_i}$'s are located {\rm on} the unite circle, and condition $(46)$ is fulfilled for those roots $z=a$ of the polynomial $(z-1)L(z)$ having the same location;

{\rm d)} $A_0\ne0$, all the $q^{\rho_i}$'s are located {\rm inside} or {\rm on} the unite circle, and condition $(46)$ is fulfilled for those roots $z=a$ of the polynomial $(z-1)L(z)$ having the same location;

{\rm e)} $A_n\ne0$, all the $q^{\rho_i}$'s are located {\rm outside} or {\rm on} the unite circle, and condition $(46)$ is fulfilled for those roots $z=a$ of the polynomial $(z-1)L(z)$ having the same location.}
\bigskip

{\bf Theorem 7} (section 6). {\it Let the generalized power series $(6)$ satisfies equation $(5)$ and
$$
F'_{y_0}(x,\varphi,\mu\varphi,\ldots,\mu^n\varphi)=A_0x^{\nu}+o(x^\nu)\in{\cal G}, \qquad A_0\ne0.
$$
Then the series $\varphi$ uniformly converges in each sector $S$ of a small radius with the vertex at the origin and of the opening less than $2\pi$, defining there a holomorphic function.}

\end{document}


\begin{center}
\bf\large Сходимость (обобщённых) степенных рядов, удовлетворяющих функциональным уравнениям
\end{center}
\bigskip

\centerline{Р.\,Р.\,Гонцов, И.\,В.\,Горючкина}
\bigskip

\begin{abstract}
В работе изучаются вопросы, связанные со сходимостью обобщённых степенных рядов (с комплексными показателями степени), формально удовлетворяющих аналитическим функциональным уравнениям, --- дифференциальному, $q$-разностному, уравнению Малера. Мы представляем как новые результаты, так и обобщения результатов, полученных нами ранее, и тем самым подводим итог наших исследований по данной теме. Статья также содержит элементы обзора классических результатов о сходимости формальных степенных рядов Тейлора, удовлетворяющих таким уравнениям, и о существовании и единственности их локальных голоморфных решений. Эти результаты вкладываются в предложенные в настоящей работе.
\end{abstract}


\section{Введение}

{Решения нелинейных функциональных уравнений в общем случае не выражаются в виде конечного числа комбинаций и композиций элементарных и известных специальных функций. Одним из методов исследования решений таких уравнений является сначала поиск решений в виде {\it формальных} функциональных рядов ({\it формальных} решений) --- рядов, члены которых можно описать и как-либо упорядочить, но про сходимость которых изначально ничего неизвестно, а затем исследование сходимости или суммируемости этих формальных решений, то есть исследование вопросов существования  и единственности функций (настоящих решений) с заданным асимптотическим разложением в некоторой области. В данной работе мы исследуем только вопрос сходимости формальных функциональных рядов, имеющих вид бесконечной суммы степенных функций с комплексными, вообще говоря, показателями степени и удовлетворяющих аналитическим функциональным уравнениям следующих типов: дифференциальному, $q$-разностному или уравнению Малера. Но начнем с напоминания о более простом и классическом случае функционального уравнения --- об алгебраическом уравнении.

Как известно, все  формальные решения алгебраического уравнения
\begin{equation}\label{alg}
F(x,y)=0, \qquad F\in{\mathbb C}[x,y], \quad F(0,0)=0,
\end{equation}
в окрестности нуля  имеют, самое сложное, вид ряда Пюизо
\begin{equation}\label{pow}
y=\sum\limits_{k=1}^{+\infty}c_k\,x^{k/p}, \qquad c_k\in\mathbb{C},
\end{equation}
то есть степенного ряда с рациональными показателями степени, кратными некоторому $1/p$, $p\in\mathbb N$.
При этом метод многоугольника Ньютона \cite{Cheb} позволяет найти все разложения вида \eqref{pow} для уравнения \eqref{alg}.
Если заменить многочлен $F$ голоморфной возле начала координат функцией и искать формальные решения получившегося уравнения в классе бесконечных сумм степенных функций с произвольными показателями степени, то оказывается, что все формальные решения в этом классе также будут не более чем рядами Пюизо \eqref{pow}, и их все также можно найти с помощью метода многоугольника Ньютона (поскольку при вычислении степенных разложений в окрестности нуля важна только ограниченность многоугольника Ньютона слева и снизу). Отметим ещё, что всякий ряд Пюизо, формально удовлетворяющий аналитическому уравнению $F(x,y)=0$, имеет ненулевой радиус сходимости, что следует из обобщений аппроксимационной теоремы Артина \cite{Art}, полученных А.\,Плоски \cite{Pl}.

Если рассмотреть аналитическое {\it дифференциальное} уравнение
\begin{equation}\label{diff}
F(x,y,\delta y,\ldots,\delta^ny)=0, \qquad \displaystyle\delta=x\frac{d}{dx},
\end{equation}
где $F(x,y_0,y_1,\ldots,y_n)\not\equiv0$ --- голоморфная около $0\in{\mathbb C}^{n+2}$ функция, то рассматриваемые функциональные ряды, формально ему удовлетворяющие, могут уже иметь вещественные и комплексные показатели степени (например, моном $y=x^{\alpha}$, $\alpha\in\mathbb C$, удовлетворяющий уравнению $xy'-\alpha y=0$).

Аналогичное наблюдение можно сделать и относительно аналитических уравнений следующих двух типов: $q$-разностного и уравнения Малера. Первое --- уравнение вида
\begin{equation}\label{qdiff}
F(x,y(x),y(qx),\ldots,y(q^nx))=0, \qquad q\in{\mathbb C}\setminus\{0,1\};
\end{equation}
второе --- уравнение вида
\begin{equation}\label{Mahler}
F(x,y(x),y(x^\ell),\ldots,y(x^{\ell^n}))=0, \qquad \ell\in{\mathbb N}_{\geqslant2}.
\end{equation}
Начальными важными примерами появления таких уравнений являются задачи
\begin{itemize}
\item[--] о приведении ростка диффеоморфизма $f:({\mathbb C},0)\rightarrow({\mathbb C},0)$, $f(x)=qx+o(x)$, к линейной части $\tilde f(x)=qx$,
\item[--] о приведении ростка голоморфной в нуле функции $g(x)=ax^\ell+o(x^\ell)$ к главной части $\tilde g(x)=x^\ell$,
\end{itemize}
с помощью сопрягающих биголоморфных преобразований: $y^{-1}\circ f\circ y=\tilde f$ или $y^{-1}\circ g\circ y=\tilde g$. В случае существования такого преобразования $y$, оно должно удовлетворять, соответственно, $q$-разностному уравнению $y(qx)=f(y(x))$ или уравнению Малера $y(x^\ell)=g(y(x))$.

Также следует отметить, что степенной ряд, формально удовлетворяющий функциональному уравнению, может иметь нулевой радиус сходимости. Например, дифференциальное уравнение $$x^2y'-y+x=0$$ (пример Эйлера) имеет формальное решение в виде степенного ряда $$\sum\limits_{k=0}^{+\infty}k!\,x^{k+1},$$ расходящегося при всех $x\ne0$.

С учётом сделанных замечаний естественно исследовать вопрос об условиях сходимости функционального ряда вида
\begin{equation}\label{powgen}
\varphi=\sum\limits_{k=1}^{+\infty} c_kx^{\lambda_k}, \qquad \lambda_k\in{\mathbb C}, \quad 0<{\rm Re}\,\lambda_1\leqslant{\rm Re}\,\lambda_2\leqslant\ldots\rightarrow +\infty,
\end{equation}
формально удовлетворяющего одному из уравнений (\ref{diff}), (\ref{qdiff}), (\ref{Mahler}). Последнее означает, что после подстановки ряда в левую часть соответствующего уравнения (корректность которой мы поясним далее отдельно, см. \S3) получается ряд такого же вида с {\it нулевыми} коэффициентами.
\medskip

Ряды вида (\ref{powgen}) принято называть {\it обобщёнными} степенными рядами.
\medskip

Обобщённые степенные ряды нередко встречаются среди формальных решений функциональных уравнений.\footnote{Например, некоторые из дифференциальных и $q$-разностных уравнений Пенлеве имеют решения в виде обобщённых степенных рядов, равномерно сходящихся в секториальных окрестностях особых точек уравнений, что говорит о довольно сложном локальном поведении решений.} Алгоритмы вычисления решений в виде формальных степенных рядов (от рядов Тейлора до обобщённых степенных рядов) для обыкновенных дифференциальных уравнений хорошо развиты (см., например, \cite{Fine}--\cite{Bruno}), а для $q$-раз\-ност\-ных уравнений и уравнений Малера алгоритмы вычисления таких решений будут им аналогичны.

Вдохновлённые следующими тремя теоремами о сходимости классических степенных рядов (Тейлора)
\begin{equation}\label{powclass}
\varphi=\sum\limits_{k=1}^{+\infty} c_kx^k\in{\mathbb C}[[x]],
\end{equation}
формально удовлетворяющих дифференциальному, $q$-разностному или уравнению Малера, в настоящей работе мы предлагаем их аналоги --- достаточные условия сходимости --- для {\it обобщённых} степенных рядов, формально удовлетворяющих таким уравнениям.
\medskip

{\bf Теорема 1} (Б.\,Мальгранж \cite{Ma}). {\it Если формальный степенной ряд $(\ref{powclass})$ удовлетворяет уравнению $(\ref{diff})$, частная производная $F'_{y_n}(x,\varphi,\delta\varphi,\ldots,\delta^n\varphi)$ вдоль него отлична от нулевого ряда и при всех $j=0,1,\ldots,n-1$ выполнено
$$
{\rm ord}_{\,0\,}F'_{y_j}(x,\varphi,\delta\varphi,\ldots,\delta^n\varphi)\geqslant{\rm ord}_{\,0\,}F'_{y_n}(x,\varphi,\delta\varphi,\ldots,\delta^n\varphi),
$$
то ряд $(\ref{powclass})$ имеет ненулевой радиус сходимости.\footnote{Напомним, что \textit{порядком} формального степенного ряда \eqref{powclass} называется число
${\rm ord}_{\,0\,}\varphi={\rm min}\{k:\,c_k\neq 0\}.$}}
\medskip

Бернард Мальгранж, выдающийся французский математик, внесший огромный вклад в теорию линейных уравнений в частных производных, дифференциальную геометрию, теорию особенностей функций и отображений, теорию нелинейных дифференциальных уравнений и другие разделы математики, умер в возрасте 95 лет 5 января 2024 года. Теорема 1, сформулированная выше, стала классической и фундаментальной в аналитической теории дифференциальных уравнений. Идеи её лаконичного доказательства, основанного на теореме о неявном отображении для банаховых пространств, широко используются при получении подобных результатов в данной и смежных областях (в частности, в доказательстве приведённой ниже теоремы 2). Мы хотим особо отметить, что работа \cite{Ma} также оказала на нас большое влияние.
\medskip

{\bf Теорема 2.} {\it Обозначим $\Phi=\bigl(\varphi(x),\varphi(qx),\ldots,\varphi(q^nx)\bigr)$.}

1) {\it Пусть $|q|>1$ $($соответственно, $|q|<1)$ и формальный степенной ряд $(\ref{powclass})$ удовлетворяет уравнению $(\ref{qdiff}): F(x,\Phi)=0$. Если $F'_{y_n}(x,\Phi)\ne0$ $($соответственно, $F'_{y_0}(x,\Phi)\ne0)$ и при всех $j=0,1,\ldots,n-1$ выполнено
$$
{\rm ord}_{\,0\,}F'_{y_j}(x,\Phi)\geqslant{\rm ord}_{\,0\,}F'_{y_n}(x,\Phi)
$$
$($соответственно, при всех $j=1,\ldots,n$ выполнено ${\rm ord}_{\,0\,}F'_{y_j}(x,\Phi)\geqslant{\rm ord}_{\,0\,}F'_{y_0}(x,\Phi))$, то ряд $(\ref{powclass})$ имеет ненулевой радиус сходимости {\rm(Ш.\,Жанг \cite{zh1}, см. также \cite{lz}).}

$2)$ Пусть $|q|=1$ и формальный степенной ряд $(\ref{powclass})$ удовлетворяет уравнению $(\ref{qdiff})$. В предположении, что $m=\min_{0\leqslant j\leqslant n}{\rm ord}_{\,0\,}F'_{y_j}(x,\Phi) <\infty$, представим каждый $F'_{y_j}(x,\Phi)$ в виде
$$
F'_{y_j}(x,\Phi)=A_jx^m+o(x^m),\footnote{Здесь $o(x^m)$ обозначает формальный степенной ряд, порядок которого больше $m$.} \qquad j=0,1,\ldots,n,
$$
и определим ненулевой многочлен $P(z)=A_nz^n+\ldots+A_1z+A_0$. Если найдутся положительные постоянные $c,\nu$ такие, что для каждого корня $z=a$ многочлена $(z-1)P(z)$ число $q$ удовлетворяет условию
\begin{equation}\label{dioph}
|q^k-a|\geqslant c\,k^{-\nu}, \qquad k\in{\mathbb N},
\end{equation}
то ряд $(\ref{powclass})$ имеет ненулевой радиус сходимости} (Л.\,Ди Визио \cite{dv}).
\medskip

{\bf Теорема 3} (Ж.-П. Безивен \cite{Bez}\footnote{В оригинальной работе Безивена рассматривается полиномиальное уравнение (функция $F$ --- полином по переменным $y_0,y_1,\ldots,y_n$). Для аналитического уравнения теорема доказана в \cite{GG2}.}).
{\it Формальный степенной ряд $(\ref{powclass})$, удовлетворяющий уравнению $(\ref{Mahler})$, имеет ненулевой радиус сходимости.}
\medskip

Отметим, что варианты теорем существования и единственности голоморфных в окрестности нуля решений уравнений (\ref{diff}), (\ref{qdiff}), (\ref{Mahler}), о которых мы поговорим подробнее в следующем параграфе, во многом пересекаются с теоремами 1, 2, 3 и где-то являются их следствиями.

Как мы уже сказали, в настоящей работе исследуется только вопрос сходимости формальных решений (\ref{powgen}) уравнений (\ref{diff}), (\ref{qdiff}), (\ref{Mahler}), однако, в качестве важного комментария добавим здесь, что для расходящихся классических степенных рядов (\ref{powclass}), которые удовлетворяют дифференциальному или $q$-разностному уравнению, существуют оценки роста их коэффициентов (теоремы типа Майе \cite{Ma}, \cite{zh1}). Более того, для таких рядов развита теория (мульти)суммируемости, имеющая завершённый вид в дифференциальном случае \cite{Bra}, \cite{RS}, и представленная в завершённом виде пока только для линейных уравнений в $q$-разностном случае \cite{zh2}, \cite{mz}.
\medskip

{\bf Основные результаты} данной работы --- расширение теорем 1, 2, 3 на случай обобщённых степенных рядов --- формулируются и доказываются в \S\S 4, 5, 6, соответственно. В \S3 излагаются общие свойства обобщённых степенных рядов, удовлетворяющих рассматриваемым функциональным уравнениям, в частности, даётся описание множества показателей степени таких рядов. Отметим, что хотя $q$-разностное уравнение и уравнение Малера являются функциональными уравнениями определённых специальных видов, в изучаемых нами вопросах уравнений этих двух типов достаточно и для того чтобы исследовать подобные вопросы в большинстве случаев функциональных уравнений более общего вида,
$$
F(x,y(x),y(f(x)),\ldots,y(f^{[n]}(x)))=0,
$$
где $f$ --- голоморфная в нуле функция, $f(0)=0$, $f^{[j]}:=\underbrace{f\circ\ldots\circ f}_j$ ($f^{[0]}={\rm id}$, то есть $f^{[0]}(x)=x$). Об этом мы говорим подробнее в заключительном \S7.

\section{Теоремы существования и единственности локальных голоморфных решений функциональных уравнений}

Этот параграф отчасти продолжает введение и всё ещё связан с классическими степенными рядами Тейлора. В нём мы показываем, как теоремы существования и единственности локального голоморфного решения дифференциального, $q$-разностного или уравнения Малера общего вида могут быть выведены из теорем 1, 2 или 3, соответственно.

\subsection{Аналитическое дифференциальное уравнение}

Рассмотрим аналитическое дифференциальное уравнение, представленное в несколько более привычном виде, чем (\ref{diff}),
\begin{equation}\label{diff2}
F(x,y,y',\ldots,y^{(n)})=0,
\end{equation}
где $F(x,y_0,y_1,\ldots,y_n)\not\equiv0$ --- голоморфная около точки $(0,a_0,a_1,\ldots,a_n)\in{\mathbb C}^{n+2}$ функция. Для того чтобы перейти от (\ref{diff2}) к уравнению вида (\ref{diff}), определим новую зависимую переменную $u$:
$$
y=a_0+a_1x+\frac{a_2}{2!}x^2+\ldots+\frac{a_n}{n!}x^n+x^n\,u=:s_n+x^n\,u.
$$
Тогда, соответственно,
$$
y'=s_n'+x^{n-1}(\delta+n)u, \quad y''=s_n''+x^{n-2}(\delta+n-1)(\delta+n)u, \quad\ldots,\quad
$$
$$
y^{(n)}=s_n^{(n)}+(\delta+1)\ldots(\delta+n)u=a_n+(\delta+1)\ldots(\delta+n)u,
$$
где, напомним, $\displaystyle \delta=x\frac{d}{dx}$. Подставив эти соотношения в (\ref{diff2}), получим уравнение
\begin{equation}\label{diff3}
G(x,u,\delta u,\ldots,\delta^nu)=0,
\end{equation}
где голоморфная около $0\in{\mathbb C}^{n+2}$ функция $G(x,u_0,u_1,\ldots,u_n)\not\equiv0$ получается из $F(x,y_0,y_1,\ldots,y_n)$ заменой переменных
$$
y_0=s_n+x^nu_0, \quad y_1=s_n'+x^{n-1}(u_1+nu_0), \quad y_2=s_n''+x^{n-2}(u_2+\alpha_1^{(2)}u_1+\alpha_0^{(2)}u_0), \quad\ldots,\quad
$$
$$
y_n=a_n+u_n+\sum_{j=0}^{n-1}\alpha_j^{(n)}u_j, \qquad \alpha_j^{(k)}\in\mathbb N.
$$
Следовательно, при всех $j=0,1,\ldots,n$ имеем
$$
G'_{u_j}(x,u_0,u_1,\ldots,u_n)=F'_{y_j}(x,y_0,y_1,\ldots,y_n)\frac{\partial y_j}{\partial u_j}+\sum_{k>j}F'_{y_k}(x,y_0,y_1,\ldots,y_n)\frac{\partial y_k}{\partial u_j}=
$$
$$
=F'_{y_j}(x,y_0,y_1,\ldots,y_n)x^{n-j}+\sum_{k>j}F'_{y_k}(x,y_0,y_1,\ldots,y_n)\alpha_j^{(k)}x^{n-k}.
$$
Тем самым, если степенной ряд $\displaystyle \psi=s_n+x^n\sum\limits_{k\geqslant1}c_kx^k$ формально удовлетворяет уравнению (\ref{diff2}), то степенной ряд $\displaystyle \varphi=\sum\limits_{k\geqslant 1}c_kx^k$ формально удовлетворяет уравнению (\ref{diff3}), при этом
$$
G'_{u_j}(x,\varphi,\delta\varphi,\ldots,\delta^n\varphi)=F'_{y_j}(x,\psi,\psi',\ldots,\psi^{(n)})x^{n-j}+\sum_{k>j}F'_{y_k}(x,\psi,\psi',\ldots,\psi^{(n)})\alpha_j^{(k)}x^{n-k}
$$
при всех $j=0,1,\ldots,n$. Таким образом, вариант Теоремы 1 для аналитического дифференциального уравнения (\ref{diff2}) имеет следующий вид.
\medskip

{\bf Теорема $\bf 1^\prime$.} {\it Если формальный степенной ряд $\displaystyle \psi=s_n+x^n\sum_{k\geqslant1}c_kx^k$ удовлетворяет уравнению $(\ref{diff2})$, $F'_{y_n}(x,\psi,\psi',\ldots,\psi^{(n)})\ne0$, и при всех $j=0,1,\ldots,n-1$ выполнено
$$
{\rm ord}_{\,0\,}F'_{y_j}(x,\psi,\psi',\ldots,\psi^{(n)})+n-j\geqslant{\rm ord}_{\,0\,}F'_{y_n}(x,\psi,\psi',\ldots,\psi^{(n)}),
$$
то ряд $\psi$ имеет ненулевой радиус сходимости.}
\medskip

Приведём также теорему существования и единственности локального голоморфного решения дифференциального уравнения (\ref{diff2}), т.~е. уравнения, не разрешённого относительно старшей производной:
\medskip

{\it Если $F(0,a_0,a_1,\ldots,a_n)=0$ и $F'_{y_n}(0,a_0,a_1,\ldots,a_n)\ne0$, то уравнение $(\ref{diff2})$ имеет и при том единственное решение $y=y(x)$, голоморфное в точке $x=0$ и такое, что $y^{(j)}(0)=a_j$, $j=0,1,\ldots,n$.}
\medskip

Данное утверждение можно доказать с помощью аналитического варианта теоремы о неявной функции \cite[гл.\,I, \S4]{Sh} и теоремы Коши существования и единственности голоморфного решения дифференциального уравнения \cite[гл.\,I, \S5]{Gol}: сначала с помощью теоремы о неявной функции переменная $y_n$ однозначно выражается как функция $y_n=f(x,y_0,$ $\ldots,$ $y_{n-1})$ переменных $x,y_0,\ldots,y_{n-1}$, голоморфная в точке $(0,a_0,\ldots,a_{n-1})\in{\mathbb C}^{n+1}$ и такая, что
$$
a_n=f(0,a_0,a_1,\ldots,a_{n-1}), \qquad F\bigl(x,y_0,\ldots,y_{n-1},f(x,y_0,\ldots,y_{n-1})\bigr)\equiv0,
$$
а затем к разрешенному относительно старшей производной дифференциальному уравнению
\begin{equation}\label{razresh}
y^{(n)}=f(x,y,y',\ldots,y^{(n-1)})
\end{equation}
применяется теорема Коши, согласно которой уравнение \eqref{razresh} обладает единственным голоморфным решением $y=y(x)$, удовлетворяющим начальным условиям $y^{(j)}(0)=a_j$, $j=0,1,\ldots,n-1$. Это решение также будет удовлетворять уравнению (\ref{diff2}) и условию $y^{(n)}(0)=a_n$, а его единственность будет следовать из однозначности разрешения уравнения (\ref{diff2}) относительно старшей производной при заданных условиях.

Отметим, что теорему существования и единственности локального голоморфного решения дифференциального уравнения (\ref{diff2}) также можно доказать с помощью теоремы~$1^\prime$. Для этого достаточно доказать существование и единственность формального ряда $\displaystyle\psi=s_n+x^n\sum_{k\geqslant1}c_kx^k$, удовлетворяющего (\ref{diff2}), а тогда его сходимость будет следовать из теоремы~$1^\prime$: поскольку ряд $F'_{y_n}(x,\psi,\psi',\ldots,\psi^{(n)})$ начинается с ненулевого свободного слагаемого $F'_{y_n}(0,a_0,a_1,\ldots,a_n)$, то ${\rm ord}_{\,0\,}F'_{y_n}(x,\psi,\psi',\ldots,\psi^{(n)})=0$.
В свою очередь, доказательство существования и единственности формального ряда $\psi$, удовлетворяющего (\ref{diff2}), сводится к доказательству существования и единственности формального ряда $\displaystyle\varphi=\sum_{k\geqslant1}c_kx^k$, $\psi=s_n+x^n\varphi$, который удовлетворяет уравнению
$$
F\bigl(x,s_n+x^n\varphi,s_n'+x^{n-1}(\delta+n)\varphi,\ldots,s^{(n)}_n+(\delta+1)\ldots(\delta+n)\varphi\bigr)=0.
$$
Применяя формулу Тейлора в точке $(0,a_0,a_1,\ldots,a_n)$ к левой части этого уравнения, приходим к соотношению
\begin{equation}\label{Teylor1}
F'_{y_n}(0,a_0,a_1,\ldots,a_n)L(\delta)\varphi=x\,H(x,\varphi,\delta\varphi,\ldots,\delta^n\varphi),
\end{equation}
где $L(\delta)=(\delta+1)\ldots(\delta+n)$ и $H$ --- голоморфная функция в точке $0\in{\mathbb C}^{n+2}$. Коэффициент $F'_{y_n}(0,a_0,a_1,\ldots,a_n)L(k)c_k$ при каждой степени $x^k$, $k\geqslant1$, степенного ряда в левой части полученного соотношения (\ref{Teylor1}) однозначно выражается через $a_0,a_1,\ldots,a_n$ и $c_1,\ldots,c_{k-1}$ --- числа, определяющие коэффициент при соответствующей степени $x^k$ степенного ряда в правой части (\ref{Teylor1}), что доказывает существование и единственность искомого ряда $\varphi$.

\subsection{Аналитическое $q$-разностное уравнение}

Рассмотрим аналитическое $q$-разностное уравнение (\ref{qdiff}) с голоморфной в точке $0\in{\mathbb C}^{n+2}$ функцией $F(x,y_0,y_1,\ldots,y_n)$, обращающейся в ней в нуль. Определим многочлен
\begin{equation}\label{polyP}
P(z)=\sum_{j=0}^nA_jz^j, \qquad A_j:=F'_{y_j}(0).
\end{equation}
Имеет место следующая теорема существования и единственности локального голоморфного решения уравнения (\ref{qdiff}):
\medskip

{\it Если
\begin{itemize}
\item[--] $|q|>1$, $A_n\ne0$ и $P(q^k)\ne0$ при всех $k\in\mathbb N$, или
\item[--] $|q|<1$, $A_0\ne0$ и $P(q^k)\ne0$ при всех $k\in\mathbb N$, или
\item[--] $|q|=1$, $P\not\equiv0$ и найдутся положительные постоянные $c,\nu$ такие, что для каждого корня $z=a$ многочлена $(z-1)P(z)$ число $q$ удовлетворяет условию
$$
|q^k-a|\geqslant c\,k^{-\nu}, \qquad k\in{\mathbb N},
$$
\end{itemize}
то уравнение $(\ref{qdiff})$ имеет и при том единственное решение $y=y(x)$, голоморфное в точке $x=0$ и такое, что $y(0)=0$.}
\medskip

Доказательство данного утверждения сводится к доказательству существования и единственности формального решения $\displaystyle\varphi=\sum_{k\geqslant1}c_kx^k$ уравнения (\ref{qdiff}) и его сходимости в окрестности нуля при заданных условиях. Представим функцию $F$ в виде
\begin{equation}\label{reprF}
F(x,y_0,y_1,\ldots,y_n)=F_0(x)+\sum_{j=0}^nA_jy_j+F_{\geqslant2}(x,y_0,y_1,\ldots,y_n),
\end{equation}
где $F_0$ --- голоморфная функция около точки $x=0$, $F_0(0)=0$, и $F_{\geqslant2}$ --- голоморфная функция около $0\in{\mathbb C}^{n+2}$, ряд Тейлора которой состоит из слагаемых однородных степеней $\geqslant2$. Обозначим через $\sigma$ $q$-разностный оператор
\begin{equation}\label{sigma}
\sigma: y(x)\mapsto y(qx),
\end{equation}
действующий на голоморфные в окрестности нуля функции и имеющий естественное продолжение на алгебру ${\mathbb C}[[x]]$ формальных степенных рядов: $\displaystyle\sigma\Bigl(\sum_{k\geqslant1}c_kx^k\Bigr)=\sum_{k\geqslant1}c_kq^kx^k$. Тогда с учётом представления (\ref{reprF}) и введённых обозначений (\ref{polyP}), (\ref{sigma}) соотношение, которому должен удовлетворять искомый формальный степенной ряд $\varphi$, примет вид
$$
P(\sigma)\varphi=-F_0(x)-F_{\geqslant2}(x,\varphi,\sigma\varphi,\ldots,\sigma^n\varphi).
$$
Коэффициент $P(q^k)c_k$ при каждой степени $x^k$, $k\geqslant1$, степенного ряда в левой части этого соотношения однозначно выражается через $c_1,\ldots,c_{k-1}$ --- числа, определяющие коэффициент при соответствующей степени $x^k$ степенного ряда в его правой части, что, с учётом условия $P(q^k)\ne0$, доказывает существование и единственность искомого ряда $\varphi$. Сходимость ряда $\varphi$ в окрестности точки $x=0$ следует из теоремы 2:
\begin{itemize}
\item[--] при $|q|>1$ имеем ${\rm ord}_{\,0\,}F'_{y_n}(x,\Phi)=0$, поскольку ряд $F'_{y_n}(x,\Phi)$ начинается с ненулевого свободного слагаемого $F'_{y_n}(0)=A_n$;
\item[--] при $|q|<1$ имеем ${\rm ord}_{\,0\,}F'_{y_0}(x,\Phi)=0$, поскольку ряд $F'_{y_0}(x,\Phi)$ начинается с ненулевого свободного слагаемого $F'_{y_0}(0)=A_0$;
\item[--] при $|q|=1$ имеем $\min_{0\leqslant j\leqslant n}{\rm ord}_{\,0\,}F'_{y_j}(x,\Phi)=0$, $F'_{y_j}(x,\Phi)=A_j+o(1)$, и для корней ненулевого многочлена $P(u)=A_nu^n+\ldots+A_1u+A_0$ выполнено требуемое условие (\ref{dioph}).
\end{itemize}

\subsubsection{Пример: линеаризация локального диффеоморфизма}

В качестве частного случая применения теоремы существования и единственности локального голоморфного решения уравнения (\ref{qdiff}) рассмотрим её применение к уравнению Шрёдера --- $q$-разностному уравнению первого порядка --- с заданным начальным условием,
$$
y(qx)=f(y(x)),\qquad y(0)=0,
$$
где $f(x)=qx+\dots$ --- росток комплексного диффеоморфизма окрестности нуля $(\mathbb{C},0)$. Как уже было отмечено во введении, существование {\it биголоморфного} решения $y(x)$ этой задачи означает, что росток $f$ линеаризуем в группе $\rm G$ диффеоморфизмов $(\mathbb{C},0)$:
$$
y^{-1}\circ f\circ y(x)=qx, \qquad y\in\rm G.
$$
Поскольку многочлен (\ref{polyP}), в данном случае равный $P(z)=z-q$, обращается в нуль при $z=q$, то сделаем замену $y=x+xu$ зависимой переменной и перейдём к $q$-разностному уравнению вида
\begin{equation}\label{ex1}
u(qx)=u(x)+h(x,u(x)), \qquad u(0)=0,
\end{equation}
где $h$ --- голоморфная функция около $0\in{\mathbb C}^2$ и $h'_u(0)=0$. Соответствующий этому уравнению многочлен $P_1(z)=z-1$ не обращается в нуль в точках $z=q^k$, $k\in\mathbb N$, если число $q$ не является корнем из единицы. Поэтому из теоремы существования и единственности голоморфного решения, применённой к уравнению (\ref{ex1}), следует, что имеется единственный элемент $y\in \rm G$, $y'(0)=1$, который линеаризует росток $f$, если
\begin{itemize}
\item[--] $|q|\ne1$ (случай, исследованный А.\,Пуанкаре) или
\item[--] $|q|=1$ и найдутся положительные постоянные $c,\nu$ такие, что число $q$ удовлетворяет условию
\begin{equation}\label{siegel}
|q^k-1|\geqslant c\,k^{-\nu}, \qquad k\in{\mathbb N}
\end{equation}
(случай, исследованный К.\,Л.\,Зигелем \cite{Siegel}).
\end{itemize}
Таким образом, в этих двух случаях множество ${\rm G}_q\subset\rm G$ диффеоморфизмов, чья производная в точке $x=0$ равна $q$, представляет собой класс эквивалентности в $\rm G$, --- всякий $f\in{\rm G}_q$ линеаризуем.

Скажем ещё несколько слов про случай $q=e^{2\pi{\rm i}\omega}$, $\omega\in(0,1)\setminus\mathbb Q$. Условие \eqref{siegel}, как несложно заметить, эквивалентно  следующему\footnote{Неравенства вида \eqref{siegeladd} иногда называют {\it диофантовыми} неравенствами \cite{Milnor}. Как самостоятельный объект они рассматриваются в теории чисел и связаны с задачами аппроксимации иррациональных чисел рациональными (см. \cite{Khinchin}, \cite{Gelfond}, \cite{Khovanskij}). Если неравенство \eqref{siegeladd} выполняется для всех $m\in\mathbb Z$ и $k\in\mathbb N$, то $\omega$ называют {\it диофантовым} числом. В частности, все алгебраические числа диофантовы. С неравенствами такого типа мы еще встретимся в \S5 данной работы.}:
\begin{equation}\label{siegeladd}
|k\omega-m|\geqslant c\,k^{-\nu}\qquad k\in{\mathbb N},\; m\in\mathbb Z.
\end{equation}
Рассмотрев последовательность подходящих дробей $p_j/q_j$ непрерывной дроби, представляющей число $\omega$, можно, в свою очередь, установить эквивалентность условий \eqref{siegeladd} и
\begin{equation}\label{siegelconv}
\ln q_{j+1}=O(\ln q_j)\qquad \mbox{ при }j\to\infty
\end{equation}
(см. \cite[теор.\,9, 13, 16]{Khinchin}). Условие \eqref{siegelconv}, эквивалентное условию Зигеля (\ref{siegel}) и достаточное для того чтобы множество ${\rm G}_q$ было классом эквивалентности в $\rm G$, впоследствии было ослаблено А.\,Д.\,Брюно \cite{Br1}, \cite[\S5.II, теор.\,6]{Br2} и Х.\,Рюссманном \cite{Russ}, получившими менее ограничительное условие
\begin{equation}\label{brunoruss}
\sum_{j=0}^\infty\frac{\ln q_{j+1}}{q_j}<+\infty
\end{equation}
(выполнение последнего условия при выполнении \eqref{siegelconv} следует из показательного роста последовательности $\{q_j\}$). Наконец, Ж.-К.\,Йоккоз \cite{Yoccoz} доказал, что условие (\ref{brunoruss}) на число $q$ является не только достаточным, но и необходимым для того чтобы множество ${\rm G}_q$ было классом эквивалентности в $\rm G$ (то есть, для линеаризуемости {\it каждого} элемента ${\rm G}_q$).
\medskip

{\bf Замечание 1.} Условие (\ref{brunoruss}) имеет также эквивалентный вид
\begin{equation}\label{brunocond}
\sum_{j=0}^\infty\frac{\ln \Omega^{-1}(2^{j+1})}{2^j}<+\infty,
\end{equation}
где $\Omega(n)=\min\limits_{1\leqslant k\leqslant n, m\in\mathbb Z}|k\omega-m|$ (см. \cite[\S5.IV]{Br2}). Иногда условие (\ref{brunoruss}), возможно, для наглядности восприятия, формулируют в виде {\it субэкспоненциальности} роста функции $\Omega^{-1}(n)$:
$$
\Omega^{-1}(n)\leqslant Ce^{n^{1-\varepsilon}},
$$
при некоторых $C>0$, $0<\varepsilon<1$. Формально последнее условие всё же более ограничительно и, вообще говоря, может не выполняться при выполнении (\ref{brunocond}). Например, если $\Omega^{-1}(n)\leqslant Ce^{n/\log_2^2n}$, то субэкспоненциальность не гарантирована, поскольку $\log_2^2n$ растёт медленнее любой степени $n^\varepsilon$, $\varepsilon>0$. В то же время,
$$
\Omega^{-1}(2^{j+1})\leqslant Ce^{2^{j+1}/(j+1)^2}\;\Longrightarrow\;\frac{\ln \Omega^{-1}(2^{j+1})}{2^j}\leqslant\frac{C_1}{(j+1)^2},
$$
что влечёт (\ref{brunocond}).
\medskip

{\bf Замечание 2.} Если $q=e^{2\pi{\rm i}p/r}$ является корнем степени $r$ из единицы, то теорема существования и единственности голоморфного решения не применима и после перехода к уравнению (\ref{ex1}), поскольку многочлен $P_1(z)=z-1$ обращается в нуль в точках $q^{jr}=1$, $j\in\mathbb N$. Данный случай исследован Ж.\,Экалем \cite{Ec} и С.\,М.\,Ворониным \cite{Vo}. Линеаризуемыми диффеоморфизмами $f$ являются только элементы множества ${\rm G}_q$ порядка $r$: необходимость условия $f^{[r]}=\rm id$ для выполнения равенства $y(qx)=f(y(x))$ очевидна; достаточность можно показать, определив линеаризующее преобразование $y(x)$ как обратное к
$$
y^{-1}(x):=\sum\limits_{k=0}^{r-1}q^{-k}f^{[k]}(x)=rx+o(x).
$$
Остальные (то есть нелинеаризуемые) элементы множества ${\rm G}_q$ разбиваются на классы эквивалентности, параметризуемые функциональными модулями (ростками диффеоморфизмов).

\subsection{Аналитическое уравнение Малера}

Рассмотрим аналитическое уравнение Малера (\ref{Mahler}) с голоморфной в точке $0\in{\mathbb C}^{n+2}$ функцией $F(x,y_0,y_1,\ldots,y_n)$, обращающейся в ней в нуль. Имеет место следующая теорема существования и единственности локального голоморфного решения уравнения (\ref{Mahler}):
\medskip

{\it Если $F'_{y_0}(0)\ne0$, то уравнение $(\ref{Mahler})$ имеет и при том единственное решение $y=y(x)$, голоморфное в точке $x=0$ и такое, что $y(0)=0$.}
\medskip

Согласно теореме 3, для доказательства данного утверждения достаточно доказать существование и единственность формального решения $\displaystyle\varphi=\sum_{k\geqslant1}c_kx^k$ уравнения (\ref{Mahler}) при заданном условии. Представим функцию $F$ (поделив, если нужно, на $F'_{y_0}(0)\ne0$) в виде
\begin{equation}\label{reprF2}
F(x,y_0,y_1,\ldots,y_n)=\alpha x+y_0+L(y_1,\ldots,y_n)+F_{\geqslant2}(x,y_0,y_1,\ldots,y_n),
\end{equation}
где $\alpha\in\mathbb C$, $L$ --- линейная функция и $F_{\geqslant2}$ --- голоморфная функция около $0\in{\mathbb C}^{n+2}$, ряд Тейлора которой состоит из слагаемых однородных степеней $\geqslant2$. Обозначим через $\mu$ оператор Малера
\begin{equation}\label{mu}
\mu: y(x)\mapsto y(x^\ell),
\end{equation}
действующий на голоморфные в окрестности нуля функции и имеющий естественное продолжение на алгебру ${\mathbb C}[[x]]$ формальных степенных рядов: $\displaystyle\mu\Bigl(\sum_{k\geqslant1}c_kx^k\Bigr)=\sum_{k\geqslant1}c_kx^{\ell k}$. Тогда с учётом представления (\ref{reprF2}) и обозначения (\ref{mu}) уравнение (\ref{Mahler}), которому должен удовлетворять искомый формальный степенной ряд $\varphi$, примет вид
$$
\alpha x+y=-L(\mu y,\ldots,\mu^ny)-F_{\geqslant2}(x,y,\mu y,\ldots,\mu^ny).
$$
Следовательно, первый коэффициент ряда $\varphi$ определяется однозначно как $c_1=-\alpha$ (напомним, что $\ell\geqslant2$). Остальные коэффициенты ряда $\varphi$ также определяются однозначно, что можно показать с помощью замены зависимой переменной
$$
y=-\alpha x+xu\; \Longrightarrow \; \mu y=-\alpha x^{\ell}+x^{\ell}\mu u,\;\ldots,\;\mu^n y=-\alpha x^{\ell^n}+x^{\ell^n}\mu^n u,
$$
приводящей последнее уравнение к виду
$$
xu=x^2\,G(x,u,\mu u,\ldots,\mu^n u),
$$
то есть,
\begin{equation}\label{Mahlerreduced}
u=x\,G(x,u,\mu u,\ldots,\mu^n u),
\end{equation}
где $G$ --- голоморфная функция около $0\in{\mathbb C}^{n+2}$. Коэффициенты ряда $\psi=\alpha+\varphi/x=\displaystyle\sum_{k\geqslant1}c_{k+1}x^k$, удовлетворяющего (\ref{Mahlerreduced}),
$$
\psi=x\,G(x,\psi,\mu\psi,\ldots,\mu^n\psi),
$$
определяются однозначно: $c_2=G(0)$ и коэффициент при каждой степени $x^k$, $k\geqslant2$, степенного ряда в правой части этого соотношения однозначно выражается через $c_2,\ldots,c_k$. Это завершает доказательство утверждения о существовании и единственности локального голоморфного решения уравнения Малера (\ref{Mahler}).

\subsubsection{Пример: приведение ростка голоморфной функции к главной части}

Приведём пример применения теоремы существования и единственности локального голоморфного решения уравнения (\ref{Mahler}) к следующему уравнению Малера первого порядка с заданным начальным условием:
$$
y(x^\ell)=g(y(x)),\qquad y(0)=0,
$$
где $g(x)=ax^\ell+\dots$, $a\in{\mathbb C}^*$, --- росток голоморфной неоднолистной ($\ell\geqslant2$) в окрестности нуля функции. Как уже было отмечено во введении, существование биголоморфного решения $y(x)$ этой задачи означает, что росток $g$ сопряжением приводится к своей главной части $\tilde g(x)=x^\ell$.
Указанная теорема не может быть применена непосредственно к данному уравнению, поскольку $F(x,y_0,y_1)=y_1-g(y_0)$ в данном случае и тогда $F'_{y_0}(0)=g'(0)=0$. Поэтому для доказательства существования сопрягающего диффеоморфизма $y$ сделаем замену $y=cx+xu$ зависимой переменной, где $c$ --- один из корней уравнения $c^{\ell-1}=a^{-1}$. Получим уравнение Малера относительно неизвестной $u$,
$$
cx^\ell+x^\ell u(x^\ell)=ax^\ell(c+u(x))^\ell+\ldots,
$$
после упрощения принимающее вид
\begin{equation}\label{ex2}
u(x^\ell)=\ell u(x)+h(x,u(x)), \qquad u(0)=0,
\end{equation}
где $h$ --- голоморфная функция около $0\in{\mathbb C}^2$ и $h'_u(0)=0$. Таким образом, из теоремы существования и единственности голоморфного решения, применённой к уравнению (\ref{ex2}), следует, что росток $g$ сопряжением приводится к своей главной части $\tilde g(x)=x^\ell$, при этом сопрягающий диффеоморфизм $y$, удовлетворяющий условию $y'(0)=c$, единственен.

\section{Обобщённые степенные ряды, формально удовлетворяющие аналитическому функциональному уравнению}

В этом параграфе мы поясним корректность подстановки обобщённого степенного ряда в функциональное уравнение какого-либо из рассматриваемых типов и, тем самым, корректность понятия {\it формального решения в виде обобщённого степенного ряда}. Также мы изучим структуру множества показателей степени такого формального решения.

\subsection{Алгебра обобщённых степенных рядов}

Прежде всего заметим, что множество
$$
{\cal G}=\Bigl\{\varphi=\sum\limits_{k=1}^{+\infty} c_kx^{\lambda_k} \mid \lambda_k\in{\mathbb C}, \quad 0<{\rm Re}\,\lambda_1\leqslant{\rm Re}\,\lambda_2\leqslant\ldots\rightarrow\infty\Bigr\}
$$
всех обобщённых степенных рядов является алгеброй. Действительно, сумма обобщённых степенных рядов является обобщённым степенным рядом, поскольку объединение множеств показателей степени таких рядов снова может быть упорядочено по неубыванию вещественных частей показателей. Произведение двух рядов
$$
\varphi=\sum\limits_{k=1}^{+\infty} c_kx^{\lambda_k}, \qquad \tilde\varphi=\sum\limits_{k=1}^{+\infty}\tilde c_kx^{\tilde\lambda_k}
$$
представляется в виде
$$
\varphi\,\tilde\varphi=\sum_{\nu\in\mathbb C}\biggl(\sum\limits_{\lambda_i+\tilde\lambda_j=\nu} c_i\tilde c_j\biggr)x^\nu,
$$
при этом внутренние суммы корректно определены, поскольку для каждого $\nu$ найдётся только конечное число пар $(\lambda_i,\tilde\lambda_j)$ таких, что $\lambda_i+\tilde\lambda_j=\nu$, в силу стремления последовательностей $\{{\rm Re}\,\lambda_k\}$, $\{{\rm Re}\,\tilde\lambda_k\}$ к бесконечности. По этой же причине для каждого $N>0$ имеется только конечное число пар $(\lambda_i,\tilde\lambda_j)$ таких, что ${\rm Re}\,\lambda_i+{\rm Re}\,\tilde\lambda_j<N$, то есть в ряду $\varphi\,\tilde\varphi$ конечное число показателей $\nu$ таких, что ${\rm Re}\,\nu<N$ и, следовательно, данный ряд также является обобщённым степенным рядом.

\subsection{Подстановка обобщённого степенного ряда в аналитическое функциональное уравнение. Формальное решение}

Напомним, что $\delta=x(d/dx)$, $\sigma$ и $\mu$ обозначают дифференциальный, $q$-разностный и оператор Малера, соответственно, действующие на пространстве голоморфных в окрестности нуля функций:
\begin{eqnarray*}
\delta: y(x) & \mapsto & xy'(x), \\ \sigma: y(x) & \mapsto & y(qx), \\ \mu: y(x) & \mapsto & y(x^\ell).
\end{eqnarray*}
Они естественным образом продолжаются на алгебру $\cal G$ обобщённых степенных рядов:
\begin{eqnarray*}
\delta\Bigl(\sum_{k=1}^{+\infty}c_kx^{\lambda_k}\Bigr)&=&\sum_{k=1}^{+\infty}c_k\lambda_kx^{\lambda_k}, \\
\sigma\Bigl(\sum_{k=1}^{+\infty}c_kx^{\lambda_k}\Bigr)&=&\sum_{k=1}^{+\infty} c_kq^{\lambda_k}x^{\lambda_k}, \\
\mu\Bigl(\sum_{k=1}^{+\infty}c_kx^{\lambda_k}\Bigr)&=&\sum_{k=1}^{+\infty}c_kx^{\ell\lambda_k},
\end{eqnarray*}
при этом продолжение оператора $\sigma$ на $\cal G$ подразумевает выбор фиксированного значения $\ln q$ (например, соответствующего выбору $0\leqslant\arg q<2\pi$). Таким образом, если $F(x,y_0,y_1,\ldots,y_n)$ --- многочлен, $\varphi\in\cal G$ и $\Delta$ --- один из операторов $\delta$, $\sigma$, $\mu$, то $F(x,\varphi,\Delta\varphi,\ldots,\Delta^n\varphi)\in\cal G$. Для корректности подстановки обобщённого степенного ряда не только в полиномиальное, но и в {\it аналитическое} функциональное уравнение, нам нужно показать, что то же самое справедливо и в случае если $F$ --- голоморфная функция в окрестности $0\in{\mathbb C}^{n+2}$.
\medskip

{\bf Лемма 1.} {\it Пусть $F(x,y_0,y_1,\ldots,y_n)$ --- голоморфная функция в окрестности $0\in{\mathbb C}^{n+2}$, $\varphi\in\cal G$ и $\Delta$ --- один из операторов $\delta$, $\sigma$, $\mu$. Тогда $F(x,\varphi,\Delta\varphi,\ldots,\Delta^n\varphi)\in\cal G$.}
\medskip

{Д о к а з а т е л ь с т в о.} Рассмотрев степенной ряд функции $F$ в окрестности точки $0\in{\mathbb C}^{n+2}$, представим её в виде
$$
F(x,y_0,y_1,\ldots,y_n)=F_1(x,y_0,y_1,\ldots,y_n)+F_2(x,y_0,y_1,\ldots,y_n)+\ldots,
$$
где $F_d$ --- многочлен, образованный слагаемыми однородной степени $d$. Тогда показатели степени каждого обобщённого степенного ряда $F_d(x,\varphi,\Delta\varphi,\ldots,\Delta^n\varphi)\in\cal G$,
$$
F_d(x,\varphi,\Delta\varphi,\ldots,\Delta^n\varphi)=\sum_{k=1}^{+\infty}c_k^{(d)}x^{\lambda_k^{(d)}}, \qquad d=1,2,\ldots,
$$
обладают следующим свойством: $0<d\min(1,{\rm Re}\,\lambda_1)<{\rm Re}\,\lambda_1^{(d)}\leqslant{\rm Re}\,\lambda_2^{(d)}\leqslant\ldots\rightarrow\infty$, где $\lambda_1$ --- первый показатель степени ряда $\varphi$. Таким образом, для каждого $N>0$ лишь конечное число показателей степени $\lambda_k^{(d)}$ в сумме $\displaystyle\sum\limits_{d=1}^{+\infty}F_d(x,\varphi,\Delta\varphi,\ldots,\Delta^n\varphi)$ удовлетворяют условию ${\rm Re}\,\lambda_k^{(d)}<N$ (те, для которых $d<N/\min(1,{\rm Re}\,\lambda_1)$) и, тем самым, $F(x,\varphi,\Delta\varphi,\ldots,\Delta^n\varphi)\in\cal G$. {\hfill$\Box$}
\medskip

Доказанная лемма делает корректным следующее определение формального решения.
\medskip

{\bf Определение.} Будем говорить, что обобщённый степенной ряд $\varphi\in\cal G$ является {\it формальным решением} одного из уравнений \eqref{diff}, \eqref{qdiff} или \eqref{Mahler}, если $F(x,\varphi,\Delta\varphi,\ldots,\Delta^n\varphi)$ --- обобщённый степенной ряд с нулевыми коэффициентами:
$$
F(x,\varphi,\Delta\varphi,\ldots,\Delta^n\varphi)=0\in\cal G,
$$
где $\Delta=\delta$, $\sigma$ или $\mu$, соответственно.
\medskip

{\bf Замечание 3.} Из доказательства леммы 1 следует, что $F$ может быть и формальным степенным рядом, $F\in{\mathbb C}[[x,y_0,y_1,\ldots,y_n]]$ (его сходимость в доказательстве не используется). Утверждение леммы также останется справедливым, если голоморфную функцию (формальный степенной ряд) $F$ заменить на ряд вида
$$
F(x,y_0,y_1,\ldots,y_n)=\sum_{{\bf p}\in{\mathbb Z}_+^{n+1}}b_{\bf p}(x)y_0^{p_0}y_1^{p_1}\ldots y_n^{p_n}, \qquad {\bf p}=(p_0,p_1,\ldots,p_n), \quad
b_{\bf p}\in{\cal G}.
$$

{\bf Обозначение 1.} Жирным шрифтом всегда будем обозначать целочисленные мультииндексы, при этом их размерность иногда указывается явно, а иногда понятна из контекста. Модуль мультииндекса, как обычно, обозначает сумму координат последнего. Например, $|{\bf p}|=p_0+p_1+\ldots+p_n$ для мультииндекса ${\bf p}$ из замечания 3.

Также будем использовать мультииндексные обозначения типа ${\bf m}>{\bf k}$, означающие, что $m_i\geqslant k_i$ при всех $i$ и $|{\bf m}|>|{\bf k}|$.

\subsection{Леммы о приведении функциональных уравнений к специальному виду}

В дальнейшем важную роль будут играть леммы о приведении функциональных уравнений \eqref{diff}, \eqref{qdiff}, \eqref{Mahler}, обладающих формальным решением в виде обобщённого степенного ряда, к специальному виду. В силу их значимости мы вынесем здесь эти леммы вместе в отдельный параграф.
\medskip

{\bf Обозначение 2.} Ниже в лемме 2 и далее в записях вида $ax^\nu+o(x^\nu)\in{\cal G}$, $o(x^\nu)$ обозначает формальный обобщённый степенной ряд, первый член которого имеет вещественную часть степени строго большую ${\rm Re}\,\nu$.
\medskip

{\bf Лемма 2.} {\it Пусть обобщённый степенной ряд \eqref{powgen} удовлетворяет уравнению \eqref{diff} или \eqref{qdiff} и
$$
F'_{y_j}(x,\varphi,\Delta\varphi,\ldots,\Delta^n\varphi)=A_jx^{\nu}+o(x^\nu)\in{\cal G},
$$
где $\Delta=\delta$ или $\sigma$, соответственно, число $\nu$ --- одинаковое для всех $j=0,1,\ldots,n$ и при этом хотя бы одно $A_j$ отлично от нуля. Тогда найдётся $N\in\mathbb N$ такое, что преобразование
$$
y=\sum_{k=1}^Nc_kx^{\lambda_k}+x^{\lambda_N}u
$$
приводит
\begin{itemize}
\item[\rm i)] уравнение \eqref{diff} к виду
\begin{eqnarray}\label{auxildiff}
L(\delta+\lambda_N)u=M(x,u,\delta u,\ldots,\delta^nu);
\end{eqnarray}
\item[\rm ii)] уравнение \eqref{qdiff} к виду
\begin{eqnarray}\label{auxilqdiff}
L(q^{\lambda_N}\sigma)u=M(x,u,\sigma u,\ldots,\sigma^nu),
\end{eqnarray}
\end{itemize}
где $L(z)=A_nz^n+\ldots+A_1z+A_0\not\equiv0$ --- многочлен степени $\leqslant n$,
$$
M(x,u_0,u_1,\ldots,u_n)=\sum_{{\bf p}\in{\mathbb Z}_+^{n+1}}a_{\bf p}(x)u_0^{p_0}u_1^{p_1}\ldots u_n^{p_n}
$$
--- голоморфная по переменным $u_j$ функция, представимая рядом c коэффициентами $a_{\bf p}\in\cal G$ --- обобщёнными степенными рядами с ненулевым радиусом сходимости, показатели степени которых принадлежат множеству $-\lambda_N-\nu+\Gamma'$, где $\Gamma'$ --- аддитивная полугруппа, порождённая числами $1,\lambda_1,\ldots,\lambda_N$.}
\medskip

{Д о к а з а т е л ь с т в о.} i) Для всякого $N\in\mathbb N$ такого, что ${\rm Re}\,\lambda_{N+1}>{\rm Re}\,\lambda_N$, представим формальный ряд $\varphi$ в виде
$$
\varphi=\sum_{k=1}^Nc_kx^{\lambda_k}+x^{\lambda_N}\psi=\varphi_N+x^{\lambda_N}\psi.
$$
Обозначая $\Phi=(\varphi,\delta\varphi,\ldots,\delta^n\varphi)=\Phi_N+x^{\lambda_N}\Psi$ и применяя (формальную) формулу Тейлора к соотношению $F(x,\Phi)=0$, получим
\begin{eqnarray}\label{diffTaylor}
0 & = & F(x,\Phi_N+x^{\lambda_N}\Psi)=F(x,\Phi_N)+x^{\lambda_N}\sum_{j=0}^n\frac{\partial F}{\partial y_j}(x,\Phi_N)\psi_j+\nonumber \\
  &  & +x^{2\lambda_N}\sum_{p,q=0}^nH_{p,q}(x,\Phi_N,x^{\lambda_N}\Psi)\psi_p\psi_q,
\end{eqnarray}
где $\psi_j=(\delta+\lambda_N)^j\psi\in\cal G$ и $H_{p,q}$ --- голоморфные функции возле $0\in{\mathbb C}^{2n+3}$. Здесь мы используем формулу Тейлора так же, как и в случае неформальных объектов. Её вывод основан на представлении функции $F$ в виде (бесконечной) суммы однородных полиномов положительных степеней и возможности применения соответствующих формул к каждому такому полиному.

Выберем число $N\in\mathbb N$ так, чтобы помимо условия ${\rm Re}\,\lambda_{N+1}>{\rm Re}\,\lambda_N$ выполнялось
\begin{equation}\label{redlemmacond1}
{\rm Re}\,\lambda_N>{\rm Re\,}\nu. 
\end{equation}
Покажем, что такое $N$ подойдёт для утверждения леммы.

Аналогично (\ref{diffTaylor}) получаем
$$
\frac{\partial F}{\partial y_j}(x,\Phi)=\frac{\partial F}{\partial y_j}(x,\Phi_N)+x^{\lambda_N}\sum_{k=0}^n\frac{\partial^2 F}
{\partial y_j\partial y_k}(x,\Phi_N)\psi_k+x^{2\lambda_N}\sum_{p,q=0}^n\widetilde H_{p,q}(x,\Phi_N,x^{\lambda_N}\Psi)\psi_p\psi_q,
$$
где $\widetilde H_{p,q}$ --- голоморфные функции возле $0\in{\mathbb C}^{2n+3}$. Следовательно, ввиду (\ref{redlemmacond1}),
$$
\frac{\partial F}{\partial y_j}(x,\Phi_N)=A_jx^{\nu}+o(x^\nu)\in{\mathbb C}\{x^{\Gamma'}\}\subset\cal G
$$
--- обобщённый степенной ряд\footnote{Обозначение ${\mathbb C}\{x^{\Gamma'}\}$ мы используем для подалгебры {\it сходящихся} (имеющих ненулевой радиус сходимости) обобщённых степенных рядов, показатели степени которых принадлежат конечно порождённой аддитивной полугруппе $\Gamma'$, в то время как ${\mathbb C}[[x^{\Gamma'}]]$ будет обозначать подалгебру {\it формальных} обобщённых степенных рядов с показателями из $\Gamma'$.}, начинающийся с того же слагаемого, что и $F'_{y_j}(x,\Phi)$. Поэтому соотношение (\ref{diffTaylor}) может быть разделено на $x^{\lambda_N+\nu}$, после чего получится равенство вида
$$
L(\delta+\lambda_N)\psi-M(x,\psi,\delta\psi,\ldots,\delta^n\psi)=0,
$$
где полином $L$ и функция $M$ --- как в утверждении леммы. Таким образом, преобразование
$$
y=\sum_{k=1}^Nc_kx^{\lambda_k}+x^{\lambda_N}u
$$
приводит уравнение (\ref{diff}) к уравнению
$$
L(\delta+\lambda_N)u=M(x,u,\delta u,\ldots,\delta^nu),
$$
обладающему формальным решением $u=\psi\in\cal G$.

ii) В случае $q$-разностного уравнения доказательство аналогично, с той лишь разницей, что здесь $\Phi=(\varphi,\sigma\varphi,\ldots,\sigma^n\varphi)$ и
$\psi_j=(q^{\lambda_N}\sigma)^j\psi$, $j=0,1,\ldots,n$. {\hfill $\Box$}
\medskip

{\bf Лемма 3.} {\it Пусть обобщённый степенной ряд \eqref{powgen} удовлетворяет уравнению \eqref{Mahler} и
$$
F'_{y_0}(x,\varphi,\mu\varphi,\ldots,\mu^n\varphi)=A_0x^{\nu_0}+o(x^{\nu_0})\in{\cal G}, \qquad A_0\ne0.
$$
Тогда найдётся $N\in\mathbb N$ такое, что преобразование
$$
y=\sum_{k=1}^Nc_kx^{\lambda_k}+x^{\lambda_N}u
$$
приводит уравнение \eqref{Mahler} к виду
\begin{eqnarray}\label{auxilMahler}
u=M(x,u,\mu u,\ldots,\mu^nu),
\end{eqnarray}
где функция $M$ --- такая же, как в лемме $2$ $($с заменой $\nu$ на $\nu_0)$.}
\medskip

{Д о к а з а т е л ь с т в о.} Как и в доказательстве леммы 2, для всякого $N\in\mathbb N$ такого, что ${\rm Re}\,\lambda_{N+1}>{\rm Re}\,\lambda_N$, представим формальный ряд $\varphi$ в виде $\varphi=\varphi_N+x^{\lambda_N}\psi$. Обозначив $\Phi=(\varphi,\mu\varphi,\ldots,\mu^n\varphi)=\Phi_N+x^{\lambda_N}\Psi$ и аналогично применяя формулу Тейлора, получим соотношение
\begin{eqnarray}\label{MahlerTaylor}
0 & = & F(x,\Phi_N+x^{\lambda_N}\Psi)=F(x,\Phi_N)+x^{\lambda_N}\sum_{j=0}^n\frac{\partial F}{\partial y_j}(x,\Phi_N)\psi_j+\nonumber \\
  &  & +x^{2\lambda_N}\sum_{p,q=0}^nH_{p,q}(x,\Phi_N,x^{\lambda_N}\Psi)\psi_p\psi_q,
\end{eqnarray}
где $\psi_j=x^{(\ell^j-1)\lambda_N}\mu^j\psi\in\cal G$ и $H_{p,q}$ --- голоморфные функции возле $0\in{\mathbb C}^{2n+3}$.

Введём обозначение $A_jx^{\nu_j}$ для слагаемого, с которого начинается обобщённый степенной ряд $F'_{y_j}(x,\Phi)$. Тем самым,
$$
+\infty>{\rm Re}\,\nu_0\geqslant\min_{j=0,1,\ldots,n}{\rm Re}\,\nu_j=:s.
$$

При достаточно большом $N$ каждый обобщённый степенной ряд $F'_{y_j}(x,\Phi_N)$ начинается с тех же слагаемых, что и $F'_{y_j}(x,\Phi)$ (см. доказательство леммы 2). Дополнительно выберем $N$ столь большим, чтобы выполнялось неравенство
$$
s+{\rm Re}(\ell-1)\lambda_N>{\rm Re}\,\nu_0.
$$
Получим соотношения
\begin{eqnarray*}
F'_{y_0}(x,\Phi_N)\psi_0&=&\bigl(A_0x^{\nu_0}+o(x^{\nu_0})\bigr)\psi,\\
F'_{y_j}(x,\Phi_N)\psi_j&=&\bigl(A_jx^{\nu_j+(\ell^j-1)\lambda_N}+\ldots\bigr)\mu^j\psi, \quad j=1,\ldots,n,
\end{eqnarray*}
при этом множители перед $\psi$, $\mu^j\psi$ суть элементы ${\mathbb C}\{x^{\Gamma'}\}$ и ${\rm Re}(\nu_j+(\ell^j-1)\lambda_N)>{\rm Re}\,\nu_0$. Поэтому равенство (\ref{MahlerTaylor}) может быть разделено на $A_0x^{\lambda_N+\nu_0}$, после чего получится соотношение вида
$$
\psi-M(x,\psi,\mu\psi,\ldots,\mu^n\psi)=0,
$$
где функция $M$ --- как в утверждении леммы. Таким образом, преобразование
$$
y=\sum_{k=1}^Nc_kx^{\lambda_k}+x^{\lambda_N}u
$$
приводит уравнение (\ref{Mahler}) к уравнению
$$
u=M(x,u,\mu u,\ldots,\mu^nu),
$$
обладающему формальным решением $u=\psi\in\cal G$. $\quad\Box$
\medskip

{\bf Замечание 4.} В случае когда показатели степени $\lambda_k$ обобщённого степенного ряда \eqref{powgen} {\it вещественны}, условие леммы 2 относительно вида частных производных $F'_{y_j}$ просто означает существование такого $j$, что данный ряд не является решением уравнения $F'_{y_j}(x,y,\Delta y,\ldots,\Delta^n y)=0$, а подобное условие леммы 3 означает, что $\varphi$ не является решением уравнения $F'_{y_0}(x,y,\mu y,\ldots,\mu^n y)=0$.

\subsection{Структура множества показателей степени обобщённого степенного ряда -- формального решения аналитического функционального уравнения}

Основная цель данного параграфа --- показать, что обобщённый степенной ряд $\varphi$, формально удовлетворяющий одному из уравнений (\ref{diff}), (\ref{qdiff}), (\ref{Mahler}), при некоторых условиях является элементом подалгебры ${\mathbb C}[[x^\Gamma]]\subset\cal G$, где $\Gamma$ --- конечно порождённая аддитивная полугруппа, образующие которой, во-первых, имеют положительные вещественные части и, во-вторых, линейно независимы над $\mathbb Z$.

Прежде всего заметим, что в доказанных выше леммах 2, 3 можно перейти от множества $-\lambda_N-\nu+\Gamma'$, участвующего в описании показателей степени коэффициентов $a_{\bf p}\in\cal G$, просто к конечно порождённой аддитивной полугруппе, образующие которой имеют положительные вещественные части.
\medskip

{\bf Лемма 4.} {\it В леммах $2$, $3$ коэффициенты $a_{\bf p}\in\cal G$ принадлежат подалгебре ${\mathbb C}\{x^{\Gamma''}\}$, где $\Gamma''$ --- конечно порождённая аддитивная полугруппа, образующие которой имеют положительные вещественные части.}
\medskip

{Д о к а з а т е л ь с т в о.} Элементы множества $-\lambda_N-\nu+\Gamma'$ суть числа вида
$$
m_0+\sum_{j=1}^Nm_j\lambda_j-\lambda_N-\nu, \qquad {\bf m}=(m_0,m_1,\ldots,m_N)\in{\mathbb Z}_+^{N+1}\setminus\{0\},
$$
с положительной вещественной частью. Множество 
$$
{\cal M}=\Bigl\{{\bf m}\in{\mathbb Z}_+^{N+1}\setminus\{0\}\mid {\rm Re}\Bigl(m_0+\sum_{j=1}^Nm_j\lambda_j-\lambda_N-\nu\Bigr)>0\Bigr\}
$$
содержит положительное конечное число минимальных элементов (относительно отношения частичного порядка ''$<$'' на множестве ${\mathbb Z}_+^{N+1}$ мультииндесов);\footnote{Данное утверждение известно под названием леммы Диксона, см., например, \cite{Dick}.} обозначим их ${\bf m}^{(1)},\ldots,{\bf m}^{(N_1)}$.

Тогда добавление к образующим $1,\lambda_1,\ldots,\lambda_N$ полугруппы $\Gamma'$ конечного набора образующих
$$
m_0^{(\alpha)}+\sum_{j=1}^Nm_j^{(\alpha)}\lambda_j-\lambda_N-\nu, \qquad \alpha=1,\ldots,N_1,
$$
приводит нас к нужной конечно порождённой аддитивной полугруппе $\Gamma''$:
$$
{\rm Re}\Bigl(m_0+\sum_{j=1}^Nm_j\lambda_j-\lambda_N-\nu\Bigr)>0\;\Longrightarrow\;m_0+\sum_{j=1}^Nm_j\lambda_j-\lambda_N-\nu\in\Gamma''.
$$
Действительно, для каждого ${\bf m}\in{\cal M}$ найдётся ${\bf m}^{(\alpha)}\in\{{\bf m}^{(1)},\ldots,{\bf m}^{(N_1)}\}$ такое, что ${\bf m}^{(\alpha)}\leqslant
{\bf m}$, и тогда
$$
m_0+\sum_{j=1}^Nm_j\lambda_j-\lambda_N-\nu=m_0^{(\alpha)}+\sum_{j=1}^Nm_j^{(\alpha)}\lambda_j-\lambda_N-\nu+(m_0-m_0^{(\alpha)})+
\sum_{j=1}^N(m_j-m_j^{(\alpha)})\lambda_j\in\Gamma''.
$$
{\hfill $\Box$}
\medskip

Дальнейший важный шаг -- переход от $\Gamma''$ к конечно порождённой аддитивной полугруппе $\Gamma$, образующие которой имеют положительные вещественные части и при этом линейно независимы над $\mathbb Z$. Это обеспечивается следующей леммой, доказанной Антоном А. Владимировым.
\medskip

{\bf Лемма 5.} {\it Пусть $\rm G$ --- аддитивная полугруппа, порождённая комплексными числами $r_1,\ldots,r_s$, вещественные части которых положительны. Найдутся комплексные числа $\rho_1,$ $\ldots,$ $\rho_{\tau}$, линейно независимые над $\mathbb Z$, такие что все ${\rm Re}\,\rho_i>0$ и аддитивная полугруппа $\Gamma$, порождённая этими числами, содержит $\rm G$.}
\medskip

{Д о к а з а т е л ь с т в о.} Пусть $r_1,\ldots,r_{\tau}$ --- максимальный набор линейно независимых над $\mathbb Z$ образующих полугруппы $\rm G$. Достаточно доказать, что при добавлении к этому набору произвольного числа $b$, ${\rm Re}\,b>0$, такого что $r_1,\ldots,r_{\tau},b$ становятся линейно зависимыми над $\mathbb Z$, полугруппа $\rm G'$, порождённая числами $r_1,\ldots,r_{\tau},b$, содержится в некоторой полугруппе $\Gamma$, порождённой $\tau$ линейно независимыми над $\mathbb Z$ комплексными числами $\rho_1,\ldots,\rho_{\tau}$ с положительными вещественными частями.

Можно считать, что
\begin{eqnarray}\label{b}
b=m_1r_1+\ldots+m_{\tau-j}r_{\tau-j}-m_{\tau-j+1}r_{\tau-j+1}-\ldots-m_{\tau}r_{\tau}, \;\;\;m_i\in{\mathbb N}, \;\; 1\leqslant j\leqslant\tau-1.
\end{eqnarray}
Вообще говоря, некоторые $m_i$ могут быть равными нулю. В таком случае соответствующие образующие $r_i$ остаются без изменений и становятся образующими новой полугруппы $\Gamma$, то есть $\rho_i=r_i$ для таких $m_i$. Читатель может сделать нужные изменения в дальнейших рассуждениях, относящихся к ситуации, когда все $m_i$ ненулевые.

Докажем существование нужной полугруппы индукцией по числу $j$ знаков ''минус'' в линейной комбинации (\ref{b}).

При $j=1$ имеем
$$
b=m_1r_1+\ldots+m_{\tau-1}r_{\tau-1}-m_{\tau}r_{\tau}.
$$
Найдутся положительные рациональные числа $p_1/q_1,\ldots,p_{\tau-1}/q_{\tau-1}$ такие, что
$$
p_1/q_1+\ldots+p_{\tau-1}/q_{\tau-1}=1,
$$
\begin{equation}\label{rational}
\frac{p_i}{q_i}<\frac{m_i{\rm Re}\,r_i}{m_{\tau}{\rm Re}\,r_{\tau}}, \qquad i=1,\ldots,\tau-1.
\end{equation}
Действительно, пересечение $(\tau-1)$-мерного параллелепипеда
$$
\Pi=\{(x_1,\ldots,x_{\tau-1})\in{\mathbb R}^{\tau-1}\mid 0<x_i<m_i{\rm Re}\,r_i/m_{\tau}{\rm Re}\,r_{\tau}\}
$$
с гиперплоскостью
$$
\pi=\{x_1+\ldots+x_{\tau-1}=1\}
$$
является непустым открытым подмножеством в $\pi$, поскольку
$$
\sum_{i=1}^{\tau-1}\frac{m_i{\rm Re}\,r_i}{m_{\tau}{\rm Re}\,r_{\tau}}>1
$$
в силу условия ${\rm Re}\,b>0$. Поэтому пересечение $\Pi\cap\pi$ содержит некоторую точку, все координаты которой рациональны. Теперь образующие $\rho_1,\ldots,\rho_{\tau}$ полугруппы $\Gamma$, линейно независимые над $\mathbb Z$, могут быть определены следующим образом:
\begin{eqnarray*}
\rho_1 & = & r_1-\frac{p_1}{q_1m_1}m_{\tau}r_{\tau}, \quad\ldots,\quad \rho_{\tau-1}=r_{\tau-1}-\frac{p_{\tau-1}}{q_{\tau-1}m_{\tau-1}}m_{\tau}r_{\tau},\\
\rho_{\tau} & = & \frac1{(q_1m_1)\ldots(q_{\tau-1}m_{\tau-1})}r_{\tau}.
\end{eqnarray*}
При этом ${\rm Re}\,\rho_i>0$ согласно неравенствам (\ref{rational}) и $\rm G'\subset\Gamma$, так как $b=m_1\rho_1+\ldots+m_{\tau-1}\rho_{\tau-1}$ и
\begin{eqnarray*}
r_1 & = & \rho_1+\frac{p_1}{q_1m_1}m_{\tau}r_{\tau}=\rho_1+n_1\rho_{\tau}, \quad n_1\in{\mathbb N},\\
\ldots & & \ldots \\
r_{\tau-1} & = & \rho_{\tau-1}+\frac{p_{\tau-1}}{q_{\tau-1}m_{\tau-1}}m_{\tau}r_{\tau}=\rho_{\tau-1}+n_{\tau-1}\rho_{\tau}, \quad n_{\tau-1}\in{\mathbb N},\\
r_{\tau} & = & (q_1m_1)\ldots(q_{\tau-1}m_{\tau-1})\rho_{\tau}.
\end{eqnarray*}

При $j>1$ представим число $b$ в виде $b=b'-m_{\tau}r_{\tau}$, где
$$
b'=m_1r_1+\ldots+m_{\tau-j}r_{\tau-j}-m_{\tau-j+1}r_{\tau-j+1}-\ldots-m_{\tau-1}r_{\tau-1}
$$
содержит $j-1$ знак ''минус'' в соответствующем разложении, аналогичном (\ref{b}). Таким образом, можем применить индуктивное предположение к числу $b'$ и представить его в виде
$$
b'=m'_1\rho_1+\ldots+m'_{\tau-1}\rho_{\tau-1}, \qquad m'_i\in{\mathbb Z}_+,
$$
где числа $\rho_1,\ldots,\rho_{\tau-1}$ линейно независимы над $\mathbb Z$ и выражаются линейными комбинациями чисел $r_1,\ldots,r_{\tau-1}$ с рациональными коэффициентами (сами $r_1,\ldots,r_{\tau-1}$, в свою очередь, представляются линейными комбинациями чисел $\rho_1,\ldots,\rho_{\tau-1}$ с целыми неотрицательными коэффициентами). Следовательно, $\rho_1,\ldots,\rho_{\tau-1}$ и $r_{\tau}$ линейно независимы над $\mathbb Z$, и в силу соотношения $b=m'_1\rho_1+\ldots+m'_{\tau-1}\rho_{\tau-1}-m_{\tau}r_{\tau}$ можно завершить доказательство как при $j=1$. {\hfill $\Box$}
\medskip

Доказанные леммы 4, 5 позволяют нам доказать основное утверждение данного параграфа --- теорему о структуре множества показателей степени обобщённого степенного ряда, формально удовлетворяющего одному из функциональных уравнений (\ref{diff}), (\ref{qdiff}), (\ref{Mahler}).
\medskip

{\bf Теорема 4.} {\it Пусть обобщённый степенной ряд \eqref{powgen} удовлетворяет
\begin{itemize}
\item[--] уравнению \eqref{diff} или \eqref{qdiff} и условию леммы $2$ или
\item[--] уравнению \eqref{Mahler} и условию леммы $3$.
\end{itemize}
Тогда $\varphi\in{\mathbb C}[[x^\Gamma]]\subset\cal G$, где $\Gamma$ --- конечно порождённая аддитивная полугруппа, образующие которой линейно независимы над $\mathbb Z$ и имеют положительные вещественные части.}
\medskip

{Д о к а з а т е л ь с т в о.} Рассмотрим сначала более подробно случай дифференциального уравнения.

Согласно леммам 2 i), 3, 5, ряд \eqref{powgen} представляется в виде $\displaystyle\varphi=\sum_{k=1}^Nc_kx^{\lambda_k}+x^{\lambda_N}\psi$ и ряд $\psi\in\cal G$
формально удовлетворяет уравнению (\ref{auxildiff}),
$$
L(\delta+\lambda_N)\psi=\sum_{{\bf p}\in{\mathbb Z}_+^{n+1}}a_{\bf p}(x)\psi^{p_0}(\delta\psi)^{p_1}\ldots(\delta^n\psi)^{p_n},
$$
где $a_{\bf p}\in{\mathbb C}\{x^\Gamma\}$ и $\Gamma$ --- конечно порождённая аддитивная полугруппа, образующие которой имеют положительные вещественные части и линейно независимы над $\mathbb Z$. Следовательно, имеем соотношение
\begin{eqnarray}\label{auxrel}
\sum_{k=N+1}^{+\infty}L(\lambda_k)c_kx^{\lambda_k-\lambda_N}=\sum_{{\bf p}\in{\mathbb Z}_+^{n+1}}a_{\bf p}(x)\psi^{p_0}(\delta\psi)^{p_1}\ldots(\delta^n\psi)^{p_n}, \\ \psi=\sum_{k=N+1}^{+\infty}c_kx^{\lambda_k-\lambda_N}. \nonumber
\end{eqnarray}
Первое слагаемое в левой части (\ref{auxrel}) --- моном $L(\lambda_{N+1})\,c_{N+1}x^{\lambda_{N+1}-\lambda_N}$, в то время как первое слагаемое в правой части (\ref{auxrel}) --- слагаемое, с которого начинается ряд $a_{\bf 0}(x)\in{\mathbb C}\{x^\Gamma\}$. Число $N$ можно выбрать столь большим, что $L(\lambda_k)\ne0$ при всех $k\geqslant N+1$. Поэтому $\lambda_{N+1}-\lambda_N\in\Gamma$.

Далее, предположим, что $\lambda_k-\lambda_N\in\Gamma$ при всех $k=N+1,\ldots,N+s-1$, и покажем, что тогда $\lambda_{N+s}-\lambda_N\in\Gamma$. Число $\lambda_{N+s}-\lambda_N$ является степенью ненулевого монома $L(\lambda_{N+s})\,c_{N+s}x^{\lambda_{N+s}-\lambda_N}$ --- слагаемого с номером $s$ в левой части соотношения (\ref{auxrel}). Степень слагаемого с тем же номером (та же степень) в правой части (\ref{auxrel}) представляется в виде
$\displaystyle\alpha+\sum_{k=N+1}^{N+s-1} m_k(\lambda_k-\lambda_N)$, $\alpha\in\Gamma$, $m_k\in{\mathbb Z}_+$, то есть по индуктивному предположению является элементом полугруппы $\Gamma$. Тем самым, $\lambda_{N+s}-\lambda_N\in\Gamma$ и, следовательно, $\psi\in{\mathbb C}[[x^\Gamma]]$. Тогда и $\varphi\in{\mathbb C} [[x^\Gamma]]$ (заметим, что $\lambda_1,\ldots,\lambda_N\in\Gamma$ по построениям из лемм 2, 4, 5).

В случае $q$-разностного уравнения и уравнения Малера подобные рассуждения приводят нас к соотношениям, аналогичным (\ref{auxrel}):
\begin{equation}\label{auxrel2}
\sum_{k=N+1}^{+\infty}L(q^{\lambda_k})c_kx^{\lambda_k-\lambda_N}=\sum_{{\bf p}\in{\mathbb Z}_+^{n+1}}a_{\bf p}(x)\psi^{p_0}(\sigma\psi)^{p_1}\ldots (\sigma^n\psi)^{p_n}
\end{equation}
для $q$-разностного уравнения и
$$
\sum_{k=N+1}^{+\infty}c_kx^{\lambda_k-\lambda_N}=\sum_{{\bf p}\in{\mathbb Z}_+^{n+1}}a_{\bf p}(x)\psi^{p_0}(\mu\psi)^{p_1}\ldots(\mu^n\psi)^{p_n}
$$
для уравнения Малера. Из последнего соотношения, так же как и в дифференциальном случае, следует принадлежность формального решения $\varphi\in\cal G$ уравнения Малера подалгебре ${\mathbb C}[[x^\Gamma]]$. Что касается $q$-разностного случая, следует отметить, что, вообще говоря, бесконечное число величин $q^{\lambda_k}$ могли бы совпадать с некоторым корнем $a$ многочлена $L$ (то есть, каким бы большим ни было $N$, тогда нашлось бы $k>N$ такое, что $L(q^{\lambda_k})=0$; см. пример ниже, в котором $L(q^{\lambda_k})=0$ для всех $\lambda_k$). Показатели $\lambda_k-\lambda_N$ соответствующих ненулевых членов $c_kx^{\lambda_k-\lambda_N}$ ряда $\psi$ при этом не обязаны выражаться через предыдущие показатели, как это происходит когда $L(q^{\lambda_k})\ne0$ в силу соотношения (\ref{auxrel2}), аналогично дифференциальному случаю. Тем не менее, в такой ситуации полугруппа $\Gamma$ может быть расширена за счёт добавления конечного числа образующих, достаточных для выражения тех $\lambda_k$, для которых $L(q^{\lambda_k})=0$. Поясним это подробными вычислениями.

Множество решений уравнения $q^\lambda=a$ имеет вид
$$
\lambda=\frac{\arg q\arg a+\ln|q|\ln|a|}{|\ln q|^2}+\frac{\arg q}{|\ln q|^2}\,2\pi n+{\rm i}\,\frac{\ln|q|\arg a-\ln|a|\arg q}{|\ln q|^2}+
{\rm i}\,\frac{\ln|q|}{|\ln q|^2}\,2\pi n, \quad n\in\mathbb Z.
$$
Для каждого корня $a$ многочлена $L$ рассмотрим минимальное значение $n=n(a)$, при котором ${\rm Re}(\lambda-\lambda_N)>0$, и обозначим через $\lambda_a$ соответствующий корень уравнения $q^\lambda=a$. Дополним множество образующих полугруппы $\Gamma$ величинами $\lambda_a-\lambda_N$, $\frac{2\pi}{|\ln q|^2}(\arg q+ {\rm i}\ln|q|)$ и вновь обозначим таким образом полученную полугруппу через $\Gamma$. (Можно считать, что $\arg q>0$, то есть вещественные части добавленных образующих положительны. В случае $\arg q=0$ все корни уравнения $q^\lambda=a$ имеют фиксированную вещественную часть, поэтому среди этих корней может содержаться лишь конечное число показателей $\lambda_k$.) Тогда $\lambda_k-\lambda_N\in\Gamma$, $k\geqslant N+1$. Действительно, если $\lambda_k$ таково, что $L(q^{\lambda_k})=0$, то найдётся корень $a$ многочлена $L$ и $n\in{\mathbb Z}_+$ такие, что
$$
\lambda_k-\lambda_N=\lambda_a-\lambda_N+n\frac{2\pi}{|\ln q|^2}(\arg q+ {\rm i}\ln|q|)\in\Gamma.
$$
Если же $L(q^{\lambda_k})\ne0$, то показатель $\lambda_k-\lambda_N$ ряда $\psi$ выражается через его предыдущие показатели и, как и в дифференциальном случае, применяем индуктивное предположение. Таким образом, $\psi\in{\mathbb C}[[x^\Gamma]]$, а тогда и $\varphi\in{\mathbb C} [[x^\Gamma]]$. {\hfill $\Box$}
\medskip

{\bf Пример 1.} Обобщённый степенной ряд $\displaystyle\varphi=\sum_{k\geqslant1}c_kx^{(1+{\rm i})k}$, $c_k\in\mathbb C$ --- {\it произвольные} постоянные, удовлетворяет $q$-разностному уравнению
$$
y(qx)=y(x), \qquad q=e^{\pi(1+{\rm i})}=-e^\pi,
$$
поскольку $q^{(1+{\rm i})k}=e^{\pi(1+{\rm i})^2k}=e^{2\pi{\rm i}k}=1$. При этом $L(x)=x-1$ и $L(q^{(1+{\rm i})k})=0$ при всех $k$. Отметим, что ряд $\varphi$ содержит {\it бесконечное} число произвольных параметров $c_k$, чего не бывает в случае дифференциального уравнения или уравнения Малера.
\medskip

{\bf Следствие 1.} {\it Пусть обобщённый степенной ряд \eqref{powgen} с {\rm вещественными} показателями степени, $\lambda_k\in\mathbb R$, удовлетворяет уравнению \eqref{diff}, \eqref{qdiff} или \eqref{Mahler}.\footnote{Выполнение условий лемм 2 или 3 здесь не требуется.} Тогда $\varphi\in{\mathbb C}[[x^\Gamma]]\subset\cal G$, где $\Gamma$ --- конечно порождённая аддитивная полугруппа, образующие которой положительны и линейно независимы над $\mathbb Z$.}
\medskip

{Д о к а з а т е л ь с т в о.} Применим индукцию\footnote{Изложенный здесь способ применения индукции заимствован из доказательства Б.\,Мальгранжа \cite{Ma} обобщения теоремы Майе для обыкновенных дифференциальных уравнений.} по порядку уравнения $n$. Случай $n=0$ аналитического уравнения $F(x,y)=0$ описан во введении: обобщённый степенной ряд, ему удовлетворяющий, является не более чем рядом Пюизо \eqref{pow}, и утверждение следствия в этом случае доказано.

Пусть теперь $n\geqslant1$ и $F(x,\varphi,\Delta\varphi,\ldots,\Delta^n\varphi)=0$, где $\Delta$ --- один из операторов $\delta,\sigma,\mu$, в зависимости от того, какому уравнению удовлетворяет ряд $\varphi$. Если $F'_{y_0}(x,\varphi,\Delta\varphi,\ldots,\Delta^n\varphi)\ne0$, то в силу вещественности показателей $\lambda_k$, очевидно, условие теоремы 4 выполнено и поэтому утверждение следствия доказано. Если же $F'_{y_0}(x,\varphi,\Delta\varphi,\ldots,\Delta^n\varphi)=0$, то возможны два случая.

{\sc Случай а}):
$$
F(x,\varphi,\Delta\varphi,\ldots,\Delta^n\varphi)=F'_{y_0}(x,\varphi,\Delta\varphi,\ldots,\Delta^n\varphi)=\ldots=F^{(k)}_{y_0}(x,\varphi,\Delta\varphi,\ldots,\Delta^n\varphi)=0
$$
и $F^{(k+1)}_{y_0}(x,\varphi,\Delta\varphi,\ldots,\Delta^n\varphi)\ne0$, при некотором $k\geqslant1$. В таком случае утверждение вновь следует из теоремы 4, применённой к уравнению $F^{(k)}_{y_0}(x,y,\Delta y,\ldots,\Delta^ny)=0$.

{\sc Случай б}): $F^{(k)}_{y_0}(x,\varphi,\Delta\varphi,\ldots,\Delta^n\varphi)=0$ при всех $k\geqslant0$. Найдётся $\lambda\in\mathbb C$ такое, что $F(x,\lambda x,y_1,$ $\ldots,y_n)\not\equiv0$ (в противном случае $F(x,y_0,y_1,\ldots,y_n)\equiv0$). В то же время, в рассматриваемом случае (формальная) формула Тейлора влечёт равенство
$$
F(x,\lambda x,\Delta\varphi,\ldots,\Delta^n\varphi)=F(x,\varphi+(\lambda x-\varphi),\Delta\varphi,\ldots,\Delta^n\varphi)=
$$
$$
=F(x,\varphi,\Delta\varphi,\ldots,\Delta^n\varphi)+F'_{y_0}(x,\varphi,\Delta\varphi,\ldots,\Delta^n\varphi)(\lambda x-\varphi)+\ldots=0,
$$
которое означает, что обобщённый степенной ряд $\Delta\varphi$ удовлетворяет аналогичному уравнению порядка $n-1$. Следовательно, по индуктивному предположению $\Delta\varphi\in{\mathbb C}[[x^\Gamma]]$, а тогда и $\varphi\in{\mathbb C}[[x^\Gamma]]$ (возможно, при другой $\Gamma$ в случае $\Delta=\mu$). {\hfill $\Box$}
\medskip

{\bf Замечание 5.} Утверждение следствия 1, возможно, без обсуждения линейной независимости над $\mathbb Z$ образующих полугруппы $\Gamma$, для {\it алгебраического} дифференциального уравнения (\ref{diff}) (когда $F$ --- полином) было доказано в \cite{GS}, а для аналитического $q$-разностного уравнения (\ref{qdiff}) --- в \cite{CA}. В данных работах при этом показано, что условие $\lim_{k\rightarrow\infty}\lambda_k=+\infty$ изначально можно не требовать, --- оно является внутренним свойством вещественных показателей ряда (\ref{powgen}), удовлетворяющего дифференциальному или $q$-разностному уравнению. Для уравнений Малера это свойство показателей формального решения, вообще говоря, может не выполняться, и утверждение следствия 1 для таких рядов (не являющихся обобщёнными степенными в смысле нашего определения) будет не верным. Например, линейное уравнение Малера
$$
y(x^\ell)-x^{\ell-1}y(x)=x^{\ell-1}
$$
первого порядка обладает формальным решением $y=\varphi$ в виде ряда
$$
\varphi=\sum\limits_{k=1}^{+\infty}x^{1-1/\ell^k}
$$
с рациональными показателями степени, который не является при этом рядом Пюизо (заметим, что $\lim_{k\rightarrow\infty}(1-1/\ell^k)=1$).

\subsection{Представление элементов алгебры ${\mathbb C} [[x^\Gamma]]$ степенными рядами (Тейлора) нескольких переменных}

Линейная независимость над $\mathbb Z$ образующих аддитивной полугруппы $\Gamma$ (обозначим их $\rho_1,\ldots,\rho_\tau$) позволяет корректно определить биективное линейное отображение
$\iota:{\mathbb C}[[x^{\Gamma}]]\rightarrow{\mathbb C}[[x_1,\ldots,x_\tau]]_*$ из алгебры ${\mathbb C}[[x^{\Gamma}]]$ в алгебру формальных рядов Тейлора $\tau$ переменных без свободного слагаемого,
$$
\iota: \sum_{{\bf m}\in{\mathbb Z}_+^\tau\setminus\{0\}}c_{\bf m}x^{m_1\rho_1+\ldots+m_\tau\rho_\tau}\mapsto
\sum_{{\bf m}\in{\mathbb Z}_+^\tau\setminus\{0\}}c_{\bf m}x_1^{m_1}\ldots x_\tau^{m_\tau}, \qquad {\bf m}=(m_1,\ldots,m_\tau).
$$
Поскольку также
$$
\iota(\eta_1\eta_2)=\iota(\eta_1)\iota(\eta_2) \qquad\forall\eta_1,\eta_2\in{\mathbb C}[[z^{\Gamma}]],
$$
данное отображение является изоморфизмом алгебр.

Как мы отметили в замечании 3, для всякого ряда вида
\begin{equation}\label{genhol}
F(x,y_0,y_1,\ldots,y_n)=\sum_{{\bf p}\in{\mathbb Z}_+^{n+1}}b_{\bf p}(x)y_0^{p_0}y_1^{p_1}\ldots y_n^{p_n}, \qquad b_{\bf p}\in{\mathbb C}[[x^\Gamma]],
\end{equation}
результат подстановки обобщённых степенных рядов $\varphi_0,\varphi_1,\ldots,\varphi_n\in{\mathbb C}[[x^\Gamma]]$ вместо переменных $y_0,y_1,\ldots,y_n$ в $F$ также является элементом алгебры
${\mathbb C}[[x^\Gamma]]$. Более того, имеет место следующее важное свойство, которое мы будем использовать в дальнейшем.
\medskip

{\bf Лемма 6.} {\it Для всякого ряда вида \eqref{genhol} и обобщённых степенных рядов $\varphi_0,\varphi_1,\ldots,\varphi_n\in{\mathbb C}[[x^\Gamma]]$ имеет место равенство
\begin{equation}\label{iotagen}
\iota\bigl(F(x,\varphi_0,\varphi_1,\ldots,\varphi_n)\bigr)=\sum_{{\bf p}\in{\mathbb Z}_+^{n+1}}\iota(b_{\bf p})\iota(\varphi_0)^{p_0}\iota(\varphi_1)^{p_1}\ldots\iota(\varphi_n)^{p_n}.
\end{equation}
}

{Д о к а з а т е л ь с т в о.} Достаточно доказать, что для каждого $N\in\mathbb N$ частичная сумма ряда Тейлора в левой части соотношения (\ref{iotagen}), соответствующая мультииндексам, модули которых не превосходят $N$, совпадает с соответствующей частичной суммой ряда Тейлора в правой части (\ref{iotagen}). Обе частичные суммы определяются конечным числом слагаемых ряда (\ref{genhol}), поскольку для его мультииндексов $\bf p$ с достаточно большими модулями $|{\bf p}|$ ряды Тейлора $\iota(b_{\bf p}\,\varphi_0^{p_0}\varphi_1^{p_1}\ldots\varphi_n^{p_n})=\iota(b_{\bf p})\iota(\varphi_0)^{p_0} \iota(\varphi_1)^{p_1}\ldots\iota(\varphi_n)^{p_n}$ содержат только мономы $x_1^{m_1}\ldots x_\tau^{m_\tau}$, для которых $|{\bf m}|>N$. Но для всякой конечной суммы $\displaystyle\sum_{\bf p}
b_{\bf p}(x)\varphi_0^{p_0}\varphi_1^{p_1}\ldots\varphi_n^{p_n}$ имеет место равенство
$$
\iota\Bigl(\sum_{\bf p}b_{\bf p}\,\varphi_0^{p_0}\varphi_1^{p_1}\ldots\varphi_n^{p_n}\Bigr)=\sum_{\bf p}\iota(b_{\bf p})\iota(\varphi_0)^{p_0}\iota(\varphi_1)^{p_1}\ldots\iota(\varphi_n)^{p_n},
$$
в силу того, что $\iota$ --- изоморфизм алгебр. {\hfill $\Box$}

\section{Условие сходимости обобщённого степенного ряда, удовлетворяющего дифференциальному уравнению}

В этом параграфе мы доказываем следующую теорему\footnote{Данная теорема доказана в \cite{GG} в случае алгебраического дифференциального уравнения (\ref{diff}), где также можно найти пример её применения к третьему уравнению Пенлеве. Здесь, помимо распространения теоремы на аналитический случай, мы также постарались сделать некоторые методические улучшения по изложению её доказательства.} о достаточном условии сходимости обобщённого степенного ряда, удовлетворяющего дифференциальному уравнению.
\medskip

{\bf Теорема 5.} {\it Пусть обобщённый степенной ряд \eqref{powgen} удовлетворяет уравнению \eqref{diff} и
$$
F'_{y_j}(x,\varphi,\delta\varphi,\ldots,\delta^n\varphi)=A_jx^{\nu}+o(x^\nu)\in{\cal G},
$$
где число $\nu$ --- одинаковое для всех $j=0,1,\ldots,n$ и при этом $A_n\ne0$. Тогда ряд $\varphi$ равномерно сходится во всяком секторе $S$ с вершиной в нуле достаточно малого радиуса и раствора, меньшего $2\pi$.}
\medskip

{Д о к а з а т е л ь с т в о.} По теореме 4 ряд $\varphi$ является элементом подалгебры ${\mathbb C}[[x^\Gamma]]\subset\cal G$, где $\Gamma$ --- конечно порождённая аддитивная полугруппа, образующие которой --- обозначим их, как и ранее, $\rho_1,\ldots,\rho_\tau$ --- имеют положительные вещественные части и линейно независимы над $\mathbb Z$.
При этом, согласно лемме 2 (см. также доказательство теоремы 4), $\varphi$ может быть представлен в виде $\displaystyle\varphi=\sum_{k=1}^Nc_kx^{\lambda_k}+x^{\lambda_N}\psi$, и ряд
$\psi\in{\mathbb C}[[x^\Gamma]]$ формально удовлетворяет соотношению
\begin{equation}\label{difforpsi}
L(\delta+\lambda_N)\psi=\sum_{{\bf p}\in{\mathbb Z}_+^{n+1}}a_{\bf p}(x)\psi^{p_0}(\delta\psi)^{p_1}\ldots(\delta^n\psi)^{p_n}, \qquad
a_{\bf p}\in{\mathbb C}\{x^\Gamma\}.
\end{equation}
Поэтому утверждение теоремы о сходимости достаточно доказать для ряда $\psi$, который мы представим в виде
$$
\psi=\sum_{{\bf m}\in{\mathbb Z}_+^\tau\setminus\{0\}}c_{\bf m}\,x^{m_1\rho_1+\ldots+m_\tau\rho_\tau}, \qquad c_{\bf m}\in\mathbb C.
$$

Докажем сначала, что степенной ряд (Тейлора) $\tilde\psi=\iota(\psi)\in{\mathbb C}[[x_1,\ldots,x_\tau]]_*$, соответствующий обобщённому степенному ряду $\psi\in{\mathbb C}[[x^\Gamma]]$ в силу изоморфизма $\iota$, имеет непустой полидиск сходимости, а из этого уже несложно будет вывести утверждение о сходимости самого ряда $\psi$.
\smallskip

{\sl 1. Соотношение для ряда $\tilde\psi$.} Дифференцирование $\delta=x(d/dx)$, действующее на ${\mathbb C}[[x^{\Gamma}]]$ по правилу
$$
\delta: \sum_{|{\bf m}|>0}c_{\bf m}\,x^{m_1\rho_1+\ldots+m_\tau\rho_\tau}\mapsto\sum_{|{\bf m}|>0}(m_1\rho_1+\ldots+m_\tau\rho_\tau)c_{\bf m}\,x^{m_1\rho_1+\ldots+m_\tau\rho_\tau},
$$
с помощью изоморфизма $\iota: {\mathbb C}[[x^{\Gamma}]]\rightarrow{\mathbb C}[[x_1,\ldots,x_\tau]]_*$ естественным образом переносится на ${\mathbb C}[[x_1,\ldots,x_\tau]]_*$,
$$
\tilde\delta: \sum_{|{\bf m}|>0}c_{\bf m}\,x_1^{m_1}\ldots x_\tau^{m_\tau}\mapsto\sum_{|{\bf m}|>0}(m_1\rho_1+\ldots+m_\tau\rho_\tau)c_{\bf m}\,x_1^{m_1}\ldots x_\tau^{m_\tau},
$$
так что $\iota\circ\delta=\tilde\delta\circ\iota$. Поэтому, с учётом леммы 6, применив отображение $\iota$ к обеим частям соотношения (\ref{difforpsi}), получим следующее равенство для $\tilde\psi$:
$$
L(\tilde\delta+\lambda_N)\tilde\psi=\sum_{{\bf p}\in{\mathbb Z}_+^{n+1}}\iota(a_{\bf p})\tilde\psi^{p_0}(\tilde\delta\tilde\psi)^{p_1}\ldots(\tilde\delta^n\tilde\psi)^{p_n}, \qquad
\iota(a_{\bf p})\in{\mathbb C}\{x_1,\ldots,x_\tau\}_*\, ,
$$
которое запишем в виде
\begin{equation}\label{diffortpsi}
L(\tilde\delta+\lambda_N)\tilde\psi=\sum_{{\bf k}\in{\mathbb Z}_+^\tau\setminus\{0\},\,{\bf p}\in{\mathbb Z}_+^{n+1}}A_{\bf k, p}\,x_1^{k_1}\ldots x_\tau^{k_\tau} \tilde\psi^{p_0}(\tilde\delta\tilde\psi)^{p_1}\ldots(\tilde\delta^n\tilde\psi)^{p_n}.
\end{equation}

{\sl 2. Мажорантное уравнение для ряда $\tilde\psi$.} Поскольку по условию теоремы старший коэффициент $A_n$ многочлена $L(z)=A_nz^n+\ldots+A_1z+A_0$ отличен от нуля и, напомним, число $N$ изначально выбрано столь большим, что $L(\lambda_N+\rho)\ne0$ при всех $\rho\in\mathbb C$, имеющих положительную вещественную часть, то найдётся такое $\alpha_1>0$, что
$$
|L(\lambda_N+\rho)|\geqslant\alpha_1\,|\rho|^n
$$
при всех $\rho\in\mathbb C$: ${\rm Re}\,\rho>0$. Также найдётся такое $\alpha_2>0$, что
$$
|L(\lambda_N+\rho)|\geqslant\alpha_2
$$
при всех $\rho\in\mathbb C$: ${\rm Re}\,\rho\geqslant0$, $|\rho|\leqslant1$. Тогда, положив $\alpha=\min(\alpha_1,\alpha_2)$, будем иметь оценки
\begin{equation}\label{est0}
|L(\lambda_N+\rho)|\geqslant\alpha\,|\rho|^j, \qquad j=0,1,\ldots,n,
\end{equation}
при всех $\rho\in\mathbb C$: ${\rm Re}\,\rho>0$.

Для доказательства сходимости ряда $\tilde\psi$ в окрестности $0\in{\mathbb C}^\tau$ рассмотрим уравнение
\begin{equation}\label{majeq}
\alpha\,W=\sum_{{\bf k}\in{\mathbb Z}_+^\tau\setminus\{0\},\,{\bf p}\in{\mathbb Z}_+^{n+1}}|A_{\bf k, p}|\,x_1^{k_1}\ldots x_\tau^{k_\tau} W^{p_0}W^{p_1}\ldots W^{p_n},
\end{equation}
имеющее вид $\alpha\,W-\widetilde M(x_1,\ldots,x_\tau,W)=0$, где $\widetilde M$ --- голоморфная функция в окрестности $0\in{\mathbb C}^{\tau+1}$, $\widetilde M(0,\ldots,0,W)\equiv0$.
По теореме о неявной функции данное уравнение имеет единственное голоморфное решение $W=W(x_1\ldots,x_\tau)$, обращающееся в нуль в начале координат,
\begin{equation}\label{majseries}
W=\sum_{|{\bf m}|>0}C_{\bf m}\,x_1^{m_1}\ldots x_\tau^{m_\tau}.
\end{equation}
Наша дальнейшая задача --- доказать, что этот (сходящийся) степенной ряд является мажорантным для $\tilde\psi$, то есть
$$
C_{\bf m}\geqslant0, \qquad |c_{\bf m}|\leqslant C_{\bf m} \quad \forall {\bf m}\in{\mathbb Z}_+^\tau\setminus\{0\},
$$
откуда и будет следовать сходимость $\tilde\psi$ в окрестности $0\in{\mathbb C}^\tau$.
\smallskip

{\sl 3. Рекуррентные соотношения для коэффициентов $c_{\bf m}$.} Используем соотношение (\ref{diffortpsi}), чтобы выразить каждый коэффициент $c_{\bf m}$ ряда $\tilde\psi$ через его предыдущие коэффициенты. Коэффициент при мономе $x_1^{m_1}\ldots x_\tau^{m_\tau}$ ряда в левой части данного соотношения равен
$$
L(\lambda_N+m_1\rho_1+\ldots+m_\tau\rho_\tau)c_{\bf m}.
$$
Найдём соответствующий коэффициент ряда в правой части соотношения (\ref{diffortpsi}). Для этого нужно найти коэффициент при мономе $x_1^{m_1-k_1}\ldots x_\tau^{m_\tau-k_\tau}$ в произведении
$\tilde\psi^{p_0}(\tilde\delta\tilde\psi)^{p_1}\ldots(\tilde\delta^n\tilde\psi)^{p_n}$. Согласно правилу произведения рядов такой коэффициент равен
$$
\sum_{{\bf l}^{(0)}+{\bf l}^{(1)}+\ldots+{\bf l}^{(n)}={\bf m}-{\bf k}}d_{{\bf l}^{(0)}}d_{{\bf l}^{(1)}}\ldots d_{{\bf l}^{(n)}}
$$
(при $p_j=0$ соответствующий мультииндекс ${\bf l}^{(j)}$ отсутствует под знаком суммы, а множитель $d_{{\bf l}^{(j)}}$ --- в произведении), где каждый множитель $d_{{\bf l}^{(j)}}$ --- коэффициент при мономе $x_1^{l_1^{(j)}}\ldots x_\tau^{l_\tau^{(j)}}$ ряда $(\tilde\delta^j\tilde\psi)^{p_j}$ --- выражается формулой
\begin{equation}\label{d}
d_{{\bf l}^{(j)}}=\sum_{\boldsymbol{\lambda}^{(1)}+\ldots+\boldsymbol{\lambda}^{(p_j)}={\bf l}^{(j)}}
(\lambda_1^{(1)}\rho_1+\ldots+\lambda_\tau^{(1)}\rho_\tau)^jc_{\boldsymbol{\lambda}^{(1)}}\ldots (\lambda_1^{(p_j)}\rho_1+\ldots+\lambda_\tau^{(p_j)}\rho_\tau)^jc_{\boldsymbol{\lambda}^{(p_j)}}, \quad j=0,1,\ldots,n.
\end{equation}
Тем самым, искомое рекуррентное соотношение для $c_{\bf m}$ имеет вид
\begin{equation}\label{crec}
L(\lambda_N+m_1\rho_1+\ldots+m_\tau\rho_\tau)c_{\bf m}=\sum_{0<{\bf k}<{\bf m},\,|{\bf p}|>0}A_{{\bf k},\,{\bf p}}\sum_{{\bf l}^{(0)}+\ldots+{\bf l}^{(n)}=
{\bf m}-{\bf k}}d_{{\bf l}^{(0)}}d_{{\bf l}^{(1)}}\ldots d_{{\bf l}^{(n)}}+A_{{\bf m},\,{\bf 0}},
\end{equation}
и множители $d_{{\bf l}^{(j)}}$ определяются формулой (\ref{d}), то есть каждый коэффициент $c_{\bf m}$ выражается через предыдущие.
\smallskip

{\sl 4. Рекуррентные соотношения для коэффициентов $C_{\bf m}$.} Используя уравнение (\ref{majeq}), которому удовлетворяет ряд (\ref{majseries}), аналогично получаем рекуррентные соотношения для его коэффициентов $C_{\bf m}$:
\begin{equation}\label{Crec}
\alpha\,C_{\bf m}=\sum_{0<{\bf k}<{\bf m},\,|{\bf p}|>0}|A_{{\bf k},\,{\bf p}}|\sum_{{\bf l}^{(0)}+\ldots+{\bf l}^{(n)}=
{\bf m}-{\bf k}}D_{{\bf l}^{(0)}}D_{{\bf l}^{(1)}}\ldots D_{{\bf l}^{(n)}}+|A_{{\bf m},\,{\bf 0}}|,
\end{equation}
где каждый множитель $D_{{\bf l}^{(j)}}$ --- коэффициент при мономе $x_1^{l_1^{(j)}}\ldots x_\tau^{l_\tau^{(j)}}$ ряда $W^{p_j}$, если $p_j\geqslant1$, --- выражается формулой
$$
D_{{\bf l}^{(j)}}=\sum_{\boldsymbol{\lambda}^{(1)}+\ldots+\boldsymbol{\lambda}^{(p_j)}={\bf l}^{(j)}}C_{\boldsymbol{\lambda}^{(1)}}\ldots C_{\boldsymbol{\lambda}^{(p_j)}}, \qquad j=0,1,\ldots,n
$$
(при $p_j=0$ соответствующий мультииндекс ${\bf l}^{(j)}$ отсутствует под знаком внутренней суммы в (\ref{Crec}), а множитель $D_{{\bf l}^{(j)}}$ --- в произведении).

Следовательно, при $|{\bf m}|=1$ имеем $\alpha\,C_{\bf m}=|A_{{\bf m},\,{\bf 0}}|$ и, тем самым, $C_{\bf m}$ --- неотрицательное вещественное число. Тогда и все остальные $C_{\bf m}$, при любом
${\bf m}\in{\mathbb Z}_+^\tau\setminus\{0\}$, в силу рекуррентного соотношения (\ref{Crec}) являются неотрицательными вещественными числами.
\smallskip

{\sl 5. Сходимость ряда $\tilde\psi$ $($оценка $|c_{\bf m}|\leqslant C_{\bf m})$.} Докажем неравенство
\begin{equation}\label{est}
|L(\lambda_N+m_1\rho_1+\ldots+m_\tau\rho_\tau)c_{\bf m}|\leqslant\alpha\,C_{\bf m} \qquad \forall{\bf m}\in{\mathbb Z}_+^\tau\setminus\{0\}
\end{equation}
(из которого, в силу (\ref{est0}) при $j=0$, будет следовать требуемая оценка $|c_{\bf m}|\leqslant C_{\bf m}$). Для этого воспользуемся индукцией по $|{\bf m}|$.

При $|{\bf m}|=1$ согласно (\ref{crec}) и (\ref{Crec}) имеем
$$
L(\lambda_N+m_1\rho_1+\ldots+m_\tau\rho_\tau)c_{\bf m}=A_{{\bf m},\,{\bf 0}}, \qquad \alpha\,C_{\bf m}=|A_{{\bf m},\,{\bf 0}}|,
$$
следовательно, (\ref{est}) выполнено. Предположив, что оценка (\ref{est}) верна при всех $|{\bf m}|<s$, докажем её для $|{\bf m}|=s$ ($s=2,3,\ldots$).

Воспользовавшись соотношением (\ref{crec}), получим
$$
|L(\lambda_N+m_1\rho_1+\ldots+m_\tau\rho_\tau)c_{\bf m}|\leqslant\sum_{0<{\bf k}<{\bf m},\,|{\bf p}|>0}|A_{{\bf k},\,{\bf p}}|\sum_{{\bf l}^{(0)}+\ldots+{\bf l}^{(n)}=
{\bf m}-{\bf k}}|d_{{\bf l}^{(0)}}d_{{\bf l}^{(1)}}\ldots d_{{\bf l}^{(n)}}|+|A_{{\bf m},\,{\bf 0}}|,
$$
при этом, в силу (\ref{d}), (\ref{est0}) и индуктивного предположения,
$$
|d_{{\bf l}^{(j)}}|\leqslant\sum_{\boldsymbol{\lambda}^{(1)}+\ldots+\boldsymbol{\lambda}^{(p_j)}={\bf l}^{(j)}}
|\lambda_1^{(1)}\rho_1+\ldots+\lambda_\tau^{(1)}\rho_\tau|^j\,|c_{\boldsymbol{\lambda}^{(1)}}|\ldots|\lambda_1^{(p_j)}\rho_1+\ldots+\lambda_\tau^{(p_j)}\rho_\tau|^j\,|c_{\boldsymbol{\lambda}^{(p_j)}}| \leqslant
$$
$$
\leqslant\sum_{\boldsymbol{\lambda}^{(1)}+\ldots+\boldsymbol{\lambda}^{(p_j)}={\bf l}^{(j)}}
\frac1\alpha|L(\lambda_N+\lambda_1^{(1)}\rho_1+\ldots+\lambda_\tau^{(1)}\rho_\tau)|\cdot|c_{\boldsymbol{\lambda}^{(1)}}|\ldots \frac1\alpha|L(\lambda_N+\lambda_1^{(p_j)}\rho_1+\ldots+\lambda_\tau^{(p_j)}\rho_\tau)|\cdot|c_{\boldsymbol{\lambda}^{(p_j)}}|\leqslant
$$
$$
\leqslant\sum_{\boldsymbol{\lambda}^{(1)}+\ldots+\boldsymbol{\lambda}^{(p_j)}={\bf l}^{(j)}}C_{\boldsymbol{\lambda}^{(1)}}\ldots C_{\boldsymbol{\lambda}^{(p_j)}}=D_{{\bf l}^{(j)}}, \qquad j=0,1,\ldots,n.
$$
Таким образом,
$$
|L(\lambda_N+m_1\rho_1+\ldots+m_\tau\rho_\tau)c_{\bf m}|\leqslant\sum_{0<{\bf k}<{\bf m},\,|{\bf p}|>0}|A_{{\bf k},\,{\bf p}}|\sum_{{\bf l}^{(0)}+\ldots+
{\bf l}^{(n)}={\bf m}-{\bf k}}D_{{\bf l}^{(0)}}D_{{\bf l}^{(1)}}\ldots D_{{\bf l}^{(n)}}+|A_{{\bf m},\,{\bf 0}}|=\alpha\,C_{\bf m}.
$$

Итак, коэффициенты $c_{\bf m}$ степенного ряда (Тейлора) $\tilde\psi$ ограничены по абсолютной величине коэффициентами $C_{\bf m}$ степенного ряда (\ref{majseries}), сходящегося в окрестности $0\in{\mathbb C}^\tau$, следовательно, ряд $\tilde\psi$ также имеет непустой полидиск сходимости.
\smallskip

Для завершения доказательства утверждения о сходимости ряда $\displaystyle\psi=\sum_{|{\bf m}|>0}c_{\bf m}x^{m_1\rho_1+\ldots+m_\tau\rho_\tau}$ остаётся заметить, что его члены являются однозначными голоморфными функциями во всяком открытом секторе $S\subset\mathbb C$ с вершиной в нуле и раствора, меньшего $2\pi$, и перейти от сходимости $\tilde\psi=\iota(\psi)$ к сходимости $\psi$. Пусть степенной ряд $\tilde\psi$ сходится в полидиске $\{|x_1|\leqslant r,\,\ldots,|x_\tau|\leqslant r\}$. Если $x\in S$ достаточно мало, так что выполнены неравенства
$$
|x^{\rho_1}|=|x|^{{\rm Re}\,\rho_1}e^{-{\rm Im}\,\rho_1\cdot\arg x}<r, \quad\ldots,\quad |x^{\rho_\tau}|=|x|^{{\rm Re}\,\rho_\tau}e^{-{\rm Im}\,\rho_\tau\cdot\arg x}<r
$$
(напомним, что все ${\rm Re}\,\rho_i>0$), то
$$
|c_{\bf m}x^{m_1\rho_1+\ldots+m_\tau\rho_\tau}|=|c_{\bf m}|\cdot|x^{\rho_1}|^{m_1}\ldots|x^{\rho_\tau}|^{m_\tau}<|c_{\bf m}|r^{m_1}\ldots r^{m_\tau}.
$$
Следовательно, обобщённый степенной ряд $\psi$ сходится равномерно в $S$ при достаточно малых $|x|$. {\hfill $\Box$}

\section{Условие сходимости обобщённого степенного ряда, удовлетворяющего $q$-разностному уравнению}

В этом параграфе мы доказываем теорему о достаточном условии сходимости обобщённого степенного ряда, удовлетворяющего $q$-разностному уравнению. Это условие, вообще говоря, может быть сложнее соответствующего условия сходимости в дифференциальном случае и дополнительно содержать некоторые арифметические требования (как и для классических степенных рядов, что можно видеть в теореме 2, случай $|q|=1$).

Пусть обобщённый степенной ряд \eqref{powgen} удовлетворяет уравнению \eqref{qdiff} и
$$
F'_{y_j}(x,\varphi,\sigma\varphi,\ldots,\sigma^n\varphi)=A_jx^{\nu}+o(x^\nu)\in{\cal G},
$$
где число $\nu$ --- одинаковое для всех $j=0,1,\ldots,n$ и при этом не все $A_j$ нулевые.
Тогда по теореме 4 ряд $\varphi$ является элементом подалгебры ${\mathbb C}[[x^\Gamma]]\subset\cal G$, где $\Gamma$ --- конечно порождённая аддитивная полугруппа, образующие которой --- обозначим их, как и ранее, $\rho_1,\ldots,\rho_\tau$ --- линейно независимы над $\mathbb Z$ и имеют положительные вещественные части.
\medskip

{\bf Теорема 6.}\footnote{Данная теорема обобщает результаты нашей совместной работы с Альберто Ластра, изложенные в препринте arXiv:2209.09365.} {\it Пусть $A_n\ne0$ и $A_0\ne0$, а также для каждого корня $z=a$ полинома $(z-1)L(z)$, $L(z)=A_nz^n+\ldots+A_1z+A_0$, выполнено следующее диофантово условие:
\begin{equation}\label{smalldiv}
|(m_1\rho_1+\ldots+m_\tau\rho_\tau)\ln q-\ln a-2\pi m{\rm i}|>c\,|m_1+\ldots+m_\tau|^{-\gamma} \quad \forall\; m_i\in{\mathbb Z}_+,\; m\in{\mathbb Z}
\end{equation}
$($за исключением случая $m_1=\ldots=m_\tau=0)$, где $c$ и $\gamma$ --- некоторые положительные постоянные. Тогда ряд $\varphi$ равномерно сходится во всяком секторе $S$ с вершиной в нуле достаточно малого радиуса и раствора, меньшего $2\pi$.}
\medskip

{\bf Замечание 6.} Диофантово условие (\ref{smalldiv}), в частности, гарантирует, что $L(q^\lambda)\ne0$ ни для какого показателя степени $\lambda=m_1\rho_1+\ldots+m_\tau\rho_\tau$ ряда $\varphi$. Это диофантово условие выполнено в случае образующих $\rho_1,\ldots,\rho_\tau$ общего положения. Для получения конкретных примеров можно воспользоваться, в частности, результатом В.\,Шмидта \cite{Schmidt2}, из которого следует, что (\ref{smalldiv}) выполнено, если
\smallskip

-- ненулевые из вещественных частей величин $\frac1{2\pi\rm i}\rho_1\ln q,\ldots,\frac1{2\pi\rm i}\rho_\tau\ln q,\frac1{2\pi\rm i}\ln a$ суть алгебраические числа и вместе с $1$ линейно независимы над $\mathbb Z$ или
\smallskip

-- ненулевые из мнимых частей величин $\frac1{2\pi\rm i}\rho_1\ln q,\ldots,\frac1{2\pi\rm i}\rho_\tau\ln q,\frac1{2\pi\rm i}\ln a$ суть алгебраические числа, линейно независимые над $\mathbb Z$.
\medskip

{Д о к а з а т е л ь с т в о.} Как и в начале доказательства теоремы 5, представим $\varphi$ в виде $\displaystyle\varphi=\sum_{k=1}^Nc_kx^{\lambda_k}+x^{\lambda_N}\psi$ ($\lambda_N\in\Gamma$) и докажем утверждение о сходимости ряда
$$
\psi=\sum_{{\bf m}\in{\mathbb Z}_+^\tau\setminus\{0\}}c_{\bf m}\,x^{m_1\rho_1+\ldots+m_\tau\rho_\tau}, \qquad c_{\bf m}\in\mathbb C,
$$
согласно лемме 2 формально удовлетворяющего соотношению
\begin{equation}\label{qdifforpsi}
L(q^{\lambda_N}\sigma)\psi=\sum_{{\bf p}\in{\mathbb Z}_+^{n+1}}a_{\bf p}(x)\psi^{p_0}(\sigma\psi)^{p_1}\ldots(\sigma^n\psi)^{p_n}, \qquad
a_{\bf p}\in{\mathbb C}\{x^\Gamma\}.
\end{equation}
Для этого, опять-таки, докажем сначала, что степенной ряд (Тейлора) $\tilde\psi=\iota(\psi)$, представляющий обобщённый степенной ряд $\psi$ в алгебре
${\mathbb C}[[x_1,\ldots,x_\tau]]_*$, имеет непустой полидиск сходимости. Основные шаги в доказательстве этого утверждения повторяют соответствующие шаги из дифференциального случая, однако, поскольку величина $|L(q^{\lambda_N+\rho})|$, в отличие от величины $|L(\lambda_N+\rho)|$, в исключительных случаях может становиться сколь угодно близкой к нулю\footnote{Это может произойти в случае если величины $q^{\lambda_N+\rho}$, $\rho\in\Gamma$, подходят сколь угодно близко к корню полинома $L$.}
при $\rho\in\Gamma$, доказательство оценок для коэффициентов $c_{\bf m}$ степенного ряда $\tilde\psi$ потребует здесь более глубокого анализа с привлечением диофантова условия (\ref{smalldiv}). Мы вынесем эту техническую часть, основанную на идеях работы К.\,Л.\,Зигеля \cite{Siegel}, на которую мы уже ссылались в п.\,2.2.1, в отдельный раздел после доказательства теоремы, чтобы не перегружать деталями изложение его основных этапов.
\smallskip

{\sl 1. Соотношение для ряда $\tilde\psi$.} Оператор $\sigma$, действующий на ${\mathbb C}[[x^{\Gamma}]]$ по правилу
$$
\sigma: \sum_{|{\bf m}|>0}c_{\bf m}\,x^{m_1\rho_1+\ldots+m_\tau\rho_\tau}\mapsto\sum_{|{\bf m}|>0}q^{m_1\rho_1+\ldots+m_\tau\rho_\tau}
c_{\bf m}\,x^{m_1\rho_1+\ldots+m_\tau\rho_\tau},
$$
с помощью изоморфизма $\iota: {\mathbb C}[[x^{\Gamma}]]\rightarrow{\mathbb C}[[x_1,\ldots,x_\tau]]_*$ естественным образом переносится на
${\mathbb C}[[x_1,\ldots,x_\tau]]_*$,
$$
\tilde\sigma: \sum_{|{\bf m}|>0}c_{\bf m}\,x_1^{m_1}\ldots x_\tau^{m_\tau}\mapsto\sum_{|{\bf m}|>0}q^{m_1\rho_1+\ldots+m_\tau\rho_\tau}c_{\bf m}\,x_1^{m_1}\ldots x_\tau^{m_\tau},
$$
так что $\iota\circ\sigma=\tilde\sigma\circ\iota$. Поэтому, с учётом леммы 6, применив отображение $\iota$ к обеим частям соотношения (\ref{qdifforpsi}), получим следующее равенство для $\tilde\psi$:
$$
L(q^{\lambda_N}\tilde\sigma)\tilde\psi=\sum_{{\bf p}\in{\mathbb Z}_+^{n+1}}\iota(a_{\bf p})\tilde\psi^{p_0}(\tilde\sigma\tilde\psi)^{p_1} \ldots(\tilde\sigma^n\tilde\psi)^{p_n}, \qquad
\iota(a_{\bf p})\in{\mathbb C}\{x_1,\ldots,x_\tau\}_*\, ,
$$
которое запишем в виде
\begin{equation}\label{qdiffortpsi}
L(q^{\lambda_N}\tilde\sigma)\tilde\psi=\sum_{{\bf k}\in{\mathbb Z}_+^\tau\setminus\{0\},\,{\bf p}\in{\mathbb Z}_+^{n+1}}A_{\bf k, p}\,x_1^{k_1}\ldots x_\tau^{k_\tau} \tilde\psi^{p_0}(\tilde\sigma\tilde\psi)^{p_1}\ldots(\tilde\sigma^n\tilde\psi)^{p_n}.
\end{equation}

{\sl 2. Мажорантное уравнение для ряда $\tilde\psi$.} Для доказательства сходимости ряда $\tilde\psi$ в окрестности $0\in{\mathbb C}^\tau$ рассмотрим уравнение
\begin{equation}\label{majeq2}
\alpha\,W=\sum_{{\bf k}\in{\mathbb Z}_+^\tau\setminus\{0\},\,{\bf p}\in{\mathbb Z}_+^{n+1}}|A_{\bf k, p}|\,x_1^{k_1}\ldots x_\tau^{k_\tau} W^{p_0}W^{p_1}\ldots W^{p_n},
\end{equation}
имеющее вид $\alpha\,W-\widetilde M(x_1,\ldots,x_\tau,W)=0$, где $\widetilde M$ --- голоморфная функция в окрестности $0\in{\mathbb C}^{\tau+1}$, $\widetilde M(0,\ldots,0,W)\equiv0$, а положительная величина $\alpha$ будет определена далее. По теореме о неявной функции данное уравнение имеет единственное голоморфное решение $W=W(x_1\ldots,x_\tau)$, обращающееся в нуль в начале координат,
\begin{equation}\label{majseries2}
W=\sum_{|{\bf m}|>0}C_{\bf m}\,x_1^{m_1}\ldots x_\tau^{m_\tau}.
\end{equation}
Наша дальнейшая задача --- доказать, что этот (сходящийся) степенной ряд является мажорантным для $\tilde\psi$, но в данном случае под этим мы будем понимать, что
последовательность $\{|c_{\bf m}/C_{\bf m}|\}$, ${\bf m}\in{\mathbb Z}_+^\tau\setminus\{0\}$, имеет не более чем показательный рост. Этого, конечно, будет достаточно для сходимости $\tilde\psi$ в окрестности $0\in{\mathbb C}^\tau$.
\smallskip

{\sl 3. Рекуррентные соотношения для коэффициентов $c_{\bf m}$.} Используем соотношение (\ref{qdiffortpsi}), чтобы выразить каждый коэффициент $c_{\bf m}$ ряда $\tilde\psi$ через его предыдущие коэффициенты. Коэффициент при мономе $x_1^{m_1}\ldots x_\tau^{m_\tau}$ ряда в левой части данного соотношения равен
$$
L\bigl(q^{\lambda_N+m_1\rho_1+\ldots+m_\tau\rho_\tau}\bigr)c_{\bf m}.
$$
Найдём соответствующий коэффициент ряда в правой части соотношения (\ref{qdiffortpsi}). Для этого нужно найти коэффициент при мономе $x_1^{m_1-k_1}\ldots x_\tau^{m_\tau-k_\tau}$ в произведении
$\tilde\psi^{p_0}(\tilde\sigma\tilde\psi)^{p_1}\ldots(\tilde\sigma^n\tilde\psi)^{p_n}$. Согласно правилу произведения рядов такой коэффициент равен
$$
\sum_{{\bf l}^{(0)}+{\bf l}^{(1)}+\ldots+{\bf l}^{(n)}={\bf m}-{\bf k}}d_{{\bf l}^{(0)}}d_{{\bf l}^{(1)}}\ldots d_{{\bf l}^{(n)}}
$$
(при $p_j=0$ соответствующий мультииндекс ${\bf l}^{(j)}$ отсутствует под знаком суммы, а множитель $d_{{\bf l}^{(j)}}$ --- в произведении),
где каждый множитель $d_{{\bf l}^{(j)}}$ --- коэффициент при мономе $x_1^{l_1^{(j)}}\ldots x_\tau^{l_\tau^{(j)}}$ ряда $(\tilde\sigma^j\tilde\psi)^{p_j}$ --- выражается формулой
\begin{equation}\label{d2}
d_{{\bf l}^{(j)}}=\sum_{\boldsymbol{\lambda}^{(1)}+\ldots+\boldsymbol{\lambda}^{(p_j)}={\bf l}^{(j)}}
q^{j(\lambda_1^{(1)}\rho_1+\ldots+\lambda_\tau^{(1)}\rho_\tau)}c_{\boldsymbol{\lambda}^{(1)}}\ldots q^{j(\lambda_1^{(p_j)}\rho_1+\ldots+\lambda_\tau^{(p_j)}\rho_\tau)}c_{\boldsymbol{\lambda}^{(p_j)}}, \quad j=0,1,\ldots,n.
\end{equation}
Тем самым, искомое рекуррентное соотношение для $c_{\bf m}$ имеет вид
\begin{equation}\label{crec2}
L\bigl(q^{\lambda_N+m_1\rho_1+\ldots+m_\tau\rho_\tau}\bigl)c_{\bf m}=\sum_{0<{\bf k}<{\bf m},\,|{\bf p}|>0}A_{{\bf k},\,{\bf p}}
\sum_{{\bf l}^{(0)}+\ldots+{\bf l}^{(n)}={\bf m}-{\bf k}}d_{{\bf l}^{(0)}}d_{{\bf l}^{(1)}}\ldots d_{{\bf l}^{(n)}}+A_{{\bf m},\,{\bf 0}},
\end{equation}
и множители $d_{{\bf l}^{(j)}}$ определяются формулой (\ref{d2}), то есть каждый коэффициент $c_{\bf m}$ выражается через предыдущие.
\smallskip

{\sl 4. Рекуррентные соотношения для коэффициентов $C_{\bf m}$.} Используя уравнение (\ref{majeq2}), которому удовлетворяет ряд (\ref{majseries2}), аналогично получаем рекуррентные соотношения для его коэффициентов $C_{\bf m}$:
\begin{equation}\label{Crec2}
\alpha\,C_{\bf m}=\sum_{0<{\bf k}<{\bf m},\,|{\bf p}|>0}|A_{{\bf k},\,{\bf p}}|\sum_{{\bf l}^{(0)}+\ldots+{\bf l}^{(n)}=
{\bf m}-{\bf k}}D_{{\bf l}^{(0)}}D_{{\bf l}^{(1)}}\ldots D_{{\bf l}^{(n)}}+|A_{{\bf m},\,{\bf 0}}|,
\end{equation}
где каждый множитель $D_{{\bf l}^{(j)}}$ --- коэффициент при мономе $x_1^{l_1^{(j)}}\ldots x_\tau^{l_\tau^{(j)}}$ ряда $W^{p_j}$, если $p_j\geqslant1$, --- выражается формулой
$$
D_{{\bf l}^{(j)}}=\sum_{\boldsymbol{\lambda}^{(1)}+\ldots+\boldsymbol{\lambda}^{(p_j)}={\bf l}^{(j)}}C_{\boldsymbol{\lambda}^{(1)}}\ldots C_{\boldsymbol{\lambda}^{(p_j)}}, \qquad j=0,1,\ldots,n
$$
(при $p_j=0$ соответствующий мультииндекс ${\bf l}^{(j)}$ отсутствует под знаком внутренней суммы в (\ref{Crec2}), а множитель $D_{{\bf l}^{(j)}}$ --- в произведении).

Следовательно, при $|{\bf m}|=1$ имеем $\alpha\,C_{\bf m}=|A_{{\bf m},\,{\bf 0}}|$ и, тем самым, $C_{\bf m}$ --- неотрицательное вещественное число. Тогда и все остальные $C_{\bf m}$, при любом ${\bf m}\in{\mathbb Z}_+^\tau\setminus\{0\}$, в силу рекуррентного соотношения (\ref{Crec2}) являются неотрицательными вещественными числами.
\smallskip

{\sl 5. Сходимость ряда $\tilde\psi$ $($показательный рост последовательности $\{|c_{\bf m}|/C_{\bf m}\})$.} При каждом ${\bf m}\in{\mathbb Z}_+^\tau\setminus\{0\}$ определим величины
$$
\epsilon_{\bf m}=|L\bigl(q^{\lambda_N+m_1\rho_1+\ldots+m_\tau\rho_\tau}\bigr)|^{-1}, \qquad s_{\bf m}=\max(1, |q^{m_1\rho_1+\ldots+m_\tau\rho_\tau}|)
$$
и рассмотрим последовательность положительных чисел $\delta_{\bf m}$, определённую рекуррентно следующим образом. При $|{\bf m}|=1$ положим $\delta_{\bf m}=1$, а при всяком $|{\bf m}|>1$ обозначим через $\mu_{\bf m}$ максимальное среди всех произведений
\begin{eqnarray*}
s^n_{{\bf m}^{(1)}}\delta_{{\bf m}^{(1)}}\ldots s^n_{{\bf m}^{(k)}}\delta_{{\bf m}^{(k)}} & \mbox{таких, что} & {{\bf m}^{(1)}}+\ldots+{{\bf m}^{(k)}}={\bf m},
                                                                                                                  \quad k\geqslant2,\\
                                                                                          &           & {{\bf m}^{(1)}}<{\bf m},\ldots,{{\bf m}^{(k)}}<{\bf m},
\end{eqnarray*}
и определим
$$
\delta_{\bf m}=\epsilon_{\bf m}\mu_{\bf m}.
$$
Выберем также число $\alpha\leqslant1$, участвующее в мажорантном уравнении (\ref{majeq2}), так, чтобы выполнялось условие
$$
\alpha\,\max_{|{\bf m}|=1}\epsilon_{\bf m}\leqslant1.
$$
Тогда для коэффициентов $c_{\bf m}$ ряда $\tilde\psi$ и $C_{\bf m}$ ряда (\ref{majseries2}), удовлетворяющего уравнению (\ref{majeq2}), справедлива оценка
$|c_{\bf m}|\leqslant\delta_{\bf m}C_{\bf m}$. Докажем её индукцией по $|{\bf m}|$.

При $|{\bf m}|=1$ согласно (\ref{crec2}) и (\ref{Crec2}) имеем
$$
\epsilon_{\bf m}^{-1}|c_{\bf m}|=|A_{\bf m,0}|=\alpha\,C_{\bf m}\quad\Longrightarrow\quad
|c_{\bf m}|=\epsilon_{\bf m}\alpha\,C_{\bf m}\leqslant\delta_{\bf m}C_{\bf m}.
$$
Предположив, что оценка $|c_{\bf m}|\leqslant\delta_{\bf m}C_{\bf m}$ верна при всех $|{\bf m}|<s$, докажем её для $|{\bf m}|=s$ ($s=2,3,\ldots$).

Воспользовавшись соотношением (\ref{crec2}), получим
\begin{equation}\label{est2}
\epsilon^{-1}_{\bf m}|c_{\bf m}|\leqslant\sum_{0<{\bf k}<{\bf m},\,|{\bf p}|>0}|A_{{\bf k},\,{\bf p}}|\sum_{{\bf l}^{(0)}+\ldots+{\bf l}^{(n)}=
{\bf m}-{\bf k}}|d_{{\bf l}^{(0)}}d_{{\bf l}^{(1)}}\ldots d_{{\bf l}^{(n)}}|+|A_{{\bf m},\,{\bf 0}}|,
\end{equation}
при этом, в силу (\ref{d2}) и индуктивного предположения,
$$
|d_{{\bf l}^{(j)}}|\leqslant\sum_{\boldsymbol{\lambda}^{(1)}+\ldots+\boldsymbol{\lambda}^{(p_j)}={\bf l}^{(j)}}
s_{\boldsymbol{\lambda}^{(1)}}^j\,|c_{\boldsymbol{\lambda}^{(1)}}|\ldots s_{\boldsymbol{\lambda}^{(p_j)}}^j\,|c_{\boldsymbol{\lambda}^{(p_j)}}|\leqslant
$$
$$
\leqslant\sum_{\boldsymbol{\lambda}^{(1)}+\ldots+\boldsymbol{\lambda}^{(p_j)}={\bf l}^{(j)}} s_{\boldsymbol{\lambda}^{(1)}}^n\delta_{\boldsymbol{\lambda}^{(1)}}\ldots s_{\boldsymbol{\lambda}^{(p_j)}}^n\delta_{\boldsymbol{\lambda}^{(p_j)}}\,
C_{\boldsymbol{\lambda}^{(1)}}\ldots C_{\boldsymbol{\lambda}^{(p_j)}}, \qquad j=0,1,\ldots,n.
$$
Относительно последней суммы имеется следующая альтернатива: если $p_j=1$, то она равна $s_{{\bf l}^{(j)}}^n\,\delta_{{\bf l}^{(j)}}C_{{\bf l}^{(j)}}=
s_{{\bf l}^{(j)}}^n\,\delta_{{\bf l}^{(j)}}D_{{\bf l}^{(j)}}$, в то время как при $p_j\geqslant2$, ввиду определения последовательности чисел $\mu_{\bf m}$, эта сумма не превосходит величины
$$
\mu_{{\bf l}^{(j)}}\sum\limits_{\boldsymbol{\lambda}^{(1)}+\ldots+\boldsymbol{\lambda}^{(p_j)}={\bf l}^{(j)}}C_{\boldsymbol{\lambda}^{(1)}}\ldots
C_{\boldsymbol{\lambda}^{(p_j)}}=\mu_{{\bf l}^{(j)}}D_{{\bf l}^{(j)}}.
$$
Таким образом,
$$
|d_{{\bf l}^{(j)}}|\leqslant\left\{\begin{array}{ll}s_{{\bf l}^{(j)}}^n\,\delta_{{\bf l}^{(j)}}D_{{\bf l}^{(j)}}, & \mbox{если }p_j=1,\\
                                                   \mu_{{\bf l}^{(j)}}D_{{\bf l}^{(j)}}, & \mbox{если }p_j\geqslant2,
                                   \end{array}\right.
$$
и мы приходим к оценке
$$
|d_{{\bf l}^{(0)}}d_{{\bf l}^{(1)}}\ldots d_{{\bf l}^{(n)}}|\leqslant\mu_{\bf m-k}\,D_{{\bf l}^{(0)}}D_{{\bf l}^{(1)}}\ldots D_{{\bf l}^{(n)}}
$$
для каждого набора мультииндексов ${\bf l}^{(0)},{\bf l}^{(1)},\ldots,{\bf l}^{(n)}$ таких, что ${\bf l}^{(0)}+\ldots+{\bf l}^{(n)}={\bf m-k}$ и каждый
${\bf l}^{(j)}<{\bf m-k}$. Из этой оценки, с учётом \eqref{est2}, следует, что
\begin{equation}\label{est3}
\epsilon_{\bf m}^{-1}|c_{\bf m}|\leqslant\sum_{0<{\bf k}<{\bf m},\,|{\bf p}|>0}|A_{{\bf k},\,{\bf p}}|\,\mu_{\bf m-k}
\sum_{{\bf l}^{(0)}+\ldots+{\bf l}^{(n)}={\bf m}-{\bf k}}D_{{\bf l}^{(0)}}D_{{\bf l}^{(1)}}\ldots D_{{\bf l}^{(n)}}+|A_{{\bf m},\,{\bf 0}}|.
\end{equation}
Заметим, что $\mu_{{\bf m}-{\bf k}}\leqslant \mu_{\bf m}$ для всякого ${\bf 0}<{\bf k}<{\bf m}$: при $|{\bf k}|=1$ это следует из цепочки неравенств
$$
\mu_{\bf m-k}\leqslant s_{\bf k}^n\,\mu_{\bf m-k}=s_{\bf k}^n\,\delta_{\bf k}\,\mu_{\bf m-k}\leqslant\mu_{\bf m},
$$
и этого достаточно, чтобы данная оценка выполнялась и при остальных $\bf k$. Поэтому, ввиду \eqref{est3}, \eqref{Crec2} и того, что $\mu_{\bf m}\geqslant1$ (по определению), окончательно получаем
$$
\epsilon_{\bf m}^{-1}|c_{\bf m}|\leqslant\mu_{\bf m}\sum_{0<{\bf k}<{\bf m},\,|{\bf p}|>0}|A_{{\bf k},\,{\bf p}}|
\sum_{{\bf l}^{(0)}+\ldots+{\bf l}^{(n)}={\bf m}-{\bf k}}D_{{\bf l}^{(0)}}D_{{\bf l}^{(1)}}\ldots D_{{\bf l}^{(n)}}+|A_{{\bf m},\,{\bf 0}}|\leqslant
\mu_{\bf m}\,\alpha\,C_{\bf m}\leqslant\mu_{\bf m}C_{\bf m},
$$
откуда следует требуемое неравенство $|c_{\bf m}|\leqslant\delta_{\bf m}C_{\bf m}$.

Ключевой для окончания доказательства сходимости ряда $\tilde\psi$ является следующая лемма.
\medskip

{\bf Лемма 7.} {\it Последовательность $\{\delta_{\bf m}\}$ имеет не более чем показательный рост, если полином $L$, участвующий в её определении, удовлетворяет условиям теоремы 6.}
\medskip

Из этой леммы следует, что отношение $|c_{\bf m}|/C_{\bf m}$ абсолютной величины коэффициентов степенного ряда (Тейлора) $\tilde\psi$ к коэффициентам степенного ряда (\ref{majseries2}), сходящегося в окрестности $0\in{\mathbb C}^\tau$, имеет не более чем показательную скорость роста, следовательно, ряд $\tilde\psi$ также имеет непустой полидиск сходимости.
\smallskip

Завершение доказательства утверждения о сходимости ряда $\displaystyle\psi=\sum_{|{\bf m}|>0}c_{\bf m} x^{m_1\rho_1+\ldots +m_\tau\rho_\tau}$ повторяет соответствующий фрагмент в конце доказательства теоремы 5. {\hfill $\Box$}

\subsection{Доказательство леммы 7}

{\sl 1. Диофантово условие в мультипликативной форме.} Диофантово условие \eqref{smalldiv} влечёт следующие две оценки для каждого корня $z=a$ полинома $L(z)$:
\begin{eqnarray}
|q^{\lambda+m_1\rho_1+\ldots+m_\tau\rho_\tau}-a| & > & 2^{-\gamma}\,|m_1+\ldots+m_\tau|^{-\gamma} \quad \forall\; m_i\in{\mathbb Z}_+, \;\lambda\in\Gamma \label{smalldiv1}, \\
|a\,q^{m_1\rho_1+\ldots+m_\tau\rho_\tau}-q^{\lambda}| & > & 2^{-\gamma}\,|m_1+\ldots+m_\tau|^{-\gamma} \quad \forall\; m_i\in{\mathbb Z}_-, \;\lambda\in\Gamma, \label{smalldiv2}
\end{eqnarray}
а также оценку
\begin{eqnarray}
|q^{m_1\rho_1+\ldots+m_\tau\rho_\tau}-1|\;>\;2^{-\gamma}\,|m_1+\ldots+m_\tau|^{-\gamma} & & \forall\; m_i\in{\mathbb Z}_+ \; \mbox{\rm или} \label{smalldiv3} \\
                                                                                        & & \forall\; m_i\in{\mathbb Z}_- \nonumber
\end{eqnarray}
$($за исключением случая $m_1=\ldots=m_\tau=0)$, где $\gamma$ --- положительная постоянная, возможно, б\'ольшая чем в \eqref{smalldiv}. Докажем эти оценки, воспользовавшись неравенством $|e^z-1|>C|z|$, $C>0$, справедливым при достаточно малых $|z|<\varepsilon$.

1) Допустим, $m_1,\ldots,m_\tau\in{\mathbb Z}_+$ таковы, что
$$
\Bigl|\sum_{i=1}^\tau m_i\rho_i\ln q+\lambda\ln q-\ln a-2\pi m{\rm i}\Bigr|<\varepsilon
$$
при некотором $m\in\mathbb Z$. Тогда
$$
|q^{\lambda+m_1\rho_1+\ldots+m_\tau\rho_\tau}-a|=|a|\cdot|q^{\lambda+m_1\rho_1+\ldots+m_\tau\rho_\tau}a^{-1}-1|=|a|\cdot|e^{\sum_{i=1}^\tau m_i\rho_i\ln q+\lambda\ln q-\ln a-2\pi m{\rm i}}-1|>
$$
$$
>|a|\cdot C\Bigl|\sum_{i=1}^\tau m_i\rho_i\ln q+\lambda\ln q-\ln a-2\pi m{\rm i}\Bigr|>c_1\,|m_1+\ldots+m_\tau|^{-\gamma},
$$
откуда следует \eqref{smalldiv1}. Если же
$$
\Bigl|\sum_{i=1}^\tau m_i\rho_i\ln q+\lambda\ln q-\ln a-2\pi m{\rm i}\Bigr|\geqslant\varepsilon
$$
при всех $m\in\mathbb Z$, то
$$
|q^{\lambda+m_1\rho_1+\ldots+m_\tau\rho_\tau}-a|=|a|\cdot|e^{\sum_{i=1}^\tau m_i\rho_i\ln q+\lambda\ln q-\ln a}-1|>c_2>0,
$$
откуда также следует \eqref{smalldiv1}.

2) Допустим, $m_1,\ldots,m_\tau\in{\mathbb Z}_-$ таковы, что
$$
\Bigl|\sum_{i=1}^\tau m_i\rho_i\ln q-\lambda\ln q+\ln a-2\pi m{\rm i}\Bigr|<\varepsilon
$$
при некотором $m\in\mathbb Z$. Тогда
$$
|a\,q^{m_1\rho_1+\ldots+m_\tau\rho_\tau}-q^\lambda|=|q^\lambda|\cdot|a\,q^{m_1\rho_1+\ldots+m_\tau\rho_\tau}q^{-\lambda}-1|=|q^\lambda|\cdot |e^{\sum_{i=1}^\tau m_i\rho_i\ln q-\lambda\ln q+\ln a-2\pi m{\rm i}}-1|>
$$
$$
>|q^\lambda|\cdot C\Bigl|\sum_{i=1}^\tau m_i\rho_i\ln q-\lambda\ln q+\ln a-2\pi m{\rm i}\Bigr|=|q^\lambda|\cdot C\Bigl|\sum_{i=1}^\tau-m_i\rho_i\ln q+\lambda\ln q -\ln a+2\pi m{\rm i}\Bigr|>
$$
$$
>\tilde c_1\,|-m_1-\ldots-m_\tau|^{-\gamma},
$$
откуда следует \eqref{smalldiv2}. Если же
$$
\Bigl|\sum_{i=1}^\tau m_i\rho_i\ln q-\lambda\ln q+\ln a-2\pi m{\rm i}\Bigr|\geqslant\varepsilon
$$
при всех $m\in\mathbb Z$, то
$$
|a\,q^{m_1\rho_1+\ldots+m_\tau\rho_\tau}-q^\lambda|=|q^\lambda|\cdot|e^{\sum_{i=1}^\tau m_i\rho_i\ln q-\lambda\ln q+\ln a}-1|>\tilde c_2>0,
$$
откуда также следует \eqref{smalldiv2}.

3) Оценка \eqref{smalldiv3} получается абсолютно аналогично \eqref{smalldiv1} и \eqref{smalldiv2}, если использовать условие
$$
|(m_1\rho_1+\ldots+m_\tau\rho_\tau)\ln q-2\pi m{\rm i}|>c\,|m_1+\ldots+m_\tau|^{-\gamma} \quad \forall\; m_i\in{\mathbb Z}_+,\; m\in{\mathbb Z},
$$
которое следует из (\ref{smalldiv}) при $a=1$.
\smallskip

{\sl 2. Вспомогательные оценки-1.} Допустив, не ограничивая общность, что $A_n=1$, разложим $L(z)$ в произведение множителей $z-a$, соответствующих корням $a$ полинома $L$. Определим величины
$$
\varepsilon_{\bf m}(a)=\bigl|q^{\lambda_N+m_1\rho_1+\ldots+m_\tau\rho_\tau}-a\bigr|^{-1}, \qquad {\bf m}\in{\mathbb Z}_+^\tau\setminus\{0\}.
$$
Докажем, что для каждого корня $a$ полинома $L$ и каждой пары мультииндексов $\bf m$, $\bf n$ таких, что ${\bf m}>{\bf n}$, справедливы оценки
\begin{eqnarray}
s_{\bf m}\,\min(\varepsilon_{\bf m}(a),\varepsilon_{\bf n}(a))&<&2^{\gamma+1}(|{\bf m}|-|{\bf n}|)^{\gamma}, \label{auxest1} \\
s_{\bf n}\,\min(\varepsilon_{\bf m}(a),\varepsilon_{\bf n}(a))&<&2^{\gamma+1}(|{\bf m}|-|{\bf n}|)^{\gamma}, \label{auxest1'}
\end{eqnarray}
где $\gamma$ --- положительная постоянная, возможно, б\'ольшая чем в (\ref{smalldiv1})--(\ref{smalldiv3}).

Во-первых, при любых ${\bf m}>{\bf n}$ имеем
$$
a\,\bigl(q^{(m_1-n_1)\rho_1+\ldots+(m_\tau-n_\tau)\rho_\tau}-1\bigr)=\bigl(q^{\lambda_N+m_1\rho_1+\ldots+m_\tau\rho_\tau}-a\bigr)-
$$
\begin{equation}\label{id1}
-\bigl(q^{\lambda_N+n_1\rho_1+\ldots+n_\tau\rho_\tau}-a\bigr)q^{(m_1-n_1)\rho_1+\ldots+(m_\tau-n_\tau)\rho_\tau},
\end{equation}
$$
a\,\bigl(q^{(n_1-m_1)\rho_1+\ldots+(n_\tau-m_\tau)\rho_\tau}-1\bigr)=\bigl(q^{\lambda_N+n_1\rho_1+\ldots+n_\tau\rho_\tau}-a\bigr)-
$$
\begin{equation}\label{id2}
\bigl(q^{\lambda_N+m_1\rho_1+\ldots+m_\tau\rho_\tau}-a\bigr)q^{(n_1-m_1)\rho_1+\ldots+(n_\tau-m_\tau)\rho_\tau}.
\end{equation}
Следовательно, при $|q^{(m_1-n_1)\rho_1+\ldots+(m_\tau-n_\tau)\rho_\tau}|\leqslant1$ равенство \eqref{id1} влечёт оценку
$$
|a|\,|q^{(m_1-n_1)\rho_1+\ldots+(m_\tau-n_\tau)\rho_\tau}-1|\leqslant\varepsilon_{\bf m}^{-1}(a)+\varepsilon_{\bf n}^{-1}(a)\leqslant
\frac2{\min(\varepsilon_{\bf m}(a),\varepsilon_{\bf n}(a))}.
$$
Соответственно, при $|q^{(m_1-n_1)\rho_1+\ldots+(m_\tau-n_\tau)\rho_\tau}|\geqslant1$ из равенства \eqref{id2} следует, что
$$
|a|\,|q^{(n_1-m_1)\rho_1+\ldots+(n_\tau-m_\tau)\rho_\tau}-1|\leqslant\varepsilon_{\bf n}^{-1}(a)+\varepsilon_{\bf m}^{-1}(a)\leqslant
\frac2{\min(\varepsilon_{\bf m}(a),\varepsilon_{\bf n}(a))}.
$$
В любом случае, применяя оценку \eqref{smalldiv3}, для всякого корня $a$ полинома $L$ будем иметь
\begin{equation}\label{est1}
\min(\varepsilon_{\bf m}(a),\varepsilon_{\bf n}(a))<\frac{2^{\gamma+1}}{|a|}\,(|{\bf m}|-|{\bf n}|)^\gamma.
\end{equation}

Во-вторых, при любых ${\bf m}>{\bf n}$ также имеем соотношения
$$
\bigl(q^{\lambda_N+m_1\rho_1+\ldots+m_\tau\rho_\tau}-a\bigr)-\bigl(q^{\lambda_N+n_1\rho_1+\ldots+n_\tau\rho_\tau}-a\bigr)=
q^{\lambda_N+n_1\rho_1+\ldots+n_\tau\rho_\tau}\bigl(q^{(m_1-n_1)\rho_1+\ldots+(m_\tau-n_\tau)\rho_\tau}-1\bigr),
$$
$$
\bigl(q^{\lambda_N+m_1\rho_1+\ldots+m_\tau\rho_\tau}-a\bigr)-\bigl(q^{\lambda_N+n_1\rho_1+\ldots+n_\tau\rho_\tau}-a\bigr)=
q^{\lambda_N+m_1\rho_1+\ldots+m_\tau\rho_\tau}\bigl(1-q^{(n_1-m_1)\rho_1+\ldots+(n_\tau-m_\tau)\rho_\tau}\bigr),
$$
из которых следуют оценки
$$
|q^{\lambda_N}|\,|q^{n_1\rho_1+\ldots+n_\tau\rho_\tau}|\,|q^{(m_1-n_1)\rho_1+\ldots+(m_\tau-n_\tau)\rho_\tau}-1|\leqslant
\varepsilon_{\bf m}^{-1}(a)+\varepsilon_{\bf n}^{-1}(a)\leqslant\frac2{\min(\varepsilon_{\bf m}(a),\varepsilon_{\bf n}(a))},
$$
$$
|q^{\lambda_N}|\,|q^{m_1\rho_1+\ldots+m_\tau\rho_\tau}|\,|q^{(n_1-m_1)\rho_1+\ldots+(n_\tau-m_\tau)\rho_\tau}-1|\leqslant
\varepsilon_{\bf m}^{-1}(a)+\varepsilon_{\bf n}^{-1}(a)\leqslant\frac2{\min(\varepsilon_{\bf m}(a),\varepsilon_{\bf n}(a))}.
$$
Поэтому, с учётом \eqref{smalldiv3}, для всякого корня $a$ полинома $L$ будем иметь
\begin{equation}\label{est1'}
|q^{m_1\rho_1+\ldots+m_\tau\rho_\tau}|\,\min(\varepsilon_{\bf m}(a),\varepsilon_{\bf n}(a))<\frac{2^{\gamma+1}}{|q^{\lambda}|}\,(|{\bf m}|-|{\bf n}|)^\gamma,
\end{equation}
\begin{equation}\label{est1''}
|q^{n_1\rho_1+\ldots+n_\tau\rho_\tau}|\,\min(\varepsilon_{\bf m}(a),\varepsilon_{\bf n}(a))<\frac{2^{\gamma+1}}{|q^{\lambda}|}\,(|{\bf m}|-|{\bf n}|)^\gamma.
\end{equation}

Таким образом, при необходимости увеличивая величину $\gamma$, от неравенств (\ref{est1})--(\ref{est1''}) переходим к нужным оценкам (\ref{auxest1}), (\ref{auxest1'}). Далее под величиной $\gamma$ всегда будем подразумевать постоянную из полученных неравенств (\ref{auxest1}), (\ref{auxest1'}).
\smallskip

{\sl 3. Вспомогательные оценки-2.} Для произвольного корня $a$ полинома $L$ и набора мультииндексов ${\bf m}^{(0)}>{\bf m}^{(1)}>\ldots>{\bf m}^{(r)}>{\bf 0}$, $r\geqslant0$, докажем, что
\begin{equation}\label{auxest2}
\prod\limits_{i=0}^rs_{{\bf m}^{(i)}}\varepsilon_{{\bf m}^{(i)}}(a)<N_1^{r+1}|{\bf m}^{(0)}|^\gamma\prod\limits_{i=1}^r
(|{\bf m}^{(i-1)}|-|{\bf m}^{(i)}|)^\gamma, \qquad N_1=2^{2\gamma+1}.
\end{equation}

Оценка (\ref{auxest2}) верна при $r=0$, поскольку вследствие неравенств (\ref{smalldiv1}), (\ref{smalldiv2}) для ${\bf m}>{\bf 0}$ имеем:
$$
\varepsilon_{\bf m}(a)<2^\gamma\,|{\bf m}|^\gamma, \qquad |q^{m_1\rho_1+\ldots+m_\tau\rho_\tau}|\,\varepsilon_{\bf m}(a)<2^\gamma\,|{\bf m}|^\gamma.
$$
Предположим, $r>0$ и применим индукцию.

Пусть минимум величин $\varepsilon_{{\bf m}^{(i)}}(a)$, $i=0,1,\ldots,r$, достигается при $i=h$. Тогда из неравенств (\ref{auxest1}), (\ref{auxest1'}) следует, что
\begin{equation}\label{min}
s_{{\bf m}^{(h)}}\varepsilon_{{\bf m}^{(h)}}(a)<2^{\gamma+1}\min\bigl\{(|{\bf m}^{(h-1)}|-|{\bf m}^{(h)}|)^\gamma,
(|{\bf m}^{(h)}|-|{\bf m}^{(h+1)}|)^\gamma\bigr\},
\end{equation}
если дополнительно определить $|{\bf m}^{(-1)}|=\infty$, $|{\bf m}^{(r+1)}|=-\infty$. С другой стороны, в предположении, что оценка (\ref{auxest2}) верна при $r-1$ вместо $r$, получим
$$
s_{{\bf m}^{(h)}}^{-1}\varepsilon_{{\bf m}^{(h)}}^{-1}(a)\prod\limits_{i=0}^{r}s_{{\bf m}^{(i)}}\varepsilon_{{\bf m}^{(i)}}(a)<
N_1^r|{\bf m}^{(0)}|^\gamma\frac{(|{\bf m}^{(h-1)}|-|{\bf m}^{(h+1)}|)^\gamma\prod\limits_{i=1}^r(|{\bf m}^{(i-1)}|-|{\bf m}^{(i)}|)^\gamma}
{(|{\bf m}^{(h-1)}|-|{\bf m}^{(h)}|)^\gamma(|{\bf m}^{(h)}|-|{\bf m}^{(h+1)}|)^\gamma}=
$$
$$
=N_1^r|{\bf m}^{(0)}|^\gamma\Bigl(\frac{1}{|{\bf m}^{(h-1)}|-|{\bf m}^{(h)}|}+\frac{1}{|{\bf m}^{(h)}|-|{\bf m}^{(h+1)}|}\Bigr)^\gamma
\prod\limits_{i=1}^r(|{\bf m}^{(i-1)}|-|{\bf m}^{(i)}|)^\gamma\leqslant
$$
$$
\leqslant\frac{2^\gamma N_1^r}{\min\bigl\{(|{\bf m}^{(h-1)}|-|{\bf m}^{(h)}|)^\gamma,(|{\bf m}^{(h)}|-|{\bf m}^{(h+1)}|)^\gamma\bigr\}}\,|{\bf m}^{(0)}|^\gamma \prod\limits_{i=1}^r(|{\bf m}^{(i-1)}|-|{\bf m}^{(i)}|)^\gamma,
$$
и тогда неравенство (\ref{auxest2}) будет следовать из \eqref{min}.

Воспользовавшись равенством
$$
s^n_{{\bf m}^{(i)}}\epsilon_{{\bf m}^{(i)}}=s_{{\bf m}^{(i)}}\varepsilon_{{\bf m}^{(i)}}(a_1)\ldots s_{{\bf m}^{(i)}}\varepsilon_{{\bf m}^{(i)}}(a_n),
$$
где $a_1,\ldots,a_n$ --- корни полинома $L$, как следствие оценки (\ref{auxest2}) получим неравенство
\begin{equation}\label{auxestmain}
\prod\limits_{i=0}^rs^n_{{\bf m}^{(i)}}\epsilon_{{\bf m}^{(i)}}<
\Bigl(N_1^{r+1}|{\bf m}^{(0)}|^\gamma\prod\limits_{i=1}^r(|{\bf m}^{(i-1)}|-|{\bf m}^{(i)}|)^\gamma\Bigr)^n.
\end{equation}

{\sl 4. Окончание доказательства леммы 7.} Докажем, что последовательность величин $\delta_{\bf m}$ имеет не более чем показательный рост, доказав неравенство
\begin{equation}\label{lemma3}
s_{\bf m}^n\delta_{\bf m}\leqslant|{\bf m}|^{-2\gamma n}N_2^{|{\bf m}|-1}Q^{|{\bf m}|},
\end{equation}
где $N_2=8^{\gamma n}N_1^n$, $Q=\max(1,|q^{\rho_1}|^n,\ldots,|q^{\rho_\tau}|^n)$.

Оценка (\ref{lemma3}) верна при $|{\bf m}|=1$, поскольку в этом случае $s_{\bf m}^n\leqslant Q$ и $\delta_{\bf m}=1$. Предположим, что $|{\bf m}|>1$ и применим индукцию.

Последовательность величин $\alpha_{\bf m}=|{\bf m}|^{-2\gamma n}N_2^{|{\bf m}|-1}Q^{|{\bf m}|}$ удовлетворяет неравенствам
$$
\frac{\alpha_{\bf m}\alpha_{\bf m'}}{\alpha_{{\bf m}+{\bf m'}}}=(|{\bf m}|^{-1}+|{\bf m'}|^{-1})^{2\gamma n}N_2^{-1}\leqslant
4^{\gamma n}N_2^{-1}<1, \qquad |{\bf m}|\geqslant1, \; |{\bf m'}|\geqslant1,
$$
следовательно,
\begin{equation}\label{estimate}
s_{{\bf m}^{(1)}}^n\delta_{{\bf m}^{(1)}}\ldots s_{{\bf m}^{(f)}}^n\delta_{{\bf m}^{(f)}}\leqslant|{\bf j}|^{-2\gamma n}N_2^{|{\bf j}|-1}Q^{|{\bf j}|}, \qquad
\mbox{если} \quad {\bf m}^{(1)}+\ldots+{\bf m}^{(f)}={\bf j}, \; |{\bf j}|<|{\bf m}|.
\end{equation}
Согласно определению величин $\delta_{\bf m}$ найдётся разложение
$$
s_{\bf m}^n\delta_{\bf m}=\bigl(s_{\bf m}^n\epsilon_{\bf m}\bigr)s_{{\bf g}^{(1)}}^n\delta_{{\bf g}^{(1)}}\ldots s_{{\bf g}^{(\alpha)}}^n\delta_{{\bf g}^{(\alpha)}}, \qquad \alpha\geqslant2,
$$
где ${\bf g}^{(1)}+\ldots+{\bf g}^{(\alpha)}={\bf m}$ и ${\bf g}^{(1)}<{\bf m},\ldots,{\bf g}^{(\alpha)}<{\bf m}$.
В случае если $|{\bf g}^{(1)}|>|{\bf m}|/2$ (и тогда $|{\bf g}^{(2)}|+\ldots+|{\bf g}^{(\alpha)}|<|{\bf m}|/2$), мы вновь используем подобное разложение, заменив
${\bf m}$ на ${\bf g}^{(1)}$:
$$
s_{{\bf g}^{(1)}}^n\delta_{{\bf g}^{(1)}}=\bigl(s_{{\bf g}^{(1)}}^n\epsilon_{{\bf g}^{(1)}}\bigr)s_{{\bf h}^{(1)}}^n\delta_{{\bf h}^{(1)}}\ldots
s_{{\bf h}^{(\beta)}}^n\delta_{{\bf h}^{(\beta)}}, \qquad \beta\geqslant 2,
$$
где ${\bf h}^{(1)}+\ldots+{\bf h}^{(\beta)}={\bf g}^{(1)}$ и ${\bf h}^{(1)}<{\bf g}^{(1)},\ldots,{\bf h}^{(\beta)}<{\bf g}^{(1)}$.
Если $|{\bf h}^{(1)}|>|{\bf m}|/2$ (и тогда $|{\bf h}^{(2)}|+\ldots+|{\bf h}^{(\beta)}|<|{\bf m}|/2$), раскладываем снова:
$$
s_{{\bf h}^{(1)}}^n\delta_{{\bf h}^{(1)}}=\bigl(s_{{\bf h}^{(1)}}^n\epsilon_{{\bf h}^{(1)}}\bigr)s_{{\bf i}^{(1)}}^n\delta_{{\bf i}^{(1)}}\ldots
s_{{\bf i}^{(\gamma)}}^n\delta_{{\bf i}^{(\gamma)}}, \qquad \gamma\geqslant 2,
$$
где ${\bf i}^{(1)}+\ldots+{\bf i}^{(\gamma)}={\bf h}^{(1)}$ и ${\bf i}^{(1)}<{\bf h}^{(1)},\ldots,{\bf i}^{(\gamma)}<{\bf h}^{(1)}$, и так далее.

Положим теперь
$$
{\bf m}^{(0)}={\bf m},\; {\bf m}^{(1)}={\bf g}^{(1)},\; {\bf m}^{(2)}={\bf h}^{(1)},\;\ldots,\;{\bf m}^{(r)}=\dots,
$$
$$
{\bf m}^{(0)}>{\bf m}^{(1)}>{\bf m}^{(2)}>\ldots>{\bf m}^{(r)},
$$
так что ${\bf m}^{(r)}$ --- последний мультииндекс в этой цепочке, модуль которого превосходит $|{\bf m}|/2$. Тогда получим формулу
\begin{equation}\label{decomp}
s_{\bf m}^n{\delta}_{\bf m}=\bigl(s_{{\bf m}^{(0)}}^n\epsilon_{{\bf m}^{(0)}}\Delta_0\bigr)\bigl(s_{{\bf m}^{(1)}}^n\epsilon_{{\bf m}^{(1)}}\Delta_1\bigr)\ldots
\bigl(s_{{\bf m}^{(r)}}^n\epsilon_{{\bf m}^{(r)}}\Delta_r\bigr),
\end{equation}
где
$$
\begin{array}{lllll}
\Delta_0&=&s_{{\bf g}^{(2)}}^n\delta_{{\bf g}^{(2)}}\ldots s_{{\bf g}^{(\alpha)}}^n\delta_{{\bf g}^{(\alpha)}}, & & {\bf g}^{(2)}+\ldots+{\bf g}^{(\alpha)}=
                                                                                                                    {\bf m}^{(0)}-{\bf m}^{(1)},\\
\Delta_1&=&s_{{\bf h}^{(2)}}^n\delta_{{\bf h}^{(2)}}\ldots s_{{\bf h}^{(\beta)}}^n\delta_{{\bf h}^{(\beta)}},  & & {\bf h}^{(2)}+\ldots+{\bf h}^{(\beta)}=
                                                                                                                   {\bf m}^{(1)}-{\bf m}^{(2)},\\
\hdots& &\hdots\hdots\hdots\hdots& &\hdots\hdots\hdots\hdots\hdots\hdots\hdots\hdots\\
\Delta_r&=&s_{{\bf j}^{(1)}}^n\delta_{{\bf j}^{(1)}}\ldots s_{{\bf j}^{(\omega)}}^n\delta_{{\bf j}^{(\omega)}}, & & {\bf j}^{(1)}+\ldots+{\bf j}^{(\omega)}={\bf m}^{(r)}, \quad \omega\geqslant2,
\end{array}
$$
и модули всех мультииндексов ${\bf g}^{(\cdot)}, {\bf h}^{(\cdot)},\ldots,{\bf j}^{(\cdot)}$ не превосходят $|{\bf m}|/2$.

Воспользовавшись оценкой \eqref{lemma3} для каждого из $\omega$ множителей, составляющих величину $\Delta_r$, а также применив \eqref{estimate} для оценки величин  $\Delta_p$, $p=0,\,1,\ldots,\,r-1$, приходим к неравенствам
$$
\Delta_r\leqslant \left(|{\bf j}^{(1)}|\ldots|{\bf j}^{(\omega)}|\right)^{-2\gamma n}N_2^{|{\bf m}^{(r)}|-\omega}Q^{|{\bf m}^{(r)}|},
$$
$$
\Delta_p\leqslant \left(|{\bf m}^{(p)}|-|{\bf m}^{(p+1)}|\right)^{-2\gamma n}N_2^{|{\bf m}^{(p)}|-|{\bf m}^{(p+1)}|-1}Q^{|{\bf m}^{(p)}|-|{\bf m}^{(p+1)}|}, \qquad p=0,\,1,\ldots,\,r-1,
$$
из которых следует, что
$$
\Delta_0\ldots\Delta_r\leqslant N_2^{|{\bf m}|-r-\omega}Q^{|{\bf m}|}\Bigl(\prod\limits_{i=1}^\omega|{\bf j}^{(i)}|\,\prod\limits_{p=0}^{r-1}
(|{\bf m}^{(p)}|-|{\bf m}^{(p+1)}|)\Bigr)^{-2\gamma n}.
$$
В то же время, в силу (\ref{auxestmain}) имеем
$$
s_{{\bf m}^{(0)}}^n\epsilon_{{\bf m}^{(0)}}\ldots s_{{\bf m}^{(r)}}^n\epsilon_{{\bf m}^{(r)}}\leqslant N_1^{n(r+1)}|{\bf m}|^{\gamma n}\,
\prod\limits_{p=0}^{r-1}(|{\bf m}^{(p)}|-|{\bf m}^{(p+1)}|)^{\gamma n}.
$$
Таким образом, соотношение (\ref{decomp}) приводит нас к оценке
$$
s_{\bf m}^n\delta_{\bf m}\leqslant N_1^{n(r+1)}N_2^{|{\bf m}|-r-\omega}Q^{|{\bf m}|} \Bigl(|{\bf m}|^{-1}\,\prod\limits_{i=1}^\omega|{\bf j}^{(i)}|^2\, \prod\limits_{p=0}^{r-1}(|{\bf m}^{(p)}|-|{\bf m}^{(p+1)}|)\Bigr)^{-\gamma n}=
$$
$$
=N_1^{n(r+1)}N_2^{|{\bf m}|-r-\omega}Q^{|{\bf m}|}\Bigl(|{\bf m}|^{-1}\,\prod\limits_{i=1}^\omega y_i^2\,\prod\limits_{p=0}^{r-1}x_p\Bigr)^{-\gamma n},
$$
где $x_p=|{\bf m}^{(p)}|-|{\bf m}^{(p+1)}|$, $y_i=|{\bf j}^{(i)}|$. Поскольку
$$
\sum\limits_{i=1}^\omega y_i=|{\bf m}^{(r)}|>|{\bf m}|/2, \quad \sum\limits_{p=0}^{r-1}x_p+\sum\limits_{i=1}^\omega y_i=|{\bf m}|, \quad\mbox{и все }\; y_i\leqslant|{\bf m}|/2,
$$
согласно лемме 1 из \cite{Siegel} для произведения в последних скобках справедлива следующая оценка снизу:
$$
\prod\limits_{i=1}^\omega y_i^2\,\prod\limits_{p=0}^{r-1}x_p\geqslant \frac{|{\bf m}|^3}{8^{r+\omega-1}}.
$$
Учитывая её, окончательно получаем нужное неравенство (\ref{lemma3}):
$$
s_{\bf m}^n\delta_{\bf m}\leqslant N_1^{n(r+1)}N_2^{|{\bf m}|-r-\omega}Q^{|{\bf m}|}|{\bf m}|^{-2\gamma n}\;8^{(r+\omega-1)\gamma n}\leqslant |{\bf m}|^{-2\gamma n}
N_2^{|{\bf m}|-1}Q^{|{\bf m}|}\frac{N_1^{n(r+\omega-1)}}{N_2^{r+\omega-1}}\,8^{(r+\omega-1)\gamma n}=
$$
$$
=|{\bf m}|^{-2\gamma n}N_2^{|{\bf m}|-1}Q^{|{\bf m}|}.
$$

\subsection{Ослабление условий теоремы 6 в специальных случаях расположения образующих полугруппы $\Gamma$ на комплексной плоскости}

Отметим, что условия теоремы 6 предполагают общий случай расположения образующих $\rho_1,\ldots,\rho_\tau$ полугруппы $\Gamma$ на комплексной плоскости. Однако, в некоторых специальных случаях такого расположения условия теоремы можно ослабить. Поясним это подробнее.

Главной в доказательстве теоремы 6 является лемма 7 о скорости роста последовательности $\{\delta_{\bf m}\}$, основанная на оценке (\ref{auxestmain}). Именно для получения последней и нужны условия теоремы: условие $A_0\ne0$ (нуль не является корнем полинома $L$) и диофантово условие (\ref{smalldiv}) используются при получении оценок (\ref{smalldiv1})-(\ref{smalldiv3}) и (\ref{auxest1}), (\ref{auxest1'}), из которых затем следует (\ref{auxest2}), а условие $A_n\ne0$ (степень $d:=\deg L$ полинома $L$ совпадает с порядком $n$ уравнения) используется при переходе от (\ref{auxest2}) к (\ref{auxestmain}).

Несложно заметить, что если все $|q^{\rho_i}|\geqslant1$, то оценки (\ref{smalldiv1}), (\ref{smalldiv2}) и (\ref{auxest1}), (\ref{auxest1'}) справедливы и при $A_0=0$ (при $a=0$), а в таком случае справедлива и следующая из них оценка (\ref{auxest2}). С другой стороны, если все $|q^{\rho_i}|\leqslant1$, то оценка (\ref{auxestmain}) следует из (\ref{auxest2}) и при $A_n=0$ (при $d<n$), поскольку в этом случае $s^n_{\bf m}=s^d_{\bf m}=1$ при всех $\bf m$ и поэтому
$$
s^n_{\bf m}\epsilon_{\bf m}=s^d_{\bf m}\epsilon_{\bf m}=s_{\bf m}\varepsilon_{\bf m}(a_1)\ldots s_{\bf m}\varepsilon_{\bf m}(a_d),
$$
где $a_1,\ldots,a_d$ --- корни полинома $L$.

Таким образом, утверждение леммы 7 остаётся справедливым при $A_0=0$, если все $|q^{\rho_i}|\geqslant1$, и при $A_n=0$, если все $|q^{\rho_i}|\leqslant1$.
Кроме того, если все $|q^{\rho_i}|>1$ или все $|q^{\rho_i}|<1$ (строгие неравенства), то диофантово условие (\ref{smalldiv}) выполняется автоматически, поскольку
в этих случаях все ${\rm Re}(\rho_i\ln q)>\kappa>0$ или, соответственно, все ${\rm Re}(\rho_i\ln q)<-\kappa<0$ при некотором $\kappa$, поэтому
$$
|{\rm Re}\bigl((m_1\rho_1+\ldots+m_\tau\rho_\tau)\ln q\bigl)|>\kappa|{\bf m}|
$$
уходит на бесконечность при $|{\bf m}|\rightarrow\infty$.

С учётом изложенных выше наблюдений можем сформулировать дополнение к теореме 6, в котором выделены частные случаи расположения величин $q^{\rho_i}$ относительно единичной окружности 
и сделаны соответствующие ослабления условий теоремы 6.
\medskip

{\bf Теорема 6 bis.} {\it Утверждение теоремы $6$ справедливо в следующих частных случаях:
\smallskip

{\rm а)} $A_0\ne0$, и все $q^{\rho_i}$ расположены {\rm строго внутри} единичной окружности $\{|z|=1\}$;

{\rm б)} $A_n\ne0$, и все $q^{\rho_i}$ расположены {\rm строго вне} единичной окружности;

{\rm в)} все $q^{\rho_i}$ расположены {\rm на} единичной окружности, и условие \eqref{smalldiv} выполнено для корней $z=a$ полинома $(z-1)L(z)$, имеющих такое же расположение;

{\rm г)} $A_0\ne0$, все $q^{\rho_i}$ расположены {\rm внутри} или {\rm на} единичной окружности, и условие \eqref{smalldiv} выполнено для корней $z=a$ полинома $(z-1)L(z)$, имеющих такое же расположение;

{\rm д)} $A_n\ne0$, все $q^{\rho_i}$ расположены {\rm вне} или {\rm на} единичной окружности, и условие \eqref{smalldiv} выполнено для корней $z=a$ полинома $(z-1)L(z)$, имеющих такое же расположение.} {\hfill $\Box$}
\medskip

{\bf Замечание 7.} Расположение числа $q^{\rho_i}$ {\it внутри} единичной окружности эквивалентно расположению самой образующей $\rho_i$ {\it над} прямой
${\cal L}\subset\mathbb C$, проходящей через начало координат и имеющей коэффициент наклона $\ln |q|/{\rm arg}\,q$. Соответственно, расположение $q^{\rho_i}$ {\it вне} единичной окружности эквивалентно расположению образующей $\rho_i$ {\it под} прямой ${\cal L}$.
\medskip

{\bf Замечание 8.} Утверждение теоремы 6 при выполнении условий а) или б) ранее доказано отдельно в работе \cite{GGL}. В случае классического степенного ряда $\varphi\in{\mathbb C}[[x]]$ эти условия превращаются в условие 1) теоремы 2, в то время как условие в) --- в условие 2) теоремы 2 (см. введение).

Отметим также, что оригинальное доказательство теоремы 2 (случай $|q|=1$) представлено в \cite{dv} лишь в виде наброска и опирается на теорему о неявном отображении для банаховых пространств. В качестве основного аргумента в этом наброске выступает ссылка на теорему 6.1 из \cite{Bez2} --- частный случай теоремы 2 ($|q|=1$) для {\it линейного} $q$-разностного уравнения. Однако, последняя теорема формулируется в \cite{Bez2} при диофантовом условии (\ref{dioph}), но остаётся справедливой (доказывается) и при менее ограничительном условии
$$
|q^k-a|\geqslant c_1\,c_2^k, \qquad k\in{\mathbb N},
$$
$c_1$, $c_2>0$. Данное условие, вообще говоря, уже не достаточно для сходимости формального решения {\it нелинейного} $q$-разностного уравнения: например, при $|q^k-1|\geqslant c_1\,c_2^k$ рост функции $\Omega^{-1}(n)$ из замечания 1 может быть экспоненциальным и условие Брюно (\ref{brunocond}) может не выполняться, а тогда, как мы знаем из результатов Ж.-К.\,Йоккоза (см. п.\,2.2.1), формальное решение нелинейного $q$-разностного уравнения
$$
y(qx)=f(y(x)), \qquad f(x)=qx+\ldots\,,
$$
при соответствующем $q$ и определённых $f$ может быть расходящимся. Поэтому, на наш взгляд, среди прочего также важно иметь более подробное доказательство теоремы 2 в случае $|q|=1$.

\subsection{Примеры}

Рассмотрим некоторые примеры $q$-разностных уравнений, в которых мы демонстрируем применение теорем 6 или 6 bis и обсуждаем препятствия к их применению. Как следует из примеров 4, 5, наличие этих препятствий не всегда указывает на расходимость формальных решений, что ожидаемо говорит о том, что условия теорем являются достаточными, но не необходимыми. Подобного рода примеры можно привести и для иллюстрации соответствующих теорем в случае дифференциального уравнения или уравнения Малера (теоремы 5, 7). Мы уделяем внимание $q$-разностному случаю ввиду большего разнообразия эффектов, вызванного наличием дополнительного арифметического условия в формулировке теоремы в этом случае.
\medskip

{\bf Пример 2.} Уравнение
\begin{equation}\label{ex2eq1}
\sigma^2 y-q^{\rm i}\,\sigma y+q^{2\rm i}\,x(1+y^2)=0, \qquad \displaystyle q=e^{-\frac{\sqrt2}2\pi(1-\rm i)},
\end{equation}
обладает формальным решением
\begin{equation}\label{ex2eq2}
c_{0,0}\,x^{\rm i}+\sum\limits_{m_1+m_2>0}c_{m_1,m_2}\,x^{{\rm i}+m_1(1+{\rm i})+m_2(1-{\rm i})},
\end{equation}
где коэффициент $c_{0,0}\neq 0$ --- произвольное комплексное число, а остальные коэффициенты $c_{m_1,m_2}$ однозначно определяются по $c_{0,0}$.

Напрямую к формальному решению \eqref{ex2eq2} мы не можем применить теорему 6 или 6 bis, поскольку эти теоремы сформулированы для обобщённых степенных рядов, показатели степени которых имеют положительные вещественные части. Но поскольку уравнение \eqref{ex2eq1} полиномиально, в нём можно сделать замену переменной $xy=u$, после которой мы получим уравнение
$$
\frac{\sigma^2u}{q^2x}-q^{\rm i}\,\frac{\sigma u}{qx}+q^{2\rm i}\,x\left(1+\frac{u^2}{x^2}\right)=0.
$$
Домножив его на $q^2x$, получим аналитическое уравнение
\begin{equation}\label{ex2eq3}
\sigma^2u-q^{1+\rm i}\,\sigma u+q^{2(1+\rm i)}\,\left(x^2+u^2\right)=0
\end{equation}
с формальным решением нужного вида
$$
\varphi=\sum\limits_{m_1\geqslant1, m_2\geqslant0}c_{m_1-1,m_2}\,x^{m_1(1+{\rm i})+m_2(1-{\rm i})}
$$
и значениями $q^{1+\rm i}=e^{-\sqrt2\pi}<1$ и $q^{1-\rm i}=e^{\sqrt2\pi\rm i}\in\{|z|=1\}$.
Но тем не менее, теорема 6 (6 bis) не может быть применена здесь, поскольку полином $L(z)=z(z-q^{1+\rm i})$ обращается в нуль при $z=0$. Покажем, что ряд $\varphi$ в действительности имеет нулевой радиус сходимости.

Рассмотрим ряд
$$
\varphi_1=\sum_{m=1}^\infty c_{m-1,0}\,x^{m(1+{\rm i})},
$$
образованный частью слагаемых обобщённого степенного ряда $\varphi$. Поскольку $\varphi$ удовлетворяет уравнению (\ref{ex2eq3}), то для $\varphi_1$ справедливо равенство
$$
\sigma^2\varphi_1-q^{1+\rm i}\sigma\varphi_1=-q^{2(1+\rm i)}\varphi_1^2,
$$
из которого следует соотношение
$$
\sum_{m=1}^\infty q^{m(1+{\rm i})}(q^{m(1+{\rm i})}-q^{1+{\rm i}})c_{m-1,0}\,x^{m(1+{\rm i})}=-q^{2(1+\rm i)}\,x^{2(1+\rm i)}\Bigl(\sum_{m=1}^\infty c_{m-1,0}\,x^{(m-1)(1+{\rm i})}\Bigr)^2,
$$
{или, после деления на $q^{2(1+\rm i)}\,x^{1+\rm i}$,}
$$
\sum_{m=1}^\infty q^{m(1+{\rm i})}(q^{m(1+{\rm i})}-1)c_{m,0}\,x^{m(1+{\rm i})}=-x^{1+\rm i}\biggl(c_{0,0}^2+\sum_{m=1}^\infty2c_{0,0}c_{m,0}\,x^{m(1+{\rm i})}+
\Bigl(\sum_{m=1}^\infty c_{m,0}\,x^{m(1+{\rm i})}\Bigr)^2\biggr).
$$
Таким образом, получаем следующие рекуррентные выражения для коэффициентов $c_{m,0}$:
$$
c_{1,0}=\frac{c_{0,0}^2}{e^{-\sqrt2\pi}(1-e^{-\sqrt2\pi})}, \qquad c_{m,0}=\frac{2c_{0,0}c_{m-1,0}+\sum_{j=1}^{m-2}c_{j,0}c_{m-1-j,0}}
{e^{-\sqrt2\pi m}(1-e^{-\sqrt2\pi m})}, \quad m\geqslant2.
$$
Если $c_{0,0}$ --- положительное вещественное число, то такими же будут и остальные коэффициенты $c_{m,0}$. Положив $c_{0,0}\geqslant1/2$, будем иметь
$$
c_{m,0}\geqslant e^{\sqrt2\pi m}c_{m-1,0}\geqslant\ldots\geqslant e^{\sqrt2\pi m(m+1)/2}c_{0,0},
$$
откуда следует расходимость ряда $\varphi_1$ при $x\ne0$ и такая же расходимость ряда $\varphi$.
\medskip

{\bf Пример 3.} Рассмотрим следующий вариант $q$-разностного аналога уравнения Пенлеве III при $a=b=0$, $c=d=1$:
$$
y\,\sigma^2 y-(\sigma y)^2-x^2y^4-x^2=0,
$$
где $q=e^{2{\rm i}\pi\omega}$, $\omega\in\mathbb{R}\setminus\mathbb{Q}$.
Данное уравнение обладает двухпараметрическим семейством формальных решений
$$
\varphi=\sum_{m_1,m_2\in{\mathbb Z}_+}c_{m_1,m_2}\,x^{(m_1+1)r+m_2(1-r)},
$$
с произвольным комплексным коэффициентом $c_{0,0}\ne0$ и $0<{\rm Re}\,r<1$, и остальными коэффициентами $c_{m_1,m_2}$, однозначно определяемыми по $c_{0,0}$ и $r$.
Многочлен $L(z)=c_{0,0}\,(z-q^r)^2$ не обращается в нуль при $z=0$, а числа $q^{r}$ и $q^{1-r}$ расположены либо по разные стороны единичной окружности (если ${\rm Im}\,r\ne0$), либо на самой окружности (если ${\rm Im}\,r=0$). Тем самым, подобрав $r$ и $\omega$ так, чтобы выполнялось диофантово условие (\ref{smalldiv}),
$$
\left|m_1r\omega+m_2(1-r)\omega-m\right|>c\,(m_1+m_2)^{-\gamma},
$$
при некоторых положительных $c$ и $\gamma$, мы получим сходящийся ряд $\varphi$. Например, для этого достаточно, чтобы число $\omega$ было алгебраическим, а число $r$ просто имело ненулевую мнимую часть.
\medskip

{\bf Пример 4.} Уравнение третьего порядка
$$
\sigma^2y-(q^{1+\rm i}+q^{1-\rm i})\sigma y+q^2y+\frac{x^2}{q^6}(\sigma^3y-q^{1+\rm i}\sigma^2y)(\sigma^3y-q^{1-\rm i}\sigma^2y)+(q^{\rm i}-q^{-\rm i})^2x^4=0
$$
обладает решением в виде (конечного, следовательно, сходящегося) обобщённого степенного ряда $\varphi=x^{1+\rm i}+x^{1-\rm i}$. При этом вычисление частных производных $F'_{y_j}$ на этом решении показывает, что
$$
F'_{y_0}(x,\varphi,\sigma\varphi,\sigma^2\varphi,\sigma^3\varphi)=q^2, \qquad F'_{y_1}(x,\varphi,\sigma\varphi,\sigma^2\varphi,\sigma^3\varphi)=-(q^{1+\rm i}+ q^{1-\rm i}),
$$
$$
F'_{y_2}(x,\varphi,\sigma\varphi,\sigma^2\varphi,\sigma^3\varphi)=1+o(x^2)
$$
$$
F'_{y_3}(x,\varphi,\sigma\varphi,\sigma^2\varphi,\sigma^3\varphi)=o(x^2),
$$
поэтому $A_3=0$ и если хотя бы одно из чисел $q^{1\pm\rm i}$ лежит вне единичной окружности, то условие теоремы 6 (6 bis) не выполнено.
\medskip

{\bf Пример 5.} Уравнение первого порядка
$$
(q^a-1)(qx-1)(x-1)y\sigma y+(qx-1)\sigma y-q^a(x-1)y=0, \qquad a\in{\mathbb R}, \quad 0<a<1,
$$
обладает решением в виде (сходящегося) обобщённого степенного ряда
$$
\varphi=\sum_{m_1+m_2\geqslant0}x^{a+m_1a+m_2}=\frac{x^a}{(1-x)(1-x^a)}, \qquad |x|<1.
$$
При этом, как несложно заметить, $L(z)=-z+q^a$ и диофантово условие (\ref{smalldiv}) примет вид
$$
|(m_1a+m_2)\ln q-2\pi m{\rm i}|>c(m_1+m_2)^{-\gamma}.
$$
Воспользовавшись известными рассуждениями о лиувиллевом числе
$$
l=\sum\limits_{k=1}^{\infty}10^{-k!},
$$
покажем, что при $a=2l$ и $q=-1$ ($\ln q=\pi\rm i$) данное условие нарушается для бесконечного числа троек $(m_1,m_2,m)$, какими бы ни были постоянные $c,\gamma>0$. Для этого выберем $m_1=10^{N!}$, $m_2=0$, $m=\sum\limits_{k=1}^{N}10^{N!-k!}$. Тогда
$$
|(m_1a+m_2)\ln q-2\pi m{\rm i}|=2\pi|m_1l-m|=2\pi\sum\limits_{k=N+1}^\infty 10^{N!-k!}=2\pi\sum\limits_{j=1}^\infty m_1^{1-\prod\limits_{i=1}^j (N+i)}<
$$
$$
<2\pi\sum\limits_{j=1}^\infty m_1^{1-(N+j)}<2\pi m_1^{1-N}<c\,m_1^{-\gamma}
$$
при всех $N\in\mathbb N$, начиная с некоторого.

\section{Условие сходимости обобщённого степенного ряда, удовлетворяющего уравнению Малера}

В этом параграфе мы доказываем теорему о достаточном условии сходимости обобщённого степенного ряда, удовлетворяющего уравнению Малера. Выражаем благодарность Стефану Малеку за то, что он привлёк наше внимание к уравнениям Малера в контексте исследования их формальных решений в виде обобщённых степенных рядов.
\medskip

{\bf Теорема 7.} {\it Пусть обобщённый степенной ряд \eqref{powgen} удовлетворяет уравнению \eqref{Mahler} и
$$
F'_{y_0}(x,\varphi,\mu\varphi,\ldots,\mu^n\varphi)=A_0x^{\nu}+o(x^\nu)\in{\cal G}, \qquad A_0\ne0.
$$
Тогда ряд $\varphi$ равномерно сходится во всяком секторе $S$ с вершиной в нуле достаточно малого радиуса и раствора, меньшего $2\pi$.}
\medskip

{Д о к а з а т е л ь с т в о.} По теореме 4 ряд $\varphi$ является элементом подалгебры ${\mathbb C}[[x^\Gamma]]\subset\cal G$, где $\Gamma$ --- конечно порождённая аддитивная полугруппа, образующие которой (обозначаемые, как и ранее, $\rho_1,\ldots,\rho_\tau$) имеют положительные вещественные части и линейно независимы над $\mathbb Z$. При этом, согласно лемме 3 (см. также доказательство теоремы 4), $\varphi$ может быть представлен в виде $\displaystyle\varphi=\sum_{k=1}^Nc_kx^{\lambda_k}+x^{\lambda_N}\psi$, и ряд $\psi\in{\mathbb C}[[x^\Gamma]]$ формально удовлетворяет соотношению
\begin{equation}\label{mahlerforpsi}
\psi=\sum_{{\bf p}\in{\mathbb Z}_+^{n+1}}a_{\bf p}(x)\psi^{p_0}(\mu\psi)^{p_1}\ldots(\mu^n\psi)^{p_n}, \qquad
a_{\bf p}\in{\mathbb C}\{x^\Gamma\}.
\end{equation}
Вновь утверждение теоремы о сходимости достаточно доказать для ряда $\psi$, имеющего вид
$$
\psi=\sum_{{\bf m}\in{\mathbb Z}_+^\tau\setminus\{0\}}c_{\bf m}\,x^{m_1\rho_1+\ldots+m_\tau\rho_\tau}, \qquad c_{\bf m}\in\mathbb C.
$$

Докажем, как и в теоремах 5, 6, что степенной ряд (Тейлора) $\tilde\psi=\iota(\psi)$, представляющий обобщённый степенной ряд $\psi$ в алгебре ${\mathbb C}[[x_1,\ldots,x_\tau]]_*$, имеет непустой полидиск сходимости, повторив основные шаги из доказательств в дифференциальном и $q$-разностном случаях.
\smallskip

{\sl 1. Соотношение для ряда $\tilde\psi$.} Оператор $\mu$, действующий на ${\mathbb C}[[x^{\Gamma}]]$ по правилу
$$
\mu: \sum_{|{\bf m}|>0}c_{\bf m}\,x^{m_1\rho_1+\ldots+m_\tau\rho_\tau}\mapsto\sum_{|{\bf m}|>0}c_{\bf m}\,x^{\ell m_1\rho_1+\ldots+\ell m_\tau\rho_\tau},
$$
с помощью изоморфизма $\iota: {\mathbb C}[[x^{\Gamma}]]\rightarrow{\mathbb C}[[x_1,\ldots,x_\tau]]_*$ естественным образом переносится на
${\mathbb C}[[x_1,\ldots,x_\tau]]_*$,
$$
\tilde\mu: \sum_{|{\bf m}|>0}c_{\bf m}\,x_1^{m_1}\ldots x_\tau^{m_\tau}\mapsto\sum_{|{\bf m}|>0}c_{\bf m}\,x_1^{\ell m_1}\ldots x_\tau^{\ell m_\tau},
$$
так что $\iota\circ\mu=\tilde\mu\circ\iota$. Поэтому, с учётом леммы 6, применив отображение $\iota$ к обеим частям соотношения (\ref{mahlerforpsi}), получим следующее равенство для $\tilde\psi$:
$$
\tilde\psi=\sum_{{\bf p}\in{\mathbb Z}_+^{n+1}}\iota(a_{\bf p})\tilde\psi^{p_0}(\tilde\mu\tilde\psi)^{p_1} \ldots(\tilde\mu^n\tilde\psi)^{p_n}, \qquad
\iota(a_{\bf p})\in{\mathbb C}\{x_1,\ldots,x_\tau\}_*\, ,
$$
которое запишем в виде
\begin{equation}\label{mahlerfortpsi}
\tilde\psi=\sum_{{\bf k}\in{\mathbb Z}_+^\tau\setminus\{0\},\,{\bf p}\in{\mathbb Z}_+^{n+1}}A_{\bf k, p}\,x_1^{k_1}\ldots x_\tau^{k_\tau} \tilde\psi^{p_0}(\tilde\mu\tilde\psi)^{p_1}\ldots(\tilde\mu^n\tilde\psi)^{p_n}.
\end{equation}

{\sl 2. Мажорантное уравнение для ряда $\tilde\psi$.} Для доказательства сходимости ряда $\tilde\psi$ в окрестности $0\in{\mathbb C}^\tau$ рассмотрим уравнение
\begin{equation}\label{majeq3}
W=\sum_{{\bf k}\in{\mathbb Z}_+^\tau\setminus\{0\},\,{\bf p}\in{\mathbb Z}_+^{n+1}}\alpha\,|A_{\bf k, p}|\,x_1^{k_1}\ldots x_\tau^{k_\tau} W^{p_0}W^{p_1}\ldots W^{p_n},
\end{equation}
имеющее вид $W-\widetilde M(x_1,\ldots,x_\tau,W)=0$, где $\widetilde M$ --- голоморфная функция в окрестности $0\in{\mathbb C}^{\tau+1}$,
$\widetilde M(0,\ldots,0,W)\equiv0$. Выбор числа $\alpha\geqslant1$ будет сделан далее. По теореме о неявной функции данное уравнение имеет единственное голоморфное решение $W=W(x_1\ldots,x_\tau)$, обращающееся в нуль в начале координат,
\begin{equation}\label{majseries3}
W=\sum_{|{\bf m}|>0}C_{\bf m}\,x_1^{m_1}\ldots x_\tau^{m_\tau}.
\end{equation}
Наша дальнейшая задача --- доказать, что этот (сходящийся) степенной ряд является мажорантным для $\tilde\psi$, то есть
$$
C_{\bf m}\geqslant0, \qquad |c_{\bf m}|\leqslant C_{\bf m} \quad \forall {\bf m}\in{\mathbb Z}_+^\tau\setminus\{0\},
$$
откуда и будет следовать сходимость $\tilde\psi$ в окрестности $0\in{\mathbb C}^\tau$.
\smallskip

{\sl 3. Рекуррентные соотношения для коэффициентов $c_{\bf m}$.} Используем соотношение (\ref{mahlerfortpsi}), чтобы выразить каждый коэффициент $c_{\bf m}$ ряда $\tilde\psi$ через его предыдущие коэффициенты. Предварительно представим ряд $\displaystyle\tilde\mu^j\tilde\psi=
\sum_{|{\bf m}|>0}c_{\bf m}\,x_1^{\ell^j m_1}\ldots x_\tau^{\ell^j m_\tau}$, $j=0,1,\ldots,n$, в виде
$$
\tilde\mu^j\tilde\psi=\sum_{|{\bf m}|>0}c_{\bf m}^{(j)}x_1^{m_1}\ldots x_\tau^{m_\tau}, \qquad
c_{\bf m}^{(j)}:=\left\{\begin{array}{rcl} c_{{\bf m}/\ell^j}, & & \hbox{если } \ell^j \hbox{ делит } {\bf m}; \\
                                                            0, & & \hbox{в противном случае}.
                                                                    \end{array}\right.
$$
Теперь найдём коэффициент при мономе $x^{m_1}\ldots x^{m_\tau}$ ряда в правой части соотношения (\ref{mahlerfortpsi}). Для этого выразим коэффициент при мономе $x_1^{m_1-k_1}\ldots x_\tau^{m_\tau-k_\tau}$ в произведении $\tilde\psi^{p_0}(\tilde\mu\tilde\psi)^{p_1}\ldots(\tilde\mu^n\tilde\psi)^{p_n}$. Согласно правилу произведения рядов такой коэффициент равен
$$
\sum_{{\bf l}^{(0)}+{\bf l}^{(1)}+\ldots+{\bf l}^{(n)}={\bf m}-{\bf k}}d_{{\bf l}^{(0)}}d_{{\bf l}^{(1)}}\ldots d_{{\bf l}^{(n)}}
$$
(при $p_j=0$ соответствующий мультииндекс ${\bf l}^{(j)}$ отсутствует под знаком суммы, а множитель $d_{{\bf l}^{(j)}}$ --- в произведении), где каждый множитель $d_{{\bf l}^{(j)}}$ --- коэффициент при мономе $x_1^{l_1^{(j)}}\ldots x_\tau^{l_\tau^{(j)}}$ ряда $(\tilde\mu^j\tilde\psi)^{p_j}$ --- выражается формулой
\begin{equation}\label{d3}
d_{{\bf l}^{(j)}}=\sum_{\boldsymbol{\lambda}^{(1)}+\ldots+\boldsymbol{\lambda}^{(p_j)}={\bf l}^{(j)}}
c_{\boldsymbol{\lambda}^{(1)}}^{(j)}\ldots c_{\boldsymbol{\lambda}^{(p_j)}}^{(j)}, \quad j=0,1,\ldots,n.
\end{equation}
Тем самым, искомое рекуррентное соотношение для $c_{\bf m}$ имеет вид
\begin{equation}\label{crec3}
c_{\bf m}=\sum_{0<{\bf k}<{\bf m},\,|{\bf p}|>0}A_{{\bf k},\,{\bf p}}\sum_{{\bf l}^{(0)}+\ldots+{\bf l}^{(n)}=
{\bf m}-{\bf k}}d_{{\bf l}^{(0)}}d_{{\bf l}^{(1)}}\ldots d_{{\bf l}^{(n)}}+A_{{\bf m},\,{\bf 0}},
\end{equation}
и множители $d_{{\bf l}^{(j)}}$ определяются формулой (\ref{d3}), то есть каждый коэффициент $c_{\bf m}$ выражается через предыдущие.
\smallskip

{\sl 4. Рекуррентные соотношения для коэффициентов $C_{\bf m}$.} Используя уравнение (\ref{majeq3}), которому удовлетворяет ряд (\ref{majseries3}), аналогично получаем рекуррентные соотношения для его коэффициентов $C_{\bf m}$:
\begin{equation}\label{Crec3}
C_{\bf m}=\sum_{0<{\bf k}<{\bf m},\,|{\bf p}|>0}\alpha\,|A_{{\bf k},\,{\bf p}}|\sum_{{\bf l}^{(0)}+\ldots+{\bf l}^{(n)}=
{\bf m}-{\bf k}}D_{{\bf l}^{(0)}}D_{{\bf l}^{(1)}}\ldots D_{{\bf l}^{(n)}}+\alpha\,|A_{{\bf m},\,{\bf 0}}|,
\end{equation}
где каждый множитель $D_{{\bf l}^{(j)}}$ --- коэффициент при мономе $x_1^{l_1^{(j)}}\ldots x_\tau^{l_\tau^{(j)}}$ ряда $W^{p_j}$, если $p_j\geqslant1$, --- выражается формулой
$$
D_{{\bf l}^{(j)}}=\sum_{\boldsymbol{\lambda}^{(1)}+\ldots+\boldsymbol{\lambda}^{(p_j)}={\bf l}^{(j)}}C_{\boldsymbol{\lambda}^{(1)}}\ldots C_{\boldsymbol{\lambda}^{(p_j)}}, \qquad j=0,1,\ldots,n
$$
(при $p_j=0$ соответствующий мультииндекс ${\bf l}^{(j)}$ отсутствует под знаком внутренней суммы в (\ref{Crec3}), а множитель $D_{{\bf l}^{(j)}}$ --- в произведении).

Следовательно, при $|{\bf m}|=1$ имеем $C_{\bf m}=|A_{{\bf m},\,{\bf 0}}|$ и, тем самым, $C_{\bf m}$ --- неотрицательное вещественное число. Тогда и все остальные $C_{\bf m}$, при любом
${\bf m}\in{\mathbb Z}_+^\tau\setminus\{0\}$, в силу рекуррентного соотношения (\ref{Crec3}) являются неотрицательными вещественными числами. Более того, выделив в выражении (\ref{Crec3}) для $C_{\bf m}$ слагаемые с мультииндексами ${\bf k}={\bf e}_i$ ($i$-й стандартный единичный вектор в ${\mathbb C}^\tau$) и ${\bf p}=\bar{\bf e}_j$ ($j$-й стандартный единичный вектор в ${\mathbb C}^{n+1}$), получим
$$
C_{\bf m}=\sum_{i=1}^\tau C_{{\bf m}-{\bf e}_i}\sum_{j=1}^{n+1}\alpha\,|A_{{\bf e}_i,\,\bar{\bf e}_j}|+\alpha\,|A_{{\bf m},\,{\bf 0}}|+\ldots,
$$
откуда, при соответствующем выборе $\alpha\geqslant1$, заключаем, что $C_{\bf m}\geqslant C_{{\bf m}-{\bf e}_i}$, $i=1,\ldots,\tau$, и, следовательно,
\begin{equation}\label{est2m}
C_{\bf m}\geqslant C_{{\bf m}'} \qquad \forall {\bf m}>{\bf m}'.
\end{equation}
\smallskip

{\sl 5. Сходимость ряда $\tilde\psi$ $($оценка $|c_{\bf m}|\leqslant C_{\bf m})$.} Воспользуемся индукцией по $|{\bf m}|$, чтобы доказать нужную нам оценку $|c_{\bf m}|\leqslant C_{\bf m}$.
При $|{\bf m}|=1$ согласно (\ref{crec3}) и (\ref{Crec3}) имеем
$$
|c_{\bf m}|=|A_{{\bf m},\,{\bf 0}}|\leqslant\alpha\,|A_{{\bf m},\,{\bf 0}}|=C_{\bf m}.
$$

Предположив, что оценка $|c_{\bf m}|\leqslant C_{\bf m}$ верна при всех $|{\bf m}|<s$, докажем её для $|{\bf m}|=s$ ($s=2,3,\ldots$).
Воспользовавшись соотношением (\ref{crec3}), получим
$$
|c_{\bf m}|\leqslant\sum_{0<{\bf k}<{\bf m},\,|{\bf p}|>0}|A_{{\bf k},\,{\bf p}}|\sum_{{\bf l}^{(0)}+\ldots+{\bf l}^{(n)}=
{\bf m}-{\bf k}}|d_{{\bf l}^{(0)}}d_{{\bf l}^{(1)}}\ldots d_{{\bf l}^{(n)}}|+|A_{{\bf m},\,{\bf 0}}|,
$$
при этом, в силу (\ref{d3}), индуктивного предположения и (\ref{est2m}),
$$
|d_{{\bf l}^{(j)}}|\leqslant\sum_{\boldsymbol{\lambda}^{(1)}+\ldots+\boldsymbol{\lambda}^{(p_j)}={\bf l}^{(j)}}
|c_{\boldsymbol{\lambda}^{(1)}}^{(j)}|\ldots|c_{\boldsymbol{\lambda}^{(p_j)}}^{(j)}|\leqslant
\sum_{\boldsymbol{\lambda}^{(1)}+\ldots+\boldsymbol{\lambda}^{(p_j)}={\bf l}^{(j)}}C_{\boldsymbol{\lambda}^{(1)}}\ldots C_{\boldsymbol{\lambda}^{(p_j)}}=D_{{\bf l}^{(j)}}, \quad j=0,1,\ldots,n.
$$
Следовательно,
$$
|c_{\bf m}|\leqslant\sum_{0<{\bf k}<{\bf m},\,|{\bf p}|>0}|A_{{\bf k},\,{\bf p}}|\sum_{{\bf l}^{(0)}+\ldots+{\bf l}^{(n)}={\bf m}-{\bf k}}
D_{{\bf l}^{(0)}}D_{{\bf l}^{(1)}}\ldots D_{{\bf l}^{(n)}}+|A_{{\bf m},\,{\bf 0}}|\leqslant C_{\bf m}.
$$

Таким образом, коэффициенты $c_{\bf m}$ степенного ряда (Тейлора) $\tilde\psi$ ограничены по абсолютной величине коэффициентами $C_{\bf m}$ степенного ряда (\ref{majseries3}), сходящегося в окрестности $0\in{\mathbb C}^\tau$, следовательно, ряд $\tilde\psi$ также имеет непустой полидиск сходимости.
\smallskip

Завершение доказательства утверждения о сходимости ряда $\displaystyle\psi=\sum_{|{\bf m}|>0}c_{\bf m} x^{m_1\rho_1+\ldots +m_\tau\rho_\tau}$ вновь повторяет соответствующий фрагмент в конце доказательства теоремы 5. {\hfill $\Box$}
\medskip

{\bf Следствие 2.} {\it Обобщённый степенной ряд \eqref{powgen} с {\rm вещественными} показателями степени $\lambda_k\in\mathbb R$, удовлетворяющий уравнению \eqref{Mahler},
равномерно сходится во всяком секторе $S$ с вершиной в нуле достаточно малого радиуса и раствора, меньшего $2\pi$.}
\medskip

{Д о к а з а т е л ь с т в о.} Следует из теоремы 7 аналогично тому, как следствие 1 из теоремы 4. {\hfill $\Box$}
\medskip

В качестве примера применения теоремы 7 (точнее, следствия 2) можно рассмотреть уравнение Малера из \cite[Example 1]{GG2} и его формальное решение в виде обобщённого степенного ряда \cite[Remark 3]{GG2}, которое будет иметь ненулевой радиус сходимости, что было высказано в \cite{GG2} в предположительной форме.

Обобщённый степенной ряд, удовлетворяющий уравнению Малера, может быть сходящимся и при невыполнении условия Теоремы 7, что показано в следующем примере. Подобные примеры можно приводить и в дифференциальном и $q$-разностном случаях, мы ограничимся здесь уравнением Малера.
\medskip

{\bf Пример 6.} Рассмотрим нелинейное уравнение Малера первого порядка
$$
F(x,y,\mu y)=(8y^2x^2+(8x^2+4)y+2x^2+1)\mu y+(4x^2-2)y^2+4yx^2+x^2=0, \qquad \ell=2,
$$
которое имеет решение $y=\varphi$ в виде обобщённого степенного ряда
$$
\varphi=\frac12\;\sum\limits_{k=1}^{\infty}(x^{1+{\rm i}}+x^{1-{\rm i}})^k=\frac12x^{1+{\rm i}}+\frac12x^{1-{\rm i}}+\frac12x^{2+2{\rm i}}+\frac12x^{2-2{\rm i}}+x^2+\ldots,
$$
сходящегося при малых $|x|$ к значению функции
$$
\frac12\;\frac{x^{1+{\rm i}}+x^{1-{\rm i}}}{1-(x^{1+{\rm i}}+x^{1-{\rm i}})}=\frac{x\cos\ln x}{1-2x\cos\ln x}.
$$
При этом частная производная
$$
F'_{y_0}(x,y_0,y_1)=4(4x^2y_0+2x^2+1)y_1+2(4x^2-2)y_0+4x^2
$$
на этом решении не имеет требуемого в теореме 7 вида:
$$
F'_{y_0}(x,\varphi,\mu\varphi)=
-2\,(x^{1+{\rm i}}+x^{1-{\rm i}})+o(x).
$$

\section{Замечания о функциональных уравнениях более общего вида, чем $q$-разностное или Малера}

Рассмотрим аналитическое функциональное уравнение вида
\begin{equation}\label{geneq}
F(x,y(x),y(f(x)),\ldots,y(f^{[n]}(x)))=0,
\end{equation}
где функция $F(x,y_0,y_1,\ldots,y_n)\not\equiv0$ голоморфна около $0\in{\mathbb C}^{n+2}$. Функция $f$, голоморфно отображающая окрестность $({\mathbb C},0)$ нуля в себя, $f(0)=0$, может быть однолистной,
$$
f(x)=qx+o(x), \qquad q\in{\mathbb C}^*,
$$
и тогда представлять росток комплексного диффеоморфизма $({\mathbb C},0)$, либо $\ell$-листной, $\ell\geqslant2$:
$$
f(x)=ax^\ell+o(x^\ell), \qquad a\in{\mathbb C}^*.
$$
Уравнения $q$-разностное (\ref{qdiff}) и Малера (\ref{Mahler}) являются частными случаями уравнения (\ref{geneq}) при $f(x)=qx$, $q\ne1$, и $f(x)=x^\ell$, соответственно.

Сделаем некоторые комментарии относительно существования и единственности локального голоморфного решения уравнения (\ref{geneq}) и сходимости обобщённых степенных рядов, формально ему удовлетворяющих. Прежде всего заметим, что для большинства функций $f$ данное уравнение может быть сведено к уже рассмотренным нами с точки зрения этих вопросов $q$-разностному уравнению или уравнению Малера.

Действительно, если $f$ однолистна и производная $f'(0)=q$ такова, что $|q|\ne1$ или $|q|=1$ и удовлетворяет условию Зигеля (\ref{siegel}), то найдётся
{\it биголоморфная} в нуле функция $\xi$, $\xi(0)=0$, такая, что $\xi(qx)=f(\xi(x))$, см. п. 2.2.1. Если же $f$ $\ell$-листна, то аналогичная функция $\xi$ удовлетворяет соотношению $\xi(x^\ell)=f(\xi(x))$, см. п. 2.3.1. Таким образом, каждое решение $y(x)$ уравнения (\ref{geneq}) преобразуется в решение
$\tilde y(x)=y(\xi(x))$ $q$-разностного уравнения
\begin{equation}\label{qdiff1}
F(\xi(x),\tilde y(x),\tilde y(qx),\ldots,\tilde y(q^nx))=0
\end{equation}
в первом случае и уравнения Малера
\begin{equation}\label{Mahler1}
F(\xi(x),\tilde y(x),\tilde y(x^\ell),\ldots,\tilde y(x^{\ell^n}))=0
\end{equation}
--- во втором случае. Обратно, соотношение $y(x)=\tilde y(\xi^{-1}(x))$ обеспечивает переход от соответствующего из уравнений (\ref{qdiff1}), (\ref{Mahler1}) к уравнению (\ref{geneq}).

\subsection{Существование и единственность локального голоморфного решения}

Для функции $\widetilde F(x,y_0,y_1,\ldots,y_n):=F(\xi(x),y_0,y_1,\ldots,y_n)$, очевидно, $\widetilde F'_{y_j}(0)=F'_{y_j}(0)$, поэтому, рассмотрев
многочлен (\ref{polyP}),
$$
P(z)=\sum_{j=0}^nA_jz^j, \qquad A_j:=F'_{y_j}(0),
$$
из результатов пп. 2.2, 2.3 и сказанного выше получим следующее утверждение о существовании и единственности локального голоморфного решения уравнения (\ref{geneq}): \medskip

{\it Если $f$ однолистна и выполнено одно из следующих трёх условий относительно её производной $f'(0)=q$:
\begin{itemize}
\item[--] $|q|>1$, $A_n\ne0$ и $P(q^k)\ne0$ при всех $k\in\mathbb N$,
\item[--] $|q|<1$, $A_0\ne0$ и $P(q^k)\ne0$ при всех $k\in\mathbb N$,
\item[--] $|q|=1$, $P\not\equiv0$ и найдутся положительные постоянные $c,\nu$ такие, что для каждого корня $z=a$ многочлена $(z-1)P(z)$ число $q$ удовлетворяет условию
$$
|q^k-a|\geqslant c\,k^{-\nu}, \qquad k\in{\mathbb N},
$$
\end{itemize}
или если $f$ $\ell$-листна и $A_0\ne0$, то уравнение $(\ref{geneq})$ имеет и при том единственное решение $y=y(x)$, голоморфное в точке $x=0$ и такое, что $y(0)=0$.}

\subsection{Формальное решение в виде обобщённого степенного ряда и его сходимость}

Подстановка функции $f$ в обобщённый степенной ряд $\displaystyle\varphi=\sum_{k\geqslant1}c_kx^{\lambda_k}\in\cal G$ приводит к формальному выражению вида
$$
\varphi(f(x))=\sum_{k=1}^{+\infty}c_kq^{\lambda_k}x^{\lambda_k}(1+o(x))^{\lambda_k}=\varphi(qx)+\sum_{j=1}^{+\infty}\varphi_j(x), \qquad \varphi_j\in{\cal G},
$$
в случае когда $f$ однолистна, $f'(0)=q\ne0$, или, соответственно, к выражению вида
$$
\varphi(f(x))=\sum_{k=1}^{+\infty}c_ka^{\lambda_k}x^{\ell\lambda_k}(1+o(x))^{\lambda_k}=\varphi(ax^\ell)+\sum_{j=1}^{+\infty}\varphi_j(x), \qquad
\varphi_j\in{\cal G},
$$
в случае когда $f$ $\ell$-листна, $f(x)=ax^\ell+o(x^\ell)$, $a\ne0$. При этом множество показателей степени каждого обобщённого степенного ряда $\varphi_j$ в первом случае представляет собой последовательность $\{\lambda_k+j\}_{k\geqslant1}$, а во втором случае --- последовательность $\{\ell\lambda_k+j\}_{k\geqslant1}$. Таким образом, для каждого $N>0$ в ряду $\varphi(f(x))$ имеется только конечное число показателей степени $\lambda$ таких, что ${\rm Re}\,\lambda<N$ и, тем самым, данный ряд является обобщённым степенным рядом. Следовательно, корректно определено понятие формального решения $\varphi\in\cal G$ уравнения (\ref{geneq}).

Пусть обобщённый степенной ряд $\varphi\in\cal G$ формально удовлетворяет уравнению (\ref{geneq}). Рассмотрим биголоморфную в нуле функцию $\xi(x)=cx+o(x)$, $c\ne0$, удовлетворяющую соотношению $\xi(qx)=f(\xi(x))$ или $\xi(x^\ell)=f(\xi(x))$, в зависимости от того, однолистна или $\ell$-листна функция $f$. Тогда обобщённый степенной ряд $\tilde\varphi(x)=\varphi(\xi(x))\in\cal G$ будет удовлетворять, соответственно, $q$-разностному уравнению (\ref{qdiff1}) или уравнению Малера (\ref{Mahler1}). Обозначив, как и в предыдущем пункте, через $\widetilde F(x,y_0,y_1,\ldots,y_n)=F(\xi(x),y_0,y_1,\ldots,y_n)$ функцию, определяющую эти уравнения,
очевидно, имеем
$$
\widetilde F'_{y_j}(x,\widetilde\Phi)=F'_{y_j}(\xi(x),\widetilde\Phi), \quad j=0,1,\ldots,n, \qquad \mbox{где} \quad \widetilde\Phi=(\tilde\varphi,\Delta\tilde\varphi,\ldots, \Delta^n\tilde\varphi),
$$
$\Delta=\sigma$ или $\mu$. Поэтому, допустив, что
$$
F'_{y_j}(x,\Phi(x))=A_jx^\nu+o(x^\nu), \quad j=0,1,\ldots,n, \qquad \mbox{где} \quad \Phi(x)=(\varphi(x),\varphi(f(x)),\ldots,\varphi(f^{[n]}(x))),
$$
получим:
$$
\widetilde F'_{y_j}(x,\widetilde\Phi(x))=F'_{y_j}(\xi(x),\Phi(\xi(x)))=A_jc^\nu x^\nu+o(x^\nu), \qquad j=0,1,\ldots,n.
$$
При этом если множество $\{\lambda_k\}_{k\geqslant1}$ показателей степени ряда $\tilde\varphi(x)$ лежит в полугруппе $\Gamma$, $1\in\Gamma$, то и множество $\{\lambda_k+j\}_{k\geqslant1,j\geqslant0}$ показателей степени ряда $\varphi(x)=\tilde\varphi(\xi^{-1}(x))$ обладает этим же свойством. Таким образом, с учётом результатов пп. 3.3, 3.4 и теорем 6, 7 приходим к следующему утверждению о формальном решении $\varphi\in\cal G$ аналитического функционального уравнения (\ref{geneq}) и его сходимости:
\medskip

{\it Пусть обобщённый степенной ряд $\varphi\in\cal G$ формально удовлетворяет уравнению $(\ref{geneq})$;
\begin{itemize}
\item[--] если функция $f$ однолистна $($и для её производной $f'(0)=q$ выполнено условие Зигеля $(\ref{siegel})$ в случае $|q|=1)$ и
$$
F'_{y_j}(x,\Phi(x))=A_jx^\nu+o(x^\nu)\in{\cal G}, \qquad j=0,1,\ldots,n,
$$
где число $\nu$ --- одинаковое для всех $j=0,1,\ldots,n$ и при этом хотя бы одно $A_j$ отлично от нуля, или
\item[--] если функция $f$ $\ell$-листна и
$$
F'_{y_0}(x,\Phi(x))=A_0x^\nu+o(x^\nu)\in{\cal G}, \qquad A_0\ne0,
$$
\end{itemize}
то $\varphi\in{\mathbb C}[[x^\Gamma]]\subset\cal G$, где $\Gamma$ --- конечно порождённая аддитивная полугруппа, образующие которой имеют положительные вещественные части и линейно независимы над $\mathbb Z$. При этом во втором случае ряд $\varphi$ равномерно сходится во всяком секторе $S$ с вершиной в нуле достаточно малого радиуса и раствора, меньшего $2\pi$, а в первом случае сходимость имеет место, если для образующих $\rho_1,\ldots,\rho_\tau$ полугруппы $\Gamma$ выполнено диофантово условие из теоремы {\rm 6} $($или условия теоремы {\rm 6 bis)}.}